\documentclass[preprint,12pt]{elsarticle}

\usepackage{amssymb}
\usepackage{amsmath}
\usepackage{amsthm}

\usepackage{amsrefs}
\usepackage{subfig}
\usepackage{color}
\usepackage{hyperref}
\usepackage{float}
\usepackage{graphicx}
\usepackage{subfig}
\usepackage{booktabs}

\usepackage{listings}
\usepackage{xcolor}
\definecolor{codegreen}{rgb}{0,0.6,0}
\definecolor{codegray}{rgb}{0.5,0.5,0.5}
\definecolor{codepurple}{rgb}{0.58,0,0.82}
\definecolor{backcolour}{rgb}{0.95,0.95,0.92}
\lstdefinestyle{mystyle}{
    backgroundcolor=\color{backcolour},   
    commentstyle=\color{codegreen},
    keywordstyle=\color{magenta},
    numberstyle=\tiny\color{codegray},
    stringstyle=\color{codepurple},
    basicstyle=\ttfamily\footnotesize,
    breakatwhitespace=false,         
    breaklines=true,                 
    captionpos=b,                    
    keepspaces=true,                 
    numbers=left,                    
    numbersep=5pt,                  
    showspaces=false,                
    showstringspaces=false,
    showtabs=false,                  
    tabsize=2
}
\newtheorem{thm}{Theorem}[section]

\newtheorem{proposition}{Proposition}[section]
\newtheorem{assumption}{Assumption}[section]
\newtheorem{remark}{Remark}[section]

\DeclareUnicodeCharacter{223C}{\ensuremath{~}}

\journal{Journal of Computational Physics}

\begin{document}

\begin{frontmatter}



\title{{A Hybrid Two-level MCMC Framework to Accelerate Posterior Mean Estimation with Deep Learning Surrogates for Bayesian Inverse Problems}}


\author[nvidia]{Juntao Yang} 
\author[nvidia]{Jeff Adie} 
\author[nvidia]{Simon See} 
\author[cambridge,INGV]{Adriano Gualandi}
\author[nus]{Gianmarco Mengaldo} 

\affiliation[nvidia]{organization={NVIDIA AI Technology Center}, 
            addressline={07-03 Suntec Tower Three, 8 Temasek Blvd}, 
            country={Singapore}}
\affiliation[cambridge]{organization={University of Cambridge}, 
            addressline={Downing St., Cambridge CB2 3EQ}, 
            country={United Kingdom}}
\affiliation[INGV]{organization={Istituto Nazionale di Geofisica e Vulcanologia}, 
            addressline={Via di Vigna Murata, 605, 00143 Roma RM},
            country={Italy}}
\affiliation[nus]{organization={National University of Singapore},
            addressline={21 Lower Kent Ridge Rd},
            country={Singapore}}

\begin{abstract}
Bayesian inverse problems arise in various scientific and engineering domains, and solving them can be computationally demanding. 
This is especially the case for problems governed by partial differential equations, where the repeated evaluation of the forward operator is extremely expensive. 
Recent advances in Deep Learning (DL)-based surrogate models have shown promising potential to accelerate the solution of such problems. 
However, despite their ability to learn from complex data, DL-based surrogate models generally cannot match the accuracy of high-fidelity numerical models, which limits their practical applicability. 
We propose a novel hybrid two-level Markov Chain Monte Carlo (MCMC) method that combines the strengths of DL-based surrogate models and high-fidelity numerical solvers to {compute the posterior mean of Quantities of Interest (QoI) in} Bayesian inverse problems governed by partial differential equations. 
The intuition is to leverage the evaluation speed of a DL-based surrogate model as the base chain, and correct its errors using a limited number of high-fidelity numerical model evaluations in a correction chain; hence its name hybrid two-level MCMC method. 
Through a detailed theoretical analysis, we show that our approach can achieve the same accuracy as a pure numerical MCMC method while requiring only a small fraction of the computational cost.
{The theoretical analysis is further supported by several numerical experiments, namely a Poisson, a non-linear reaction-diffusion, and a Navier-Stokes equation.} 
The proposed hybrid framework can be generalized to other approaches such as the ensemble Kalman filter and sequential Monte Carlo methods.
\end{abstract}

\begin{keyword}
Markov Chain Monte Carlo \sep Deep Learning \sep Bayesian Inverse Problems
\end{keyword}
\end{frontmatter}


\section{Introduction}
Inverse problems arise in various fields of applied science, including design optimization in engineering, seismic inversion in geophysics, and data assimilation in weather forecasting~\cite{MR2102218}. 
The behavior of these systems is described by a mathematical model that frequently consists of a system of partial differential equations that depends on a set of inputs and parameters.
Inverse problems involve determining the inputs or parameters of the mathematical model based on observations or partial observations of the model solution. 
The mathematical model, in the context of inverse problems, is also known as the forward problem, and it is typically expressed as
\begin{equation}\label{eq:generic-problem}
    y = \mathcal{G}(z),
\end{equation}
where $\mathcal{G}$ is the forward operator (also referred to as the forward map), $z$ represents the inputs or parameters, and $y$ are the observed data defined in Equation~\eqref{eq:generic-problem}.

The objective of an inverse problem is to identify the inputs or parameters $z$, or some {Quantities of Interest (QoI) that depends on $z$, denoted as $Q(z)$. This can be, for instance, the permeability field of a Darcy flow's subsurface model, or the initial condition of a Navier Stokes equation.}
Optimization techniques, such as least squares optimization, are commonly employed to solve inverse problems~\cites{MR2130010, MR3285819}. 
However, inverse problems are often ill-posed, meaning they may lack uniqueness, stability, or the existence of a solution.

To address the challenges related to ill-posedness, Tikhonov regularization is frequently used~\cites{MR3285819, MR2102218}{, where the inverse problem is solved as an optimization problem 
\begin{equation}
\label{eqn:MAP}
\text{argmin}_{z \in U} \left( \frac{1}{2}\|y - \mathcal{G}(z)\|_Y^2 + \frac{1}{2}\|z-m_0\|_U^2 \right),    
\end{equation}
with norms $\|\cdot \|_Y$ and $\|\cdot\|_U$ defined on two Banach spaces, namely $Y$ and $U$ representing the data and model space, respectively.}
Although the incorporation of Tikhonov regularization might initially seem arbitrary, it can be explicitly interpreted from a Bayesian perspective as a prior distribution. 
This connection bridges the optimization approach with the probabilistic Bayesian framework, where data are considered as observations subject to noise {together with a prior belief on the parameters or inputs $z$, namely $m_0$ in Equation~\ref{eqn:MAP}}. 
In this context, a noise term $\eta$ is added to the forward operator defined in Equation~\eqref{eq:generic-problem} to account for observational noise,
\begin{equation}\label{eq:bayesian-setup}
y = \mathcal{G}(z) + \eta.
\end{equation}
{ From Equation~\eqref{eq:bayesian-setup}, we can write a posterior distribution for the model parameters given the observations as follows
\begin{equation}
\label{eq:posterior}
\gamma^y = P(z \lvert y) = \frac{P(y \lvert z)P(z)}{P(y)},
\end{equation}
where $P(y)$ is the evidence of the data, $P(z)$ is the prior probability about the model parameters, and $P(y \lvert z)$ is the likelihood of the given observations. Rather than solving for the full posterior distribution, it is often convenient to solve a maximization problem just for the numerator of the right hand side. {Under the assumption that the prior distribution of $z$ as well as the distribution of $\eta$ are Guassian (the latter with zero mean), this is equivalent to maximizing}

\begin{equation}
\label{eq:likelihoodprior}
\gamma^y \propto P(y \lvert z)P(z) \propto \exp\left(-\frac{1}{2}\|y - \mathcal{G}(z)\|_Y^2 - \frac{1}{2}\|z-m_0\|_U^2\right),
\end{equation}
where norms $\|\cdot \|_Y$ and $\|\cdot\|_U$ are defined on the covariance of the prior ($U$) and of the noise ($Y$).
Finding the maximum a posterior (MAP) of the distribution in Equation~\eqref{eq:likelihoodprior}, leads to the same optimization problem as Equation~\eqref{eqn:MAP}.

In this work, we approach inverse problems, such as the one in Equation~\eqref{eq:bayesian-setup}, leveraging Bayes' rule (Equation~\ref{eq:posterior}).}
In this Bayesian context, we can calculate the posterior distribution, which reflects the updated beliefs about the unknowns after observing the data, {even under more general assumptions than Gaussianity of $\eta$ and of the prior distribution of $z$}.
Bayes' theorem is often expressed formally with measure-theoretic terms from the mathematical framework of the Radon-Nikodym theorem, such that probability masses or densities over real numbers can be extended to probability measures over any arbitrary sets~\cite{MR2767184}. 
In this framework, the posterior measure and the prior measure are related through the Radon-Nikodym derivative~\cite{MR2652785}.

Bayesian inverse problems can be finite-~\cite{MR3285819}, or infinite-dimensional~\cite{MR2652785}. 
The former usually arises in the context of parameter estimation, whereby a finite set of parameters is of interest. 
The latter instead arises in the context of full-field inversion of partial differential equations (PDE) problems, where infinite-dimensional functions are of interest. 
In our paper, we adopt the Bayes theorem in measure-theoretic terms, which is compatible with the infinite-dimensional setup.

{Solving Bayesian inverse problems typically leads to the repeated solution of the forward problem in equation~\eqref{eq:bayesian-setup}.} 
For example, to solve PDE-based (i.e. infinite-dimensional) Bayesian inverse problems, it is necessary to discretize the continuous PDE problem via a suitable numerical method, such as the finite element method~\cite{CANTWELL2015205} or the finite volume method~\cite{MR1925043}.
This often leads to a high-dimensional linear system of equations, that are extremely expensive computationally.
These computationally expensive high-dimensional linear systems, in turn, need to be solved several times to approximate the posterior distribution, making the problem intractable due to the curse of dimensionality.

Recent developments in deep learning (DL) have provided a possible pathway to accelerate the solution of Bayesian inverse problems. 
In particular, DL-based models (i.e., deep neural networks) can be used as surrogate models to substitute the computationally expensive high-dimensional linear system of equations that arise when numerically discretizing the continuous PDE problem. 
Two of the critical advantages of DL-based surrogates are their fast differentiability (thanks to automatic differentiation) and fast evaluation. 
The first feature makes DL-based surrogates an excellent candidate for solving Bayesian inverse problems using deterministic methods, such as variational methods; see, e.g.~\cite{maulik2022efficient}.
Variational methods typically lead to finding the maximum a posteriori (MAP) with optimization techniques. 
Gradient-based optimization such as gradient descent and L-BFGS is {a} family of the most used optimization techniques for variational methods. 
Gradients can be easily computed by automatic differentiation from a differentiable DL-based surrogate model, which makes the DL-based surrogate model a great fit.
The second feature makes them an excellent candidate for sampling-based statistical methods, such as Markov Chain Monte Carlo and ensemble Kalman filter. 
These sampling techniques approximate the posterior distribution through Monte Carlo samples, which typically converge at a slow rate of $1/\sqrt{M}$ ($M$ being the number of samples). 
In this case, the fast evaluation speed of DL-based surrogates can be utilized to dramatically accelerate sampling procedures.

In the literature, several works explored the use of DL-based models for inverse problems.  
For example, physics-informed neural networks (PINNs) have demonstrated their ability to solve parameterized PDEs; feature that can be used for finite-dimensional inverse problems, such as design optimization~\cites{MR3881695, tay2025optimization}.
{It has been shown that neural operators can learn linear and non-linear mappings between function spaces~\cite{LU2022114778};} hence, they are a promising category of DL-based surrogates for infinite-dimensional inverse problems (e.g.~\cite{DBLP:conf/iclr/LiKALBSA21}). 
Indeed, the fast evaluation properties of neural operators make them extremely competitive for high-dimensional problems with respect to more traditional surrogate models, such as generalized polynomial chaos and Gaussian processes \cites{marzouk2009stochastic, MR3084684}. 
These recent advances have made the solution of otherwise intractable inverse problems a real possibility.

However, despite the advantages of DL-based surrogate models for inverse problems, there are still some key areas that need to be addressed. 
Namely, a complete mathematical framework to estimate the error bounds of a deep neural network model is still missing.
The expected error bounds (also referred interchangeably as error estimates) consist of three components: approximation error, optimization error, and generalization error~\cite{jin2020quantifying}. 
While universal approximation theorems exist~\cites{LESHNO1993861, chen1995universal}, together with the expressivity analysis of neural networks (that depends on the number of layers and nodes provided)~\cites{petersen2018optimal, aadebffb88b448c89c654fcdda52c02b}, these only address the approximation error. 
A large body of literature attempts to address the optimization error by investigating the landscape of non-convex loss functions as well as the optimization process by stochastic gradient~\cites{jin2020quantifying,pmlr-v97-allen-zhu19a,du2019gradient,10.1093/imanum/drz055}. 
Some works also attempt to quantify the generalization error~\cite{jin2020quantifying}. 
However, a general theory is still lacking as most existing analyses make several simplifying assumptions that do not hold for practical problems.

Because of this lack of theoretical error framework, DL-based models are commonly treated as a black-box, and the expected error bounds are empirically estimated {with a test dataset}, noting that {increasing test accuracy of DL models often requires exponentially more data} (property known as the power law)~\cites{hestness2017deep, 10.1073/pnas.2311878121}. 

In the context of Bayesian inverse problems, the lack of a rigorous theoretical framework on DL models' error bounds hinders the adoption of DL surrogates in critical applications, where a desired error estimated a priori may be required. 
In fact, in Bayesian inverse problems, a naive replacement of a high-fidelity numerical PDE solver with a DL-based surrogate will lead to propagation of the error of the surrogate to the posterior distribution~\cite{CAO2023112104}. 
For instance, in the context of MCMC methods, that sample the posterior distribution with an ergodic Markov Chain generated by a given algorithm (e.g., the Metropolis-Hasting algorithm~\cite{brooks2011handbook}), the { estimation} error of the posterior { mean on the QoI} depends on the surrogate error shown in Section~\ref{sec:dl-based-approx}. 

Therefore, in order to make the DL-based surrogate model practically useful in solving Bayesian inverse problems, e.g., using MCMC methods, the posterior error induced by the DL surrogate needs to be contained to a given a priori threshold, an aspect that is still lacking and that represents a key gap in the literature.

In this work, we focus on MCMC methods that are among the most widely adopted methods for solving Bayesian inverse problems, given their ability to handle high-dimensional problems and their embedded uncertainty estimates.
{ More specifically, we focus on problems that aim to compute statistical properties such as the posterior mean and variance of some QoI.}
However, MCMC methods have a critical issue: they are extremely expensive computationally. 
To address this issue, we propose a new MCMC approach to { estimate the posterior mean in} Bayesian inverse problems that we named \textit{two-level hybrid MCMC approach}. 
The new method leverages the fast evaluation properties of DL surrogates to accelerate the MCMC method, while also using a high-fidelity numerical model for accuracy.
The latter aspect allows for a priori theoretical error estimates { of the posterior mean of the QoI which can be controlled by the choice of the high-fidelity numerical model.}  
This is typically not readily available when only using DL surrogates.

{ Our method draws inspiration from numerical multilevel MCMC methods \cite{doi:10.1137/130915005}. 
Generally speaking, there are two approaches to improve the computational cost of MCMC via multilevel methods. 

The first approach uses coarser resolution numerical models as filters to pre-screen the proposed sample before going to the acceptance/rejection step in the Metropolis-Hasting algorithm with expensive high-resolution numerical models~\cites{MR4537564, Christen01122005, efendiev2005efficient, cui2011bayesian, MR2231730}. 
Hybrid MCMC algorithms were also proposed within this context, using traditional surrogate models (including generalized polynomial chaos and Gaussian processes) as the pre-screening filters~\cites{laloy2013efficient,reddy2024accelerating}.
This first approach improves the sampling efficiency of the Metropolis-Hasting algorithm by improving acceptance rates at the finer resolution level.
However, the number of effective samples required remains unchanged to reach a target error $\epsilon \leq C M_{fine}^{-1/2}$, where $C$ is a constant and $M_{fine}$ is the number of effective samples at fine resolution.
For an elliptic partial differential equation-governed multiscale Bayesian inverse problem, the overall computational complexity of such approach remains at best $\mathcal{O}(\epsilon^{-d-2})$  where $\epsilon$ is the desired approximation error for the posterior mean of QoI and $d$ is the dimension of the problem~\cite{MR3084684}. 

The second approach is based on the idea of telescoping sum, upon which our method is based. 
The telescoping sum technique was initially proposed for the multilevel Monte Carlo~\cite{MR2436856}, then it was extended to the MCMC method with some modifications~\cites{doi:10.1137/130915005,MR3084684,MR4523340}. 
Instead of focusing on sampling efficiency, this technique exploits the fact that the variance of the difference between solutions of two resolution levels $L$ and $L-1$ in a PDE-constrained Bayesian inverse problem typically decreases with larger $L$, where $L$ is the level of mesh refinement with the mesh size $h$ of scale $\mathcal{O}(2^{-L})$. 
This allows an efficient multilevel approach to achieve an accurate posterior mean approximation by the telescoping sum technique, where less expensive high resolution samples are needed when $L$ is large thanks to the smaller variance. 
In certain problems, such as the Bayesian inverse problem with elliptic PDEs with bi-hierarchical setup, the computational complexity can be reduced to $\mathcal{O}(\epsilon^{-d})$~\cite{MR3084684}. 

{Despite the attempts in~\cites{laloy2013efficient,reddy2024accelerating} to build hybrid MCMC methods using surrogate model in delayed-acceptance-like MCMC algorithms, there have been no attempts to build hybrid MCMC methods using surrogate models in the telescopic approach. The method we propose in this paper fills this gap.}

We note that in the telescoping sum approach, each MCMC chain runs independently at different levels. 
As a consequence, there is no restriction on the type of MCMC algorithm that can be used to accelerate the sampling efficiency of each MCMC chain, and several potential methods can be used, including the Delayed Acceptance algorithm, the Preconditioned Crank–Nicolson algorithm, and the Stochastic Newton MCMC method, among others~\cites{doi:10.1137/110845598, cotter2013mcmc}.  
}

\begin{figure}
    \centering
    \includegraphics[width=\linewidth]{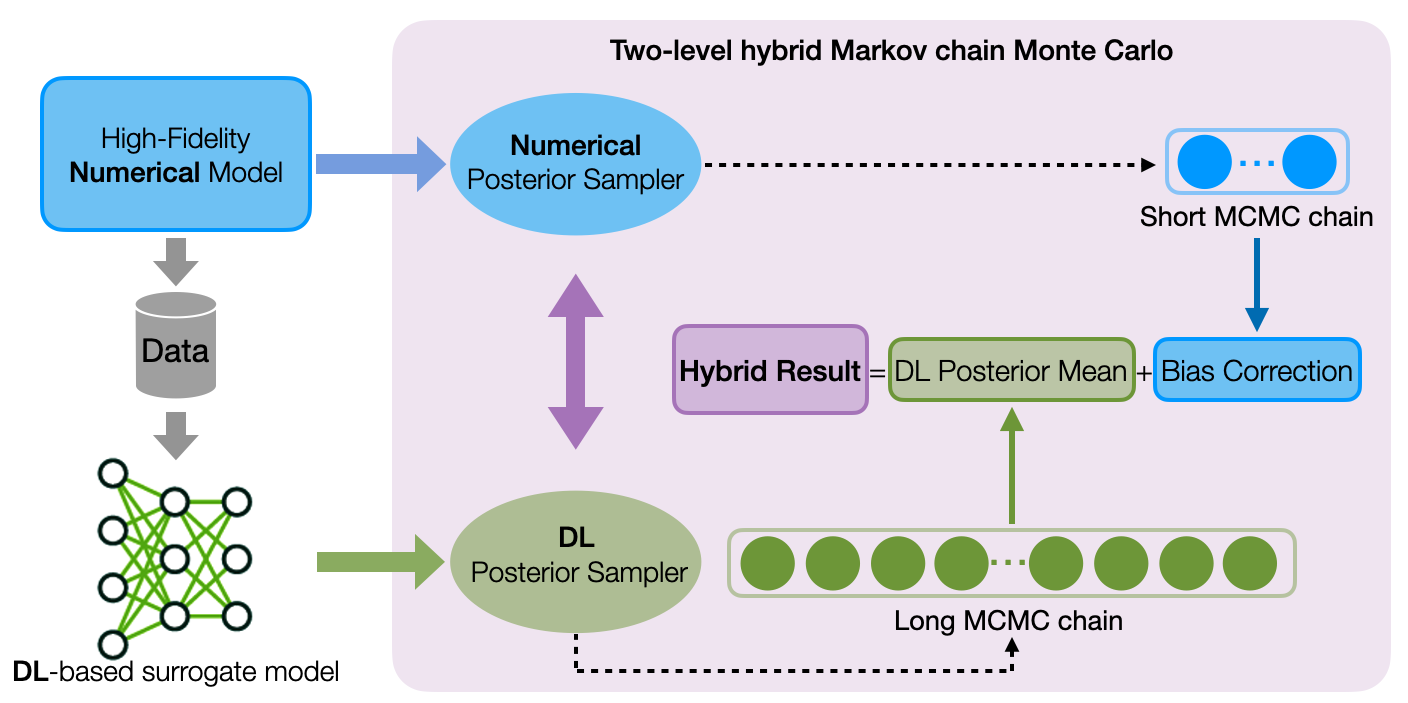}
    \caption{Hybrid two-level MCMC}
    \label{fig:overview}
\end{figure}

{With the telescoping sum technique,} our \textit{two-level hybrid MCMC approach}  (also referred simply to as Hybrid MCMC) uses a DL-based surrogate model to obtain a base MCMC chain with a large number of samples (leveraging the fast evaluation speed of the DL surrogate).
A short correction MCMC chain is then generated to sample the differences between the high-fidelity numerical model and the DL-based surrogate model.
This is done to correct for the bias introduced by the surrogate as shown in Fig~\ref{fig:overview}.
{Despite our focus on the fast developing deep learning-based surrogate models, the hybrid two level MCMC structure is generally applicable to all kinds of surrogate models including Polynomial Chaos, Gaussian Process Regressions and etc.}

We provide a detailed theoretical analysis, showing that the new method has the same a priori error bound $\mathcal{O}(h)$ {as a plain, i.e. single chain, MCMC method} that uses a high-fidelity numerical model, discretized at a known mesh size equal to $h$. 
However, our method requires a small fraction of the computational cost necessary to run a plain MCMC chain with a high-fidelity numerical model.
{ Despite the computational advantages, we shall mention that the proposed approach is limited to the computation of the posterior mean of QoI, as it is not possible to generate a histogram as can be commonly done when using single-chain MCMC algorithms. 
Yet, several useful {statistical quantities}, such as variance, cumulative distribution function, a quantile, and the associated conditional value-at-risk, may still be estimated with techniques such as those mentioned in~\cite{doi:10.1137/17M1135566}.}

We complement the theoretical findings with numerical experiments on {an elliptic, a reaction-diffusion, and a fluid dynamics problem.}
The numerical results support the theoretical findings and highlight how the new method provides a lightweight DL-based surrogate based alternative to existing MCMC approaches, with rigorously defined a priori error bounds. 
The latter aspect closes the gap in the literature regarding the lack of rigorous error bounds when using DL-based surrogates, and constitutes an important milestone for the fast solution of Bayesian inverse problems via deep learning.
The rest of the paper is organized as follows. 
Section~\ref{sec:approach}, introduces the Bayesian inverse problem setup (Section~\ref{sec:setup}), the approximation of the forward problem (Section~\ref{section:approx-forward-problem}), and the new hybrid two-level MCMC for a uniform prior (Section~\ref{sec:two-level-mcmc-uniform}), noting that the Gaussian prior case is presented in \ref{app:two-level-mcmc-gaussian}. 
Section~\ref{sec:experiments} shows the numerical experiments that validate the theoretical estimates provided in Section~\ref{sec:approach}.
Section~\ref{sec:conclusions} draws some closing remarks, including limitations and future work.

\section{A new approach to accelerate Bayesian inversion}
\label{sec:approach}

\subsection{Bayesian inverse problem setup}\label{sec:setup}

To present our new approach, we first introduce the theoretical background of Bayesian inverse problems. 
We consider inverse problems governed by a forward mathematical model as the one defined in equation~\eqref{eq:bayesian-setup}, where the underlying system is constituted of PDEs.
More formally, the PDE-based forward model predicts the states $u$ provided the inputs/parameters $\mathbf{z} = \{z_1, z_2, ..., z_n\}$. 
In order to introduce the problem setup based on Equation~\eqref{eq:bayesian-setup}, we need to define the inputs/parameters $\mathbf{z}$, the forward operator (or forward map) $\mathcal{G}(\mathbf{z})$, and the observational noise $\eta$.

We start by defining the inputs/parameters $\mathbf{z}$.
These represent a finite number of constants or functions within the governing equations, or the coefficients associated with the spectral expansion of a random field defining the initial conditions or forcing terms. 
For example, $\mathbf{z}$ can be the Lamé constants of the material in the elasticity equation of solid mechanics~\cite{MR2322235},  the coefficient of the Karhunen–Loève expansion of the porosity random field in the subsurface flow model~\cite{MR2421969}, or the initial condition and random forcing in the Navier-Stokes equations~\cite{MR2558668}. In many practical applications, the following truncated Karhunen–Loève expansion of a random field $K$ is commonly used,
\begin{equation}
\label{eq:KL-expansion}
    K(z) = \Bar{K} + \sum_{j=1}^n z_j \psi_j,
\end{equation}
where $\Bar{K}, \psi_j$ are functions in {$L^\infty(D)$, where $D$ is the physical domain.}
For simplicity, hereafter we name $\mathbf{z}$ as the \textit{parameters} of the forward problem, without lacking generality.

In the context of Bayesian inverse problems, we need to define the \textit{prior} probability distribution for the parameters $\mathbf{z}$. 
To this end, we consider a uniform prior, where $z_i$ in $\mathbf{z}$ is uniformly distributed within $[-1, 1]$. 
By denoting $\mathcal{B}$ the Borel $\sigma$-algebra, we define a measurable space $(U, \Theta)$, where $\Theta$ is a $\sigma$-algebra on $U$, defined as $\Theta = \bigotimes_{j=1}^n \mathcal{B}([-1, 1])$, and $U = [-1, 1]^n$ is the parameter space. 
Together with the prior measure $\gamma = \bigotimes_{j=1}^n \frac{d z_i}{2}$ on the measurable space $(U, \Theta)$, we have the complete probability space $(U, \Theta, \gamma)$.

With the prior measure of $\mathbf{z}$ defined, we focus on the forward operator (or map) $\mathcal{G}$. 
Within the framework just introduced, the forward operator can be written as follows,
\begin{equation}\label{eq:continuous-forward-map}
\mathcal{G} : U \rightarrow \mathbb{R}^k\;\;\; \forall \mathbf{z} \in U; \;\;\;
\mathcal{G}(\mathbf{z})= (\mathcal{F}_1(u(\mathbf{z})), \mathcal{F}_2(u(\mathbf{z})), ..., \mathcal{F}_k(u(\mathbf{z}))),
\end{equation}
where $u(\mathbf{z}) \in V$ is the state solution of the forward PDEs which depends on the input $\mathbf{z}$, and $\mathcal{F}_i(\cdot), i=1,2,...k$ are $k$ continuous bounded linear functionals.  
$V$ is a suitable vector space, e.g. a Sobolev space {over $D$}, which depends on the specific physical problem. 
$\mathcal{F}$ is included in the forward operator to better reflect real-world problems, where real-world observations are usually discrete and sparse while the state solutions of a PDE system are typically continuous functions. 
For example, in the context of weather data assimilation, $\mathcal{F}$ is known as the observation operator~\cite{skamarock2019description}.
In order to frame our Bayesian problem and guarantee the existence of the posterior, we need to formulate a key assumption on the forward operator $\mathcal{G}$. 
\begin{assumption}
\label{assumption:forward-continuity}
The forward operator $\mathcal{G}(z) : U \rightarrow \mathbb{R}^k$ is a continuous map from the measurable space $(U, \Theta)$ to $(\mathbb{R}^k, \mathcal{B}(\mathbb{R}^k))$. 
\end{assumption}
Assumption \ref{assumption:forward-continuity} is valid for most PDE-constrained systems, and it leads to the existence of the posterior in Bayesian inverse problems. 
Proofs of Assumption~\ref{assumption:forward-continuity} for elliptic and parabolic equations can be found in~\cites{MR3084684, MR4246090}, while the proofs for the elasticity and Navier-Stokes equations can be found in~\cites{MR2652785, MR2558668}. 

We finally define the observational noise, $\eta$. 
We assume it to be Gaussian and independent of the parameters $\mathbf{z}$.
Therefore, $\eta$ is a random variable with values in $\mathbb{R}^k$ and it follows a normal distribution $\mathcal{N}(0, \Sigma)$, where $\Sigma$ is a known $k \times k$ symmetric positive covariance matrix. 

Having defined the parameters $\mathbf{z}$, the forward operator $\mathcal{G}$, and the observational noise $\eta$, along with the measurable space $(U, \Theta)$, and the prior distribution of the parameters $\gamma$, we now show the existence of the posterior, which we denote by $\gamma^y$.
This is possible thanks to Assumption~\ref{assumption:forward-continuity} that leads to $\gamma^{y}$ being absolutely continuous with respect to the prior $\gamma$. 
The posterior probability measure is defined through the Radon-Nikodym derivative
\begin{equation}\label{eq:radon-nikodym}
\frac{d \gamma^y}{d \gamma} \propto \exp(-\Phi(z; y)), 
\end{equation}
where $\Phi$ is known as the potential function, for example with Gaussian noise assumption, $\Phi(z, y) = \frac{1}{2}\|y - \mathcal{G}(z)\|_{\Sigma}^2$. Detailed proof of equation~\eqref{eq:radon-nikodym} can be found in \cite{MR2652785}. 

The last step to fully setup the framework is to show that the posterior measure is well-posed. 
This can be achieved following the results in~\cites{MR2652785, MR2558668,MR4069815}, that detail how the posterior measure is Lipschitz continuous with respect to the data under a certain distance metric. 
Specifically, for every $r > 0$ and $y, y' \in \mathbb{R}^d$ with $\|y\|_{\Sigma}, \|y'\|_{\Sigma} \leq r$, there exists $C = C(r) > 0$ such that
\begin{equation}\label{eq:continuous-well-posedness}
d_{H}(\gamma', \gamma'') = \Bigg(\frac{1}{2}\int_U\Bigg(\sqrt{\frac{d \gamma'}{d \gamma}} - \sqrt{\frac{d \gamma''}{d\gamma}}\Bigg)^2d \gamma\Bigg)^{\frac{1}{2}} \leq C(r)\|y-y'\|_{\Sigma}, 
\end{equation}
where $\gamma'$ and $\gamma''$ are two measures on $U$, which are absolutely continuous with respect to the measure $\gamma$, and where we chose as a distance metric the Hellinger distance $d_{H}$.
The latter was chosen to facilitate various proofs related to our new hybrid two-level MCMC approach, leading to Theorem~\ref{thm:final_error_estimate}, in Section~\ref{sec:two-level-mcmc-uniform}.

We note that the setup considered uses a uniform prior for the sake of simplicity. 
However, we can also work with a Gaussian setup, e.g. $U=\mathbb{R}^k$, $\Theta = \bigotimes_{j=1}^n \mathcal{B}(\mathbb{R})$, and $\gamma = \bigotimes_{j=1}^n \mathcal{N}(0, 1)$. 


\subsection{Approximation of the forward problem}
\label{section:approx-forward-problem}
A particularly expensive task in the Bayesian inverse problem setup introduced in Section~\ref{sec:setup} is the solution of the forward problem, especially when the forward operator $\mathcal{G}$ is constituted of PDEs. 
We distinguish two cases: (i) when the PDEs are approximated and solved via traditional numerical methods (e.g., FEM or others), also referred to as high-fidelity numerical models/solvers, and (ii) when the PDEs are solved via DL surrogates. 
We detail these two cases next, where, leveraging the results highlighted in Section~\ref{sec:setup}, we derive theoretical error estimates for each case, and make some observations on the computational costs.

\subsubsection{Traditional approximation methods}
\label{sec:traditional-approx}

The first case considered uses traditional numerical methods, such as FEM~\cite{parolini2005mathematical}, SEM~\cite{mengaldo2021industry}, or FVM~\cite{jameson1977finite}, to discretize the PDE system.
These numerical approximations lead to large linear systems of equations that are {extremely computationally expensive}, rendering the solution of forward problems impractical in the context of sampling-based statistical techniques, such as MCMC.
Yet, they provide well-defined theoretical error estimates, that typically depend on how fine the discretization (i.e., the tessellation of the computational domain via elements or grid points, also known as mesh) is.
This property is particularly useful in the context of inverse problems, because it allows practitioners to have a clear understanding of the errors incurred within their solution framework.

In particular, when considering any of the numerical methods above, we can define a priori estimates on the error we might expect for a certain discretization level $\ell$. 
The latter is {an integer value that specifies a characteristic}, $h$, that represents the dimensions of the elements (or spacing between grid points) tessellating the computational domain where the PDEs are being solved. 
For the purpose of this work, we assume an FEM-based discretization and the following generic error estimates.
\begin{assumption}
\label{assumption:numerical-error}
Let $u$ be the solution of the PDE equations in the forward problem.
We assume that $u \in V$, where $V$ is a suitable vector space, e.g. a Sobolev space. 
The FEM approximation error is given by
\begin{equation}
\|u - u_{\ell}\|_V \leq C 2^{-\ell},
\end{equation}
%
where $\ell$ is the level of discretization (each level $\ell$ halves the mesh size of the previous level $\ell - 1$) {and the corresponding mesh size is} $h = 2^{-\ell}$.
\end{assumption}
\begin{remark}
For simplicity, we did not include the error rate of time discretization. 
However, the time discretization error typically can be controlled by the discretization scheme to scale with the same rate of the spatial discretization. 
This will lead to the same convergence rate as in Assumption~\ref{assumption:numerical-error}.
Taking the finite-time two-dimensional Navier Stokes equation as an example, the error rate is $\|u(t) - u_{\ell}(t)\|_{H^1} \leq C {|h + \Delta t|}$ with a $\mathbb{Q}_1$-iso-$\mathbb{Q}_2/\mathbb{Q}_1$ mixed Finite element discretization and Implicit/Explicit (IMEX)
Euler time discretization scheme~\cites{MR4523340, MR2429876}.
\end{remark}
More details on error estimates for FEM can be found in~\cites{MR2050138,MR1043610}, while for SEM and FVM, the interested reader may refer to~\cite{karniadakis2005spectral} and~\cite{MR2897628}. 
In analogy to Equation~\eqref{eq:continuous-forward-map}, we can define the numerical approximation of the forward map as follows
\begin{equation}
\mathcal{G}^{\ell}: U \rightarrow \mathbb{R}^k\;\;\; \forall \mathbf{z} \in U; \;\;\; \mathcal{G}^{\ell}(\mathbf{z})= (\mathcal{F}_1(u_{\ell}(\mathbf{z}),\mathcal{F}_2(u_{\ell}(\mathbf{z})),  ...,\mathcal{F}_k(u_{\ell}(\mathbf{z}))).
\end{equation}  
where $u_{\ell}$ is the solution of the discrete forward problem and $\mathcal{F}_i$ for $i = 1, \dots, k$ are $k$ continuous bounded linear functionals.
Thanks to Assumption~\ref{assumption:forward-continuity}, we can write the posterior probability measure also for the discrete problem we are considering here (in analogy with the continuous counterpart in Equation~\eqref{eq:radon-nikodym})
\begin{equation}\label{eq:discrete-radon-nikodym}
    \frac{d \gamma^{\ell,y}}{d \gamma} \propto \exp(-\Phi^{\ell}(z; y)), 
\end{equation}
where $\Phi^{\ell}(z; y)$ is the discrete potential function.
Given Assumption~\ref{assumption:numerical-error} and Equation~\eqref{eq:continuous-well-posedness}, it follows immediately that the Hellinger distance metric between the continuous posterior $\gamma^y$ and the discrete one $\gamma^{y,\ell}$ is bounded for every numerical refinement level $\ell$
\begin{equation}\label{eq:discrete-hellinger-distance}
d_{H}(\gamma^{y}, \gamma^{\ell,y}) \leq C(y) 2^{-\ell},
\end{equation}
where $C(y)$ is a positive constant, that depends only of the data $y$.

Obviously, the larger the discretization level $\ell$ (i.e. the finer the mesh), the more computationally expensive the problem. 
In fact, the number of degrees of freedom of the corresponding discrete linear system increases exponentially with respect to $\ell$. 
Hence, achieving a solution with a desired (and ideally small) error might be out of reach even with abundant computational resources. 
DL-based approximation methods (also referred to as DL-based surrogates) can come to the rescue here, and are introduced next.

\subsubsection{DL-based approximation methods}
\label{sec:dl-based-approx}

The second case considered uses DL models to accelerate the solution of Bayesian inverse problems by replacing the computationally expensive numerical approximation just introduced in Section~\ref{sec:traditional-approx} with their {faster DL-based surrogate model counterparts}.
Let us denote $\Tilde{\mathcal{G}}: U \rightarrow \mathbb{R}^k$ as a nonlinear map defined by a trained DL model. 
We assume that the DL model is trained with data generated with classical numerical methods, e.g. FEM, and that the objective is to solve the inverse problem with an error less than or equal to $\mathcal{O}(2^{-L})$. 
The procedure for solving such an inverse problem with DL-based surrogate acceleration is as follows. 

First, we use a suitable numerical method (e.g. FEM, SEM, or FVM) to discretize the problem and generate the data.
We assume that we use a characteristic mesh size equal to $h = \mathcal{O}(2^{-L})$; to achieve the target accuracy, thanks to Assumption~\ref{assumption:numerical-error}.
Second, we use the generated numerical data as training data for the DL model. 
Third, we use the trained DL model as a surrogate to quickly run an MCMC chain. 
The estimated expectation of the QoI will be within the desired error if the DL-based surrogate model is as accurate as the numerical model. 

However, empirically, the trained DL model can hardly achieve the same level of accuracy as the numerical approximation used to generate the training data, and will lead to additional error.
We can formalize this statement as follows.
\begin{assumption}
\label{assumption:ml-model-error}
Given a DL model trained with data generated by a numerical approximation of the underlying forward problem that uses a mesh size $h = 2^{-L}$, and that has the error bound defined in Assumption~\ref{assumption:numerical-error}, we can write 
\begin{equation}\label{eq:ml-model-error}
\|\Tilde{G}(\mathbf{z}) - u\|_V \leq C 2^{-L + \epsilon},
\end{equation}
where $\mathbf{z}$ is the input, and $\epsilon$ accounts for the error of the DL model. 
In order for the DL error to be small, we require a small value of $\epsilon$. 
\end{assumption}
In practice, we expect $\epsilon$ to be small when we have a reasonably good DL model trained with sufficient data.
However, in general, direct replacement of the numerical solver with a DL-based surrogate will lead to a posterior distribution estimation error of $\mathcal{O}(2^{-L+\epsilon})$ given by Equation~\eqref{eq:discrete-hellinger-distance}. 
Hence, producing an error gap to the targeted $\mathcal{O}(2^{-L})$.

In order to mitigate the shortcomings of DL-based surrogate models, one can refine the mesh of the numerical model used for data generation and increase the size of the training data to produce a possibly more accurate DL model that can reach the desired accuracy. 
However, that will increase the computational cost by many folds. 
In general, to solve a two-dimensional problem, the minimum increment of the computational cost of the numerical model is 4 times, and 8 times for three-dimensional problems, not to mention more challenging problems whose computational cost does not scale linearly with the degrees of freedom.
In addition to the cost of finer numerical solvers, the larger amount of data will also increase the DL training cost. 
Even if we are willing to pay the cost, it was shown that there is an empirical limit to the accuracy that certain DL models can reach~\cite{CAO2023112104}. 
Therefore, in some cases, the desired accuracy may be unreachable by direct substitution of the numerical model with a DL surrogate.
In the next section, we propose a different approach to correct the error statistically with the MCMC method. 

\subsection{Hybrid two-level MCMC with uniform prior}
\label{sec:two-level-mcmc-uniform}
In this section, we propose the hybrid two-level MCMC method for error correction of the DL-based surrogate model for Bayesian inverse problems. 
The new approach is inspired by the multilevel version of MCMC, which was shown to reduce the computational cost of standard MCMC for various problems by two orders of magnitude~\cites{MR3372290, MR3084684, MR4246090, MR4523340}. 
Our hybrid two-level method utilizes both the DL-based surrogate model and the high-fidelity numerical model to sample the posterior probability of Bayesian inverse problems.
In particular, we run a base MCMC chain with a DL-based surrogate, and a {short} correction MCMC chain with a numerical model with known accuracy.
The latter is deployed to correct the statistical error of the MCMC chain generated by the DL-based surrogate.

Numerical multi-level approaches have been very successful for multi-scale physical problems. 
However, implementing multi-level algorithm for generic engineering or scientific problems can be challenging due to complex meshes and instability of coarse numerical models. 
We see the potential to avoid those challenges by hybridizing the DL-based surrogate and numerical solvers under the same mathematical framework.
Hence, we propose a two-level hybrid approach inspired by the telescoping argument of the multilevel Monte Carlo algorithm \cite{MR2436856}. 
We note that another potential approach is to use the DL-based surrogate model as a filter for the MCMC sampler, as inspired by~\cite{MR2231730}.
{However, our proposed approach can be an alternative approach and potentially generalizable beyond the MCMC algorithm (i.e., our approach may also be applied to other methods such as the ensemble Kalman filter and sequential Monte Carlo methods, among others where a telescoping sum structure is applicable).}

We now start introducing our hybrid two-level MCMC method. 
To this end, we denote the $Q(z)$ as $Q$ for simplicity. We further indicate the posterior distribution approximated by the DL-based surrogate model as $\gamma^{\rm DL}$, and the numerically approximated posterior distribution as $\gamma^{\rm num}$. 
With the target precision of $\mathcal{O}(2^{-L})$ and Assumption~\ref{assumption:ml-model-error}, $\gamma^{\rm num}$ and $\gamma^{\rm DL}$ are equivalent to $\gamma^{L, y}$ and $\gamma^{L-\epsilon, y}$ in Section \ref{section:approx-forward-problem}. 
In our two-level approach, we can rewrite the numerical approximation of the expected QoI $Q$ as follows
\begin{align}
\label{eq:two-level-expansion}
\mathbb{E}^{\gamma^{\rm num}}[Q] &= \mathbb{E}^{\gamma^{\rm num}}[Q] - \mathbb{E}^{\gamma^{\rm DL}}[Q]  + \mathbb{E}^{\gamma^{DL}}[Q] \nonumber \\
&= \left(\mathbb{E}^{\gamma^{\rm num}}-\mathbb{E}^{\gamma^{\rm DL}}\right)[Q] + \mathbb{E}^{\gamma^{DL}}[Q]. 
\end{align}  
To derive a computable estimator with MCMC chains, we observe that the first term on the {right} hand side in~\eqref{eq:two-level-expansion} can be transformed as follows
\begin{align}
\label{equation: two level sum}
\left(\mathbb{E}^{\gamma^{\rm num}}-\mathbb{E}^{\gamma^{\rm DL}}\right)[Q] &= \frac{1}{N^{\rm num}}\int_{U} \exp(-\Phi^{\rm num}) Q d \gamma - \frac{1}{N^{\rm DL}}\int_{U} \exp(-\Phi^{\rm DL}) Q d \gamma \nonumber \\
&= \frac{1}{N^{\rm num}}\int_{U} (\exp(-\Phi^{\rm num}) - \exp(-\Phi^{\rm DL})) Q d \gamma \nonumber \\
&+\left(\frac{1}{N^{\rm num}} - \frac{1}{N^{\rm DL}}\right)\int_{U} \exp(-\Phi^{\rm DL}) Q d\gamma \nonumber \\
&= \frac{1}{N^{\rm num}}\int_U \exp(-\Phi^{\rm num})(1-\exp(\Phi^{\rm num}-\Phi^{\rm DL})) Q d\gamma \nonumber \\
&+ \left(\frac{N^{\rm DL}}{N^{\rm num}}-1\right) \frac{1}{N^{\rm DL}} \int_U \exp(-\Phi^{\rm DL}) Q d \gamma,
\end{align} 
{where $N^{\rm num} = \int_U \exp(-\Phi^{\rm num})d\gamma$ and $N^{\rm DL} = \int_U \exp(-\Phi^{\rm DL})d\gamma$ are the normalization constants.}
The constant $(N^{\rm DL}/N^{\rm num}-1)$ can be expanded as 
\begin{equation}
\label{equation: mcmc for additional constant}
\bigg(\frac{N^{\rm DL}}{N^{\rm num}}-1\bigg) = \frac{1}{N^{\rm num}} \int_U (\exp(\Phi^{\rm num}-\Phi^{\rm DL})-1) \exp(-\Phi^{\rm num}) d \gamma. 
\end{equation} 
{We note that the integral $\frac{1}{N^{\rm num}}\int_U (\cdot) \exp(-\Phi^{\rm num}) d\gamma$ and $\frac{1}{N^{\rm DL}}\int_U (\cdot) \exp(-\Phi^{\rm DL}) d\gamma$ can be estimated with an MCMC estimator $\mathbf{E}_{M_{\rm num}}^{\gamma^{\rm num}}[\cdot]$ and $\mathbf{E}_{M_{\rm DL}}^{\gamma^{\rm DL}}[\cdot]$, where the $M_{\rm num}$ and $M_{\rm DL}$ is the number of numerical MCMC samples and DL surrogate MCMC samples.} Having defined equations~\eqref{eq:two-level-expansion}, \eqref{equation: two level sum}, and \eqref{equation: mcmc for additional constant}, we can write

\begin{align}
\label{equation: two level expansion explicit}
\mathbb{E}^{\gamma^{\rm num}}[Q] &= \mathbb{E}^{\gamma^{\rm num}}[(1-\exp(\Phi^{\rm num}-\Phi^{\rm DL}))Q] \nonumber \\
&+ \mathbb{E}^{\gamma^{\rm num}}[\exp(\Phi^{\rm num}-\Phi^{\rm DL})-1] \cdot \mathbb{E}^{\gamma^{\rm DL}}[Q] + \mathbb{E}^{\gamma^{\rm DL}}[Q], 
\end{align}  
and we can now define the hybrid two-level MCMC estimator $\mathbf{E}^{\rm hybrid}[Q]$ of $\mathbb{E}^{\gamma^y}[Q]$ as follows
\begin{align}
\label{eq:hybrid-estimator}
\mathbf{E}^{\rm hybrid}[Q] &= \mathbf{E}^{\gamma^{\rm num}}_{M_{\rm num}}[(1-\exp(\Phi^{\rm num}-\Phi^{\rm DL}))Q] \nonumber \\
&+ \mathbf{E}^{\gamma^{\rm num}}_{M_{\rm num}}[\exp(\Phi^{\rm num}-\Phi^{\rm DL})-1] \cdot \mathbf{E}^{\gamma^{\rm DL}}_{M_{\rm DL}}[Q] + \mathbf{E}_{M_{\rm DL}}^{\gamma^{\rm DL}}[Q].
\end{align}  
The new hybrid two-level MCMC approach for Bayesian inverse problems just introduced is simple, and it can therefore be adopted easily with legacy numerical models without too many changes in the code base.
An important aspect of the new hybrid two-level MCMC introduced in equation~\eqref{eq:hybrid-estimator} is its error analysis. 
In particular, we aim to show how the correction chain effectively corrects the estimator error caused by DL-based surrogate model.

\begin{thm}
\label{thm:error decomposition}

The hybrid two-level MCMC estimator error can be decomposed into the following three components:  
\begin{subequations}\label{eq:error-analysis}
\begin{align}
& \mathbb{E}^{\gamma^y}[Q] - \mathbf{E}^{\rm hybrid}[Q] = \mathrm{I} + \mathrm{II} + \mathrm{III}, \label{eq:error} \\[0.6em]
\text{where}\;\;\;  
& \mathrm{I} := \mathbb{E}^{\gamma^y}[Q] - \mathbb{E}^{\gamma^{\rm num}}[Q] \label{eq:error-I} \\[0.5em]
& \mathrm{II} := \mathbb{E}^{\gamma^{DL}}[Q] - \mathbf{E}^{\gamma^{DL}}_{M_{\rm DL}}[Q] \label{eq:error-II} \\[0.5em]
& \mathrm{III} := \mathbb{E}^{\gamma^{\rm num}}[(1-\exp(\Phi^{\rm num}-\Phi^{\rm DL}))Q] \nonumber \\[0.1em]
& \;\;\;\;\; - \mathbf{E}^{\gamma^{\rm num}}_{M_{\rm num}}[(1-\exp(\Phi^{\rm num}-\Phi^{\rm DL}))Q] \nonumber \\[0.1em]
& \;\;\;\;\; + \mathbb{E}^{\gamma^{\rm num}}[\exp(\Phi^{\rm num}-\Phi^{\rm DL})-1] \cdot \mathbb{E}^{\gamma^{\rm DL}}[Q] \nonumber \\[0.1em] 
& \;\;\;\;\; - \mathbf{E}^{\gamma^{\rm num}}_{M_{\rm num}}[\exp(\Phi^{\rm num}-\Phi^{\rm DL})-1] \cdot \mathbf{E}^{\gamma^{\rm DL}}_{M_{\rm DL}}[Q] \label{eq:error-III}
\end{align}
\end{subequations}
\end{thm}
\begin{proof}
Given the estimator error, we observe
\begin{align}
\mathbb{E}^{\gamma^y}[Q] - \mathbf{E}^{\rm hybrid}[Q] &= \mathbb{E}^{\gamma^y}[Q] - \mathbb{E}^{\gamma^{\rm num}}[Q] + \mathbb{E}^{\gamma^{\rm num}}[Q]- \mathbf{E}^{\rm hybrid}[Q] \nonumber \\ 
&= \mathrm{I} + \mathbb{E}^{\gamma^{\rm num}}[Q] - \mathbf{E}^{\rm hybrid}[Q]\label{eq:proof-1}
\end{align}
With equation \eqref{equation: two level expansion explicit} and \eqref{eq:hybrid-estimator}, we have
\begin{align}
\mathbb{E}^{\gamma^y}[Q] - \mathbf{E}^{\rm hybrid}[Q] &= \mathrm{I} + \mathbb{E}^{\gamma^{\rm num}}[Q] - \mathbf{E}^{\rm hybrid}[Q] \nonumber \\
&= \mathrm{I} + \mathbb{E}^{\gamma^{DL}}[Q] + \mathbb{E}^{\gamma^{\rm num}}[(1-\exp(\Phi^{\rm num}-\Phi^{\rm DL}))Q] \nonumber \\
&+ \mathbb{E}^{\gamma^{\rm num}}[\exp(\Phi^{\rm num}-\Phi^{\rm DL})-1] \cdot \mathbb{E}^{\gamma^{\rm DL}}[Q] - \mathbf{E}^{\gamma^{DL}}_{M_{\rm DL}}[Q] \nonumber \\
&- \mathbf{E}^{\gamma^{\rm num}}_{M_{\rm num}}[(1-\exp(\Phi^{\rm num} - \Phi^{\rm DL}))Q] \nonumber \\
&- \mathbf{E}^{\gamma^{\rm num}}_{M_{\rm num}}[\exp(\Phi^{\rm num} - \Phi^{\rm DL})-1] \cdot \mathbf{E}^{\gamma^{\rm DL}}_{M_{\rm DL}}[Q] \nonumber \\
& = \mathrm{I} + \mathrm{II} + \mathrm{III}. \label{eq:proof-2}
\end{align}
\end{proof}
As shown, the overall error bound for our hybrid two-level MCMC method is composed of three error terms in equation~\eqref{eq:error-analysis}.
We analyse each error term individually, and assemble the overall error result as a conclusion to this analysis.
\begin{proposition}
\label{prop:sampling mean square error bound}

Let $\mathcal{C} = \{\mathbf{z}^{(n)}\}_{n\in\mathbb{N}}$ be a Markov chain, $\mathcal{P}$ be the probability measure of the Markov chain. For every bounded $Q$ and every $M \in \mathbb{N}$, we have the following mean square error bound:
\begin{equation*}
    (\mathcal{E}[|\mathbb{E}[Q(z)] - \frac{1}{M}\sum_{n=1}^{M}Q(z^{(n)})|^2])^{1/2} \leq C \sup_{z\in U} |Q(z)| M^{-1/2},
\end{equation*}
where $\mathcal{E}$ is the expectation over all realizations of $\mathcal{C}$ with respect to the measure $\mathcal{{P}}$. 
\end{proposition}
{This is a standard result from Markov chain theory, detailed proof can be found in~\cites{levin2017markov, MR3084684}.
For the flow of the paper, we include the proposition here without proof.}
\begin{thm}
\label{thm:final_error_estimate}
{We denote by $\mathbf{C}_{\rm hybrid} = \{\mathcal{C}_{\rm num}, \mathcal{C}_{\rm DL}\}$ the collection of Markov chains obtained with numerical forward solver and DL-based surrogate solver 
Let $\mathcal{P}^{\rm num}$ and $\mathcal{P}^{\rm DL}$ be the probability measure of respective Markov chains, we denote $\mathbf{P}_{\rm hybrid} = \mathcal{P}^{\rm num} \bigotimes \mathcal{P}^{\rm DL}$.} 
With $M_{\rm DL} = {C_{\rm DL}} 2^{2L}$ and $M_{\rm num} = {C_{\rm num}}(1+2^{\epsilon})^2$, we have the following theoretical error estimate  of our hybrid two-level MCMC approach under uniform priors
\begin{equation}
\mathcal{E}_{\rm hybrid}[|\mathbb{E}^{\gamma^y}[Q] - \mathbf{E}^{\rm hybrid}[Q]|] \leq {C_{\rm hybrid}} 2^{-L}, 
\end{equation}
{where $\mathcal{E}_{\rm hybrid}$ is the expectation over all realizations of the collection $\mathbf{C}_{\rm hybrid}$ with respect to the product meansure $\mathbf{P}_{\rm hybrid}$.}
\end{thm}
\begin{proof}

From Theorem \ref{thm:error decomposition}, we decompose the overall error into three components. 
For error term I, from equation~\eqref{eq:error-analysis} and equation~\eqref{eq:discrete-hellinger-distance}, we can obtain the following error bound
\begin{equation}
    |\mathrm{I}| := |\mathbb{E}^{\gamma^y}[Q] - \mathbb{E}^{\gamma^{\rm num}}[Q]| {\leq 2(\mathbb{E}^{\gamma^y}[Q^2] + \mathbb{E}^{\gamma^{\rm num}}[Q^2])^{1/2}d_H(\gamma^y, \gamma^{\rm num})}\leq C 2^{-L},
\end{equation}

\noindent {where the details of the first inequality can be found in~\cite{MR2652785}.} For error term II, from {Proposition~\ref{prop:sampling mean square error bound}} we can obtain the following error bound
\begin{equation}
\label{eqn:error_term_2}
    \mathcal{E}_{\rm DL}[|\mathrm{II}|] \leq (\mathcal{E}_{\rm DL}[|\mathrm{II}|^2])^{1/2} := (\mathcal{E}_{\rm DL}[|\mathbb{E}^{\gamma^{DL}}[Q] - \mathbf{E}^{\gamma^{DL}}_{M_{\rm DL}}[Q]|^2])^{1/2} \leq C M_{DL}^{-1/2},
\end{equation}
{where $\mathcal{E}_{\rm DL}$ is the expectation over all realizations of Markov chain $\mathcal{C_{\rm DL}}$ with respect to the probability measure $\mathcal{P}^{DL}$. }

\noindent Finally, for error term III, we can use the inequality $|\exp(x) - \exp(y)| \leq |x-y|(\exp(x)+\exp(y))$, to obtain
\begin{align}
\sup_{z \in U} |1 - \exp(\Phi^{\rm num} - \Phi^{\rm DL})| &\leq \sup_{z \in U} |\Phi^{\rm num} - \Phi^{\rm DL}|(1+\exp(\Phi^{\rm num} - \Phi^{\rm DL})) \nonumber \\
&  \leq C \sup_{z \in U}(||y - \mathcal{G}^{\rm num}|^2 - |y - \mathcal{G}^{\rm DL}|^2|) \nonumber \\
&  \leq C \sup_{z \in U}(2|y| + |\mathcal{G}^{\rm num}| + |\mathcal{G}^{\rm DL}|)|\mathcal{G}^{\rm num} - \mathcal{G}^{DL}| \nonumber \\
&  \leq C \sup_{z \in U} (|\mathcal{G}^{\rm num} - \mathcal{G}| + |\mathcal{G}^{\rm DL} - \mathcal{G}|) \nonumber \\
&  \leq C (2^{-L} + 2^{-L+\epsilon}) \nonumber \\
& \leq C (1+2^\epsilon)2^{-L},
\end{align}  
that leads to
\begin{align}
\label{eq:first_part}
&\mathcal{E}_{\rm num}[\{\mathbb{E}^{\gamma^{\rm num}}[(1-\exp(\Phi^{\rm num}-\Phi^{\rm DL}))Q] -\mathbf{E}^{\gamma^{\rm num}}_{M_{\rm num}}[(1-\exp(\Phi^{\rm num}-\Phi^{\rm DL}))Q]\}^2] \nonumber \\
&  \leq C M_{\rm num}^{-1}\sup_{z\in U}(|1-\exp(\Phi^{\rm num}-\Phi^{\rm DL})|^2) \nonumber \\
& \leq C M_{\rm num}^{-1} (1+2^\epsilon)^2 2^{-2L}.
\end{align}  
Similarly, we have 
\begin{align}
\label{eq:second_part}
&\mathcal{E}_{\rm hybrid}[\{\mathbb{E}^{\gamma^{\rm num}}[\exp(\Phi^{\rm num} - \Phi^{\rm DL})-1] \cdot \mathbb{E}^{\gamma^{\rm DL}}[Q] \nonumber\\ 
& - \mathbf{E}^{\gamma^{\rm num}}_{M_{\rm num}}[\exp(\Phi^{\rm num}-\Phi^{\rm DL})-1] \cdot \mathbf{E}^{\gamma^{\rm DL}}_{M_{\rm DL}}[Q]\}^2] \nonumber \\[0.5em]
\leq & \;{C} \mathcal{E}_{\rm num}[\{\mathbb{E}^{\gamma^{\rm num}}[\exp(\Phi^{\rm num}-\Phi^{\rm DL})-1] \nonumber\\ 
& - \mathbf{E}^{\gamma^{\rm num}}_{M_{\rm num}}[\exp(\Phi^{\rm num}-\Phi^{\rm DL})-1]\}^2] \cdot \sup_{z\in U}|Q|^2 \nonumber\\
& + {C} \sup_{z \in U}|\exp(\Phi^{\rm num}-\Phi^{\rm DL})-1|^2 \cdot \mathcal{E}_{\rm DL}[\{ \mathbb{E}^{\gamma^{\rm DL}}[Q] -\mathbf{E}^{\gamma^{\rm DL}}_{M_{\rm DL}}[Q]  \}^2] \nonumber \\[0.5em]
\leq & {C} M_{\rm num}^{-1} \sup_{z\in U}|\exp(\Phi^{\rm num} - \Phi^{\rm DL}) -1|^2+ {C} M_{\rm DL}^{-1} \sup_{z\in U}|\exp(\Phi^{\rm num} - \Phi^{\rm DL}) -1|^2. \nonumber \\
\leq & \;{C} M_{\rm num}^{-1} (1+2^{\epsilon})^2 2^{-2L} + {C} M_{\rm DL}^{-1} (1+2^{\epsilon})^2 2^{-2L}.
\end{align}
By combining equations~\eqref{eq:first_part} and  \eqref{eq:second_part}, we have the overall error estimate for III, that is
\begin{align}
    \mathcal{E}_{\rm hybrid}|\mathrm{III}|^2 &:= \mathcal{E}_{\rm hybrid}[|\mathbb{E}^{\gamma^{\rm num}}[(1-\exp(\Phi^{\rm num} -\Phi^{\rm DL}))Q] \nonumber\\ 
     &\;\; -  \mathbf{E}^{\gamma^{\rm num}}_{M_{\rm num}}[(1-\exp(\Phi^{\rm num}-\Phi^{\rm DL}))Q] \nonumber\\
    &\;\;+ \mathbb{E}^{\gamma^{\rm num}}[\exp(\Phi^{\rm num}-\Phi^{\rm DL})-1] \cdot \mathbb{E}^{\gamma^{\rm DL}}[Q]  \nonumber\\ 
    &\;\;- \mathbf{E}^{\gamma^{\rm num}}_{M_{\rm num}}[\exp(\Phi^{\rm num}-\Phi^{\rm DL})-1] \cdot \mathbf{E}^{\gamma^{\rm DL}}_{M_{\rm DL}}[Q]|^2] \nonumber \\[0.5em]
    & \hspace{-0.4cm}\leq {C} M_{\rm num}^{-1} (1+2^{\epsilon})^2 2^{-2L} + {} M_{\rm DL}^{-1} (1+2^{\epsilon})^2 2^{-2L}.
\end{align}  
Until now, we are still free to choose the number of samples for $M_{\rm num}$ and $M_{\rm DL}$. 
To balance errors I, II and III, we can choose the sampling number $M_{\rm DL} = {C_{DL}} 2^{2L}$ (see Equation~\ref{eqn:error_term_2}) and $M_{\rm num} = {C_{\rm num}} (1+2^\epsilon)^2$, that allows us to write the overall error estimate. 
{Here, we differentiate the two constants of different values by $C_{\rm num}$ and $C_{\rm DL}$.} 
Finally we have a final numerical to DL sample ratio of ${\frac{C_{\rm num}}{C_{\rm DL}}} (1+2^\epsilon)^2/2^{2L}$.
\end{proof}
{With the DL-based surrogate sample number $M_{\rm DL} = C_{\rm DL} 2 ^{2L}$ and numerical sample number $M_{\rm num} = C_{\rm num} (1+2^{\epsilon})$, where both $C$ does not depend on $L$ and $\epsilon$, 
Theorem~\ref{thm:final_error_estimate} shows that there is an optimal ratio $C_{\rm num} (1+2^\epsilon)^2/C_{\rm DL}2^{2L}$ for our hybrid method to reach the same accuracy as plain MCMC with a numerical model. However it is theoretically non-trivial to work out both constants' values. 
For high fidelity problems which benefit the most from the adoption of a DL-based surrogate model, a large L value is expected in which the $(1+2^\epsilon)^2/2^{2L}$ term will be the dominant factor leading to a small overall ratio. 
Nevertheless, the theory is only good enough to give a general guideline. 
Exact optimal ratio is not available analytically with dependencies on factors other than $L$ and $\epsilon$, such as regularity of specific problems. 
We propose a sweep test for our numerical experiments and potential practical applications, where we empirically test a range of different ratios to get the best hybrid configuration. 
In addition, if we assume the computational speedup rate of the DL-based surrogate model as $s=t_{\rm num}/t_{\rm DL}$, where $t_{\rm num}$ is the average computational time for one forward solve with numerical model and $t_{\rm DL}$ is the computational time of one DL-based surrogate evaluation time}, we can estimate the overall speedup of our new hybrid two-level approach compared to the plain numerical MCMC. 
In particular, with the choice of samples $M_{\rm DL}$ and $M_{\rm num}$, the overall speedup is {$\mathcal{O}( 2^{2L}/(\frac{1}{s}2^{2L} + \frac{C_{\rm num}}{C_{\rm DL}}(1+2^{\epsilon})^2))$} compared to the plain numerical MCMC. 
We note that the hybrid two-level approach can be easily parallelized with two processes, thus the actual speedup is {$\mathcal{O}( 2^{2L}/\max(\frac{1}{s}2^{2L}, \frac{C_{\rm num}}{C_{\rm DL}}(1+2^{\epsilon})^2) )$}. 
Even though the discussion is based on the setup of a uniform prior, the same conclusions hold for Gaussian priors. 
Details of the latter are provided in~\ref{app:two-level-mcmc-gaussian}. 
In Section~\ref{sec:experiments}, we show several experiments that validate the theoretical findings.

As a final note, we shall point out that, even though the hybrid method is able to provide accurate posterior expectations of the quantities of interest and the variance (which can be computed from the expectations), the method will not generate a large number of actual numerical samples from the highly accurate numerically approximated posterior distribution (the exact reason why computational cost is saved). 
This limits the method from producing a histogram like a conventional MCMC chain. However, a less accurate histogram can still be generated from the DL surrogate samples.

\section{Numerical experiments}
\label{sec:experiments}
After having introduced the new hybrid two-level MCMC approach, we now present some numerical experiments to validate the theoretical error estimates and to understand the computational performance of the new approach.
{The numerical experiments span three different PDE problems, namely a Poisson equation, a nonlinear reaction-diffusion  equation, and a Navier-Stokes equation.}
These three problems include both elliptic and parabolic differential problems, with different levels of complexity, and constitute an established benchmark for testing novel methods in the context of Bayesian inverse problems~\cite{VillaPetraGhattas21}.

In the numerical experiments, we consider both uniform and Gaussian priors, to demonstrate the theoretical error estimates presented in Sections~\ref{sec:two-level-mcmc-uniform}, and~\ref{app:two-level-mcmc-gaussian}, respectively. 
\begin{figure}[H]
    \centering
    \subfloat[Mesh A: 32 $\times$ 32 resolution]{
        \includegraphics[width=0.44\textwidth]{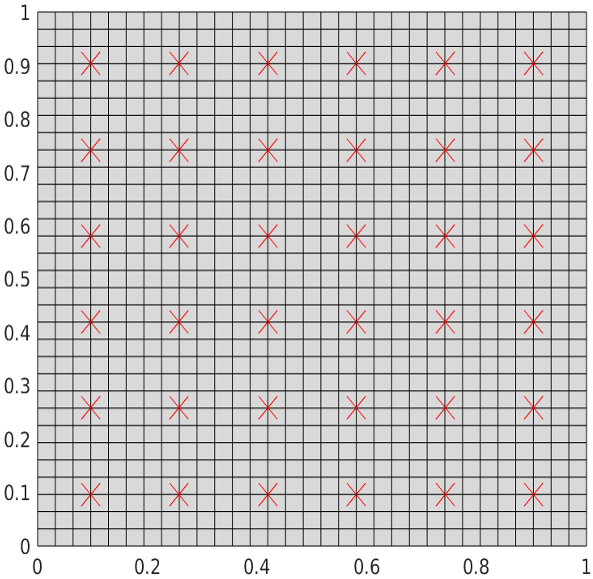}
        \label{fig:mesh1}
    }
    \hfill
    \subfloat[Mesh B: 64 $\times$ 64 resolution]{
        \includegraphics[width=0.44\textwidth]{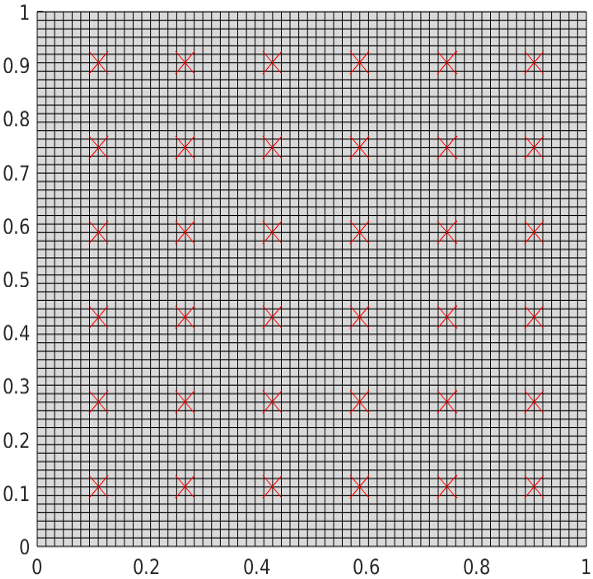}
        \label{fig:mesh2}
    }
    \caption{Meshes adopted and points where the solution is observed (red crosses). Mesh A is used for the numerical experiments on a Poisson and a nonlinear reaction-diffusion  equations, while Mesh B is used for the numerical experiments on a Navier-Stokes equation.}
    \label{fig:meshes}
\end{figure}
All results are obtained in the two-dimensional squared computational domains (or meshes) $D_h \in [0, 1] \times [0, 1]$, Mesh A and Mesh B, depicted in Figure~\ref{fig:meshes}, where we observe the solution $u$ in equally spaced fixed positions, highlighted as red crosses.
Mesh A in Figure~\ref{fig:mesh1} is used for the numerical experiments with a Poisson equation and a nonlinear reaction-diffusion equation, and has a number of cells equal to $32 \times 32 = 1024$ (mesh level $\ell=5$).
Mesh B in Figure~\ref{fig:mesh2} is used for the numerical experiment with a Navier-Stokes equation, and has a number of cells $64 \times 64 = 4096$ (mesh level $\ell=6$).

For the different problems considered, we used a range of different neural network architectures, to show that the method proposed here is agnostic to the choice of DL surrogate. 
The list of experiments carried out in the following subsections is summarized in Table~\ref{tab:list_of_experiments}. 
\begin{table}[ht]
    \centering
    {
    \footnotesize
    \begin{tabular}{|c|c|c|c|c|c|c|}
    \hline
    {}  & 1 & 2 & 3 & 4 & 5 & 6\\ \hline
    Problem  & \multicolumn{2}{|c|}{Poisson Equation} & \multicolumn{2}{|c|}{Reaction Diffusion} & \multicolumn{2}{|c|}{Navier Stokes} \\ \hline
    Prior & Uniform & Gaussian & Uniform & Gaussian & Uniform & Gaussian \\ \hline
    DL Model & FCN & CNN & GNN & U-Net & DeepONet & FNO \\ \hline
    Section & \ref{sec:poisson-uniform} & \ref{sec:poisson-gauss} & \ref{sec:reaction-diffusion-uniform} & \ref{sec:reaction-diffusion-uniform} & \ref{sec:navier-stokes-uniform} & \ref{sec:navier-stokes-uniform} \\ \hline
    \end{tabular}
    }
    \caption{List of numerical experiments}
    \label{tab:list_of_experiments}
\end{table}
{We use Metropolis-Hasting algorithm with prior distribution as our proposal density for all our MCMC chains in the numerical experiments.}

\subsection{Poisson equation}
\label{subsec:poisson}
We first consider a Bayesian inverse problem with a forward model governed by the Poisson equation in the two-dimensional computational domain depicted in Figure~\ref{fig:mesh1}.
\begin{equation}
\begin{cases}
\nabla \cdot (K(z) \nabla u(x)) =  \cos(2\pi x_1) \sin(2 \pi x_2), \\[0.4em]
u(x_1=0) = 0, \\[0.4em]
u(x_1=1) = 1, \\[0.4em]
\displaystyle \frac{\partial u}{\partial x_2} (x_2 = 0) = \frac{\partial u}{\partial x_2} (x_2 = 1) = 0.
\end{cases} 
\label{eq:poisson}
\end{equation}
The data $u$ is observed at thirty-six equally-distanced positions, as shown in figure~\ref{fig:mesh1}, using a random realization of the forward model with additive Gaussian noise $\delta$ that has zero mean and variance $\sigma^2 = 0.001$. 
We next present the uniform and Gaussian prior cases.

\subsubsection{Uniform prior}
\label{sec:poisson-uniform}
\sloppy For this first case, we set the QoI to be the random field $K = z\cos(2\pi x_1)\sin(2\pi x_2)+2.0$, whereby its prior distribution is uniform, namely
$z \sim U[0, 1]$.  
We solve equation~\eqref{eq:poisson} using FEM on Mesh A (depicted in Figure~\ref{fig:mesh1}), and randomly generate 8000 solution samples. 
The 8000 samples are partitioned into 4000 training samples, 2000 validation samples, and 2000 test samples to train a fully connected ReLU neural network.
For details on the FEM solver and the neural network model, the interested reader may refer to \ref{section:poisson_solver}.

\begin{figure}
    \centering
    \includegraphics[width=0.8\linewidth]{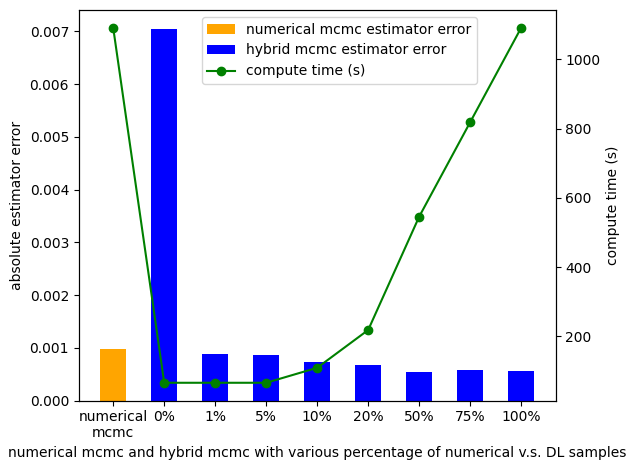}
    \caption{Estimator error with different percentage of numerical sample against DL-based surrogate samples for the elliptic experiment with uniform prior}
    \label{fig:elliptic uniform percentage num samples}
\end{figure}

We performed numerical experiments on three setups: (i) a plain MCMC chain with numerical solver, denoted as \textbf{Numerical MCMC} in Table~\ref{tab:elliptic_uniform_results}, (ii) a plain MCMC chain with DL-based surrogate model, denoted as \textbf{DL MCMC} in Table~\ref{tab:elliptic_uniform_results}, (iii) and the proposed hybrid approach, denoted as \textbf{Hybrid MCMC} in Table~\ref{tab:elliptic_uniform_results}, where we tested different lengths of the numerical samples chain, ranging from 1\% to 100\% of the total number of DL samples (see  Fig~\ref{fig:elliptic uniform percentage num samples}). 
We show in detail two cases in Table~\ref{tab:elliptic_uniform_results}, one with $1\%$ numerical samples and the other with $5\%$ numerical samples (compared to the total number of DL samples). 
All values showed in Table~\ref{tab:elliptic_uniform_results} are average results of eight MCMC runs.

\begin{table}[H]
\footnotesize
\centering
\begin{tabular}{l|c|c|c|c}
\toprule
\textbf{Method} &  
\begin{tabular}{@{}c@{}} \textbf{Numerical}  \\ \textbf{MCMC} \end{tabular} & 
\begin{tabular}{@{}c@{}} \textbf{DL}         \\ \textbf{MCMC} \end{tabular} & 
\begin{tabular}{@{}c@{}} \textbf{Hybrid}     \\ \textbf{MCMC} \end{tabular} & 
\begin{tabular}{@{}c@{}} \textbf{Hybrid}     \\ \textbf{MCMC} \end{tabular} \\ 
\hline
\textbf{Samples}  & 100,000 & 100,000 & 
\begin{tabular}{@{}c@{}} 100,000 \\ +1,000 \end{tabular} & 
\begin{tabular}{@{}c@{}} 100,000 \\ +5,000 \end{tabular} \\ 
\midrule
\begin{tabular}{@{}l@{}} \textbf{Estimator} \\ \textbf{error} \end{tabular} 
& 9.72E-4 & 7.524E-3 & 8.75E-4 & 8.69E-4 \\ 
\midrule
\begin{tabular}{@{}l@{}} \textbf{Compute} \\ \textbf{time [s]} \end{tabular} 
& 1090.95 & 65.77 & \begin{tabular}{@{}l@{}} 66.68 (serial) \\ 65.77 (parallel) \end{tabular}  & \begin{tabular}{@{}l@{}} 120.32 (serial) \\ 65.77 (parallel) \end{tabular}  \\ 
\midrule
\textbf{Speed up} & - & 16.59x & \begin{tabular}{@{}l@{}} 16.36 (serial) \\ 16.59x (parallel) \end{tabular}  & \begin{tabular}{@{}l@{}} 9.07x (serial) \\ 16.59x (parallel) \end{tabular}  \\ 
\bottomrule
\end{tabular}
\caption{Estimator error and compute speedup for the elliptic problem with uniform prior}
\label{tab:elliptic_uniform_results}
\end{table}
{The reference value computed employing the FEM model and approximating the expected value of the posterior with a quadrature scheme with 32 Gaussian Legendre quadrature over the parametric domain $z\in[0, 1]$} is also included.
Quadrature estimation of the posterior mean shows a value of 0.313237 with mesh {level $\ell = 10 (1024 \times 1024 \, \rm{cells})$}. 
This provides a highly accurate posterior expectation which can be considered as a true reference. 
Hence, all the estimator errors are calculated by comparing to this reference.
We observe how the proposed hybrid two-level MCMC approach provides results that are comparable to the numerical MCMC method.
The comparison between the two hybrid MCMC experiments also qualitatively validates Theorem~\ref{thm:final_error_estimate}, where an optimal ratio exists and additional numerical samples do not further improve the results.
In addition to the results presented in Table~\ref{tab:elliptic_uniform_results}, we also experimentally estimate the mean of the surrogate model error and numerical model error with respect to z again by a Gaussian Legendre quadrature (this time with 64 quadrature points), that in turn allows us to approximate $\epsilon$ in Assumption~\ref{assumption:numerical-error}. 
Using again the mesh at level $\ell=10$ as reference, we obtain $\mathbf{E}\|u_{L=10}-u_{L=5}\|_{L^{\infty}(D)} \approx 2.0E-4$, and of
$\mathbf{E}\|u_{L=10}-u_{L=5}\|_{L^{\infty}(D)} \approx 1.0E-3$ for the numerical forward solver and for the DL surrogate, respectively.
With such, we have a rough estimate of $\epsilon=2.3$. 
This is not a rigorous error bound estimate; nevertheless, it still provides a useful indication of the accuracy of the DL-based surrogate model.
Follow Theorem~\ref{thm:final_error_estimate}, our estimate provides the ratio {$(1+2^\epsilon)^2/2^{2L}$}. 
In this specific case we obtain a ratio of approximately $4\%$. It means that we require 4\% numerical simulations in the correction chain to get an error with approximately the same order of magnitude of a full numerical MCMC chain, {if $C_{\rm num} \approx C_{\rm DL}$}.
In practice, obtaining $\epsilon$ may be challenging {and constants $C_{\rm num}$ and $C_{\rm DL}$ are generally unknown}, it may be more advantageous to compute $M_{\rm num}/M_{\rm DL}$ empirically. 
Indeed, this ratio is application and user dependent. 
For instance, if we can afford more numerical forward simulations, we could potentially run more of them, albeit at  higher computational costs. If we instead have a limited number of computational resources, and this is our primary constraint, we may want to limit the number of numerical forward simulations, while potentially accepting a higher error. 
This trade-off is depicted in Fig.~\ref{fig:elliptic uniform percentage num samples}, where we show the error with respect to the percentage of numerical samples against DL-based surrogate samples. 
We observe that the error reduces significantly with small number of numerical samples. 
However the error reduction by further increasing the ratio of numerical samples, even to 100\%, is not as significant.

{In our experiments, we ran very long (100,000 samples) base MCMC chains with DL-based surrogate models to ensure convergence. 
For completeness, we include additional diagnostic details of the MCMC experiments in ~\ref{sec:diagnostics}. 
First we show the trace plots and histogram plots of the MCMC chains with DL-based surrogate models in Fig~\ref{fig:elliptic_coef_uniform_MLchain_traceplot} and Fig~\ref{fig:elliptic_coef_uniform_MLchain_histogram}. 
The trace plots shows great mixing of samples from the posterior distribution. 
The histograms from different experiments show similar posterior distribution.
We also show a between chains comparison of sample mean with the eight independent experiments starting from random initial sample in Fig~\ref{fig:sample mean of elliptic uniform}. 
The plots show clear convergence of the sample mean. 
We also calculated the Effective Sample Size (ESS). The eight experiments show an average ESS of 30,842 out of 100,000 MCMC samples. 
The sample mean convergence between eight MCMC chain and the Potential Scale Reduction Factor (PSRF) is shown in Fig~\ref{fig:PSRF elliptic uniform}. 
The final PSRF, also known as R-hat, from Gelman-Rubin diagnostic is 1.00003. 
A value smaller than 1.2 is often considered as a good indication of convergence \cites{MR1665662,gelman1992inference}. 
We also show the trace plot, autocorrelation function and histogram of the QoI in the correction chain in Fig~\ref{fig:elliptic_uniform_qoia_traceplot_acf}, Fig~\ref{fig:qoia elliptic uniform hist}, Fig~\ref{fig:elliptic_uniform_qoib_traceplot_acf} and Fig~\ref{fig:qoib elliptic uniform hist}, where we observe a smaller variance compared with the base MCMC samples. 
Sample mean and PSRF plots are shown in Fig~\ref{fig:qoia sample mean and psrf uniform elliptic} and Fig~\ref{fig:qoib sample mean and psrf uniform elliptic}, where the PSRF value for both QoI are 1.00209 and 1.00238.
All the diagnostics show good indication of convergence. }

Finally, we estimated the overall computational time of the hybrid two-level MCMC approach compared to the numerical MCMC and DL MCMC. 
The estimation is based on the average runtime for one numerical sample and one DL-based surrogate sample. 
The FEM solver uses Intel Xeon E5-2620 CPU, while the DL-based surrogate model evaluations used one NVIDIA RTX A6000 GPU (we use the same CPU and GPU specs for all the subsequent experiments).
Due to the fact that the two MCMC chains in the hybrid MCMC method can be run concurrently, the hybrid MCMC with 100,000 surrogate samples and 5,000 numerical samples achieved the same speed-up as the plain MCMC with purely surrogate samples, but achieved a smaller error.

\subsubsection{Gaussian (log-normal) prior}
\label{sec:poisson-gauss}
\begin{figure}[H]
    \centering
    \subfloat[Sample 1.]{
        \includegraphics[width=0.25\textwidth]{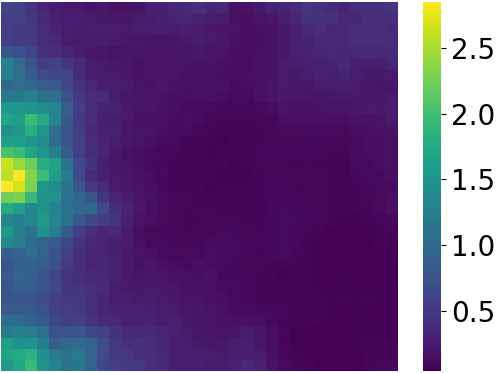}
        \label{fig:poisson-samples-sample1}
    }
    \subfloat[Sample 2.]{
        \includegraphics[width=0.25\textwidth]{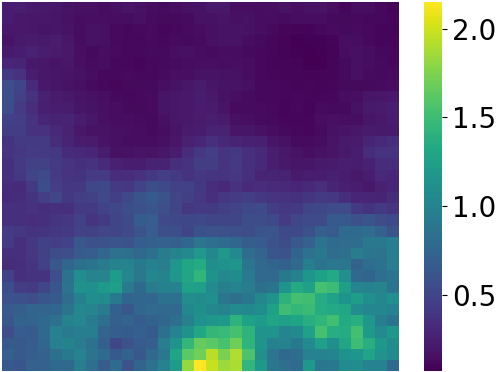}
        \label{fig:poisson-samples-sample2}
    }
    \subfloat[Sample 3.]{
        \includegraphics[width=0.25\textwidth]{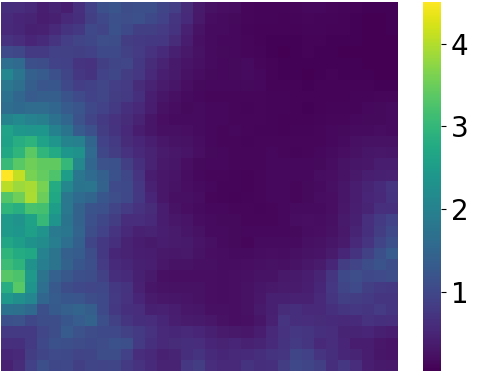}
        \label{fig:poisson-samples-sample3}
    }
    \subfloat[Sample 4.]{
        \includegraphics[width=0.25\textwidth]{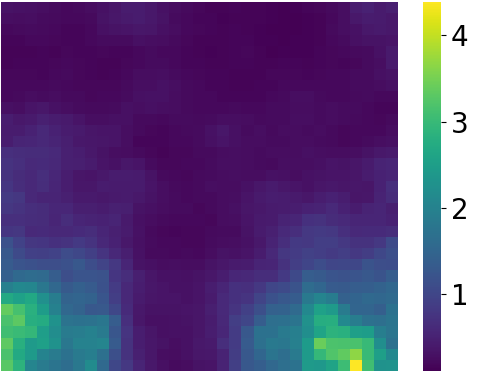}
        \label{fig:poisson-samples-sample4}
    }
    \caption{Samples obtained from the bi-Laplacian random field in equation~\eqref{eq:bi-laplacian} for the Poisson equation with Gaussian prior.}
    \label{fig:poisson-samples}
\end{figure}

\begin{figure}
    \centering
    \includegraphics[width=0.8\linewidth]{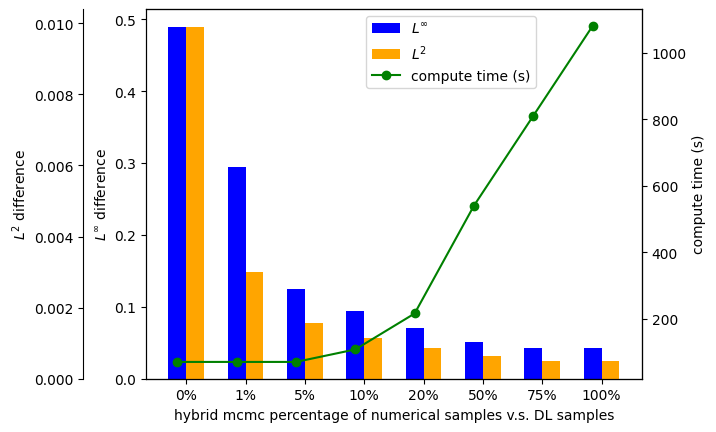}
    \caption{Error in comparison with numerical MCMC results at different percentage of numerical samples against DL-based surrogate samples for the elliptic experiment with Gaussian prior}
    \label{fig:elliptic gaussian percentage num sample}
\end{figure}

In contrast to the uniform prior case just shown, for this experiment, the QoI $K$ is spatially varying, that is: $K = K(x)$, and we assume its prior distribution to be sampled from the following bi-Laplacian Gaussian prior
\begin{equation}
\mathcal{A} m = \begin{cases}
\gamma \nabla \cdot(\Xi \nabla m)+\delta m & \text { in } D \\[0.4em]
(\Xi \nabla m) \cdot \boldsymbol{n}+\beta m & \text { on } \partial D,
\end{cases}
\label{eq:bi-laplacian}
\end{equation}
where $\boldsymbol{n}$ is the normal direction with respect to the boundary, $\gamma=0.1$, $\delta=0.5$, $\beta = \sqrt{\gamma \delta}$ and $\Xi$ controls the anisotropicity where we use an identity matrix. 
{In the numerical experiments, we sample $m$ to obtain the non-negative random field sample $K(x) = \exp(m)$ by solving the equation $\mathcal{A}m = f$ with FEM method on mesh~\ref{fig:mesh1}, where $f$ is sampled from white noise. This is equivalent to sample from the bi-Laplacian covariance operator $\mathcal{A}^{-2}$.
In this numerical experiment, the discretized $K(x)$ with finite dimension $32 \times 32 = 1024$ can be considered as the $z$ in $Q(z)$ in the preceding discussions.}
For more details of this particular prior setup, one can refer to Section 5.1.1 in~\cite{VillaPetraGhattas21} and references therein.
We show four examples of random samples generated in this experiment in Figure~\ref{fig:poisson-samples}.
Similarly to the uniform-prior case, we solve the forward problem in equation~\eqref{eq:poisson} via FEM on Mesh A (Figure~\ref{fig:mesh1}), and generate 8000 random samples.
Again, the 8000 samples are partitioned into 4000 training samples, 2000 validation samples, and 2000 test samples to train a convolutional neural network. 
The convolutional neural network consists of 3 encoding layer, 1 fully connected layer and 3 decoding layers, and it is trained using the Adam optimizer for 10000 epochs.
{We have a rough estimate of $\epsilon=0.228$.
With reference to Theorem~\ref{thm:final_error_estimate}{ assuming $C_{\rm num} \approx C_{\rm DL}$}, the correction chain needs around 1\% numerical samples compared to the long DL-based MCMC chain.
Again, this is not a rigorous estimate of the DL model error, {where the ratio is an underestimation}.} 
For a more accurate estimation of percentage of numerical samples needed, we include Fig~\ref{fig:elliptic gaussian percentage num sample} to show the error with respect to the percentage of numerical samples against DL-based surrogate samples. 
Fig~\ref{fig:elliptic gaussian percentage num sample} shows the error reduced significantly with 5\% or more numerical samples.
{We remark that the error used here is the $L^2$ and $L^{\infty}$ difference with respect to the classical numerical MCMC results, which is different from the error computed against a highly accurate quadrature in Fig~\ref{fig:elliptic uniform percentage num samples}. This is due to the computational challenge to compute quadrature of a high dimensional problem like this experiment. In the subsequent numerical experiment, we stick to the same approach.}

In analogy with the uniform-prior case, we run three experiments: (i) a plain MCMC chain with numerical solver, denoted as \textbf{Numerical MCMC} in Table~\ref{tab:elliptic_random_field_results}, (ii) a plain MCMC chain with DL-based surrogate model, denoted as \textbf{DL MCMC} in Table~\ref{tab:elliptic_random_field_results}, and the proposed hybrid approach, denoted as \textbf{Hybrid MCMC} in Table~\ref{tab:elliptic_random_field_results}. 
We show in detail three cases in Table~\ref{tab:elliptic_random_field_results}, namely $M_{\rm num} / M_{\rm DL} = 1\%, 5\%, 10\%$.
{Due to the high dimensional nature of the random field samples, we will not include all the trace plots and histogram of each pixel to show the diagnostics of each MCMC chain. 
We simply report that the max PSRF of the DL-based surrogate MCMC chain of 100,000 samples converged, being $1.00000714 < 1.2$.
The max PSRF for the QoI of A1, A2, ..., A8 {(refer to the hybrid method for Gaussian prior in~\ref{app:two-level-mcmc-gaussian})} in the correction chain of 1000 samples are 1.00794, 1.00744, 1.00558, 1.00413, 1.00231, 1.01205, 1.00962 and 1.00007. 
All PSRF values for the correction MCMC chain show good indication of convergence. 
}

\begin{figure}[H]
    \captionsetup[subfloat]{labelfont=scriptsize,textfont=scriptsize}
    \subfloat[Numerical-only \\ (100000).]{
        \includegraphics[width=0.19\textwidth]{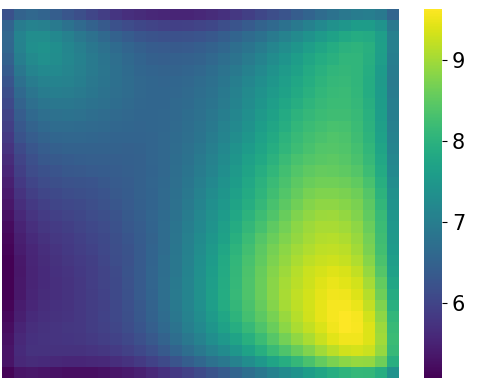}
        \label{fig:poisson-reference}
    }
    \subfloat[DL-only \\ (100000).]{
        \includegraphics[width=0.19\textwidth]{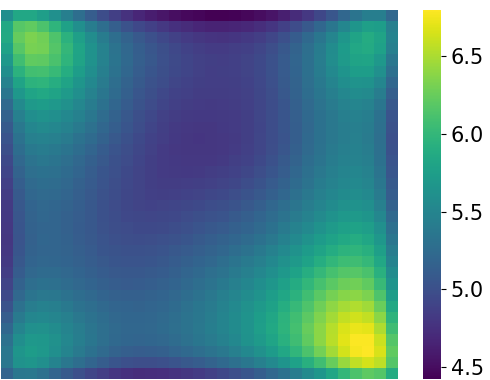}
        \label{fig:poisson-ml-posterior}
    }
    \subfloat[Hybrid \\ (100000+1000).]{
        \includegraphics[width=0.19\textwidth]{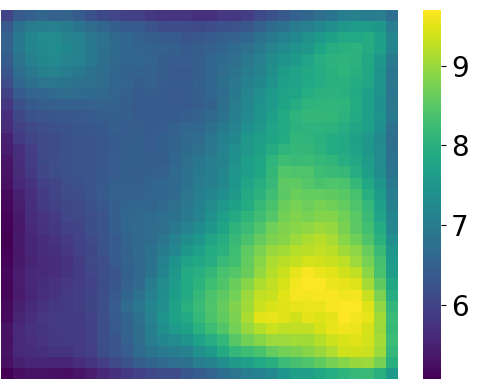}
        \label{fig:poisson-hybrid-posterior-1000}
    }
    \subfloat[Hybrid \\ (100000+5000).]{
        \includegraphics[width=0.19\textwidth]{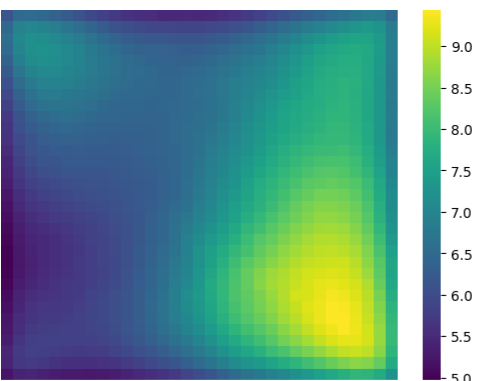}
        \label{fig:poisson-hybrid-posterior-5000}
    }
    \subfloat[Hybrid \\ (100000+10000).]{
        \includegraphics[width=0.19\textwidth]{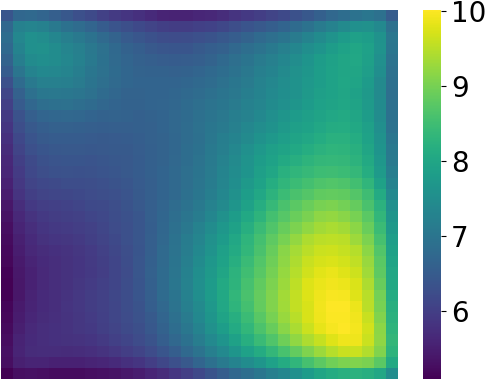}
        \label{fig:poisson-hybrid-posterior-10000}
    }
    \caption{Expected mean of $K(x)$ from eight runs of MCMC with elliptic equation with Gaussian prior.}
    \label{fig:poisson-expected-posterior}
\end{figure}

The average results of eight MCMC runs are depicted in Figure~\ref{fig:poisson-expected-posterior}. 
The expectation of the posterior from the MCMC chains generated solely with a DL-based surrogate model has obvious discrepancies with the plain MCMC chains generated solely with a numerical solver (that constitute the reference). 
Our hybrid two-level MCMC approach, with the addition of only few numerical samples, is able to significantly improve the results making it comparable to the reference, at a fraction of the computational cost.

We summarize the results in Table~\ref{tab:elliptic_random_field_results}, where the $L^2$ and $L^{\infty}$ difference between the DL-based surrogate accelerated MCMC results and classical numerical MCMC results are presented. 
From the results presented, there are significant improvements in terms of accuracy with the hybrid approach. 
With only $1\%$ additional numerical samples on top of the DL-based MCMC chain, the results get much closer to the one from the numerical MCMC.
With 5\% and 10\% additional numerical samples, the results of hybrid approach get even closer to the single chain numerical MCMC result. 
\begin{table}[H]
\scriptsize
\centering
\begin{tabular}{l|c|c|c|c|c}
\toprule
\textbf{Method} & 
  \begin{tabular}{@{}c@{}} \textbf{Numerical} \\ \textbf{MCMC} \end{tabular}   &   \begin{tabular}{@{}c@{}} \textbf{DL}        \\ \textbf{MCMC} \end{tabular} & \begin{tabular}{@{}c@{}} \textbf{Hybrid}    \\ \textbf{MCMC} \end{tabular}   &   \begin{tabular}{@{}c@{}} \textbf{Hybrid}    \\ \textbf{MCMC} \end{tabular}   &   \begin{tabular}{@{}c@{}} \textbf{Hybrid}    \\ \textbf{MCMC}   \end{tabular} \\ 
\midrule
\textbf{Samples} & 100,000 & 100,000 & 
\begin{tabular}{@{}c@{}}100,000 \\ + 1,000  \end{tabular} & 
\begin{tabular}{@{}c@{}}100,000 \\ + 5,000  \end{tabular} & 
\begin{tabular}{@{}c@{}}100,000 \\ + 10,000 \end{tabular} \\ 
\midrule
\textbf{$L^2$ difference} & 
- & 9.374E-3 & 2.997E-3 & 1.562E-3 & 1.114E-3 \\ 
\midrule
\textbf{$L^\infty$ difference} & 
- & 4.834E-1 & 2.944E-1 & 1.244E-1 & 9.407E-2 \\ 
\midrule
\begin{tabular}{@{}l@{}} \textbf{Compute} \\ \textbf{time [s]} \end{tabular} & 
1080.67 & 70.54 & \begin{tabular}{@{}l@{}} 81.35 (serial) \\ 70.54 (parallel) \end{tabular} & \begin{tabular}{@{}l@{}} 124.57 (serial) \\ 70.54 (parallel) \end{tabular} & \begin{tabular}{@{}l@{}} 178.61 (serial) \\ 108.07 (parallel) \end{tabular} \\ 
\midrule
\textbf{Speed up} & - & 15.32x & \begin{tabular}{@{}l@{}} 13.28x (serial) \\ 15.32x (parallel) \end{tabular} & \begin{tabular}{@{}l@{}} 8.67x (serial) \\ 15.32x (parallel) \end{tabular} & \begin{tabular}{@{}l@{}} 6.05x (serial) \\ 10x (parallel) \end{tabular} \\ 
\bottomrule
\end{tabular}
\caption{Difference between the posterior expectation results of plain numerical MCMC and the posterior expectation results of the plain and hybrid DL-based MCMC with elliptic equation with Gaussian prior}
\label{tab:elliptic_random_field_results}
\end{table}

\subsection{Nonlinear reaction-diffusion equation} 
\label{sec:reaction-diffusion}

We consider the Bayesian inverse problem with the forward model governed by the following nonlinear reaction-diffusion  equation in a two-dimensional unit square domain $D$
\begin{equation}
\label{eqn:reaction_diffusion}
\begin{cases}
\nabla \cdot (K(x) \nabla u(x)) + u^3 = 0, \\[0.5em]
u(x_1=0) = 0, \\[0.2em]
u(x_1=1) = 1, \\[0.5em]
\displaystyle \frac{\partial u}{\partial x_2} (x_2 = 0) = \frac{\partial u}{\partial x_2} (x_2 = 1) = 0.
\end{cases} 
\end{equation}  
Thirty-six equally distanced observations are captured from a random realization of the forward model with additional Gaussian noise $\delta$ with zero mean and variance $\sigma^2 = 0.1$.

\subsubsection{Uniform prior}
\label{sec:reaction-diffusion-uniform}

\begin{figure}
    \centering
    \includegraphics[width=0.8\linewidth]{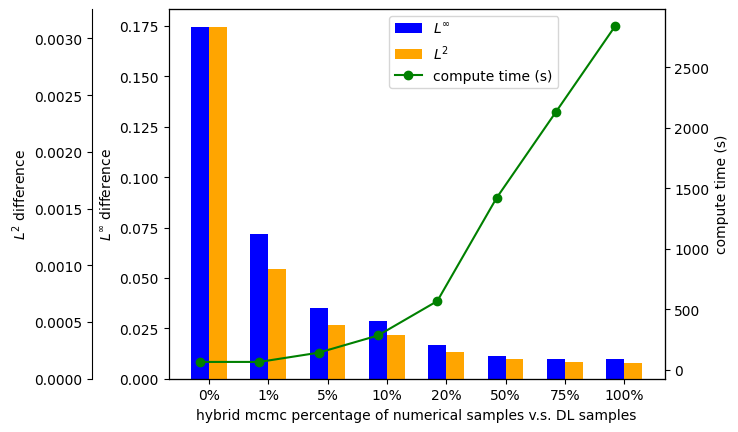}
    \caption{Error in comparison with numerical MCMC results at different percentage of numerical samples against DL-based surrogate samples}
    \label{fig:nonlinear_reaction_diffusion_uniform_sample_percentage}
\end{figure}

In this section, we consider the uniform prior case with a random field $K(x)$ that depends on uniformly distributed coefficients $z_i, i=0, 1, ..., 4$
\begin{align*}
    \ln(K(x)) = z_0 &+ z_1\cos(2\pi x_1)\sin(2\pi x_2) + z_2\sin(2\pi x_1)\cos(2\pi x_2)) \\
    &+ z_3\cos(2\pi x_1)\cos(2\pi x_2) + z_4\sin(2\pi x_1)\sin(2\pi x_2)),
\end{align*}
where $z_i \sim U[-1, 1], i=0,1, ..., 4$. 
We solve the reaction diffusion equation~\eqref{eqn:reaction_diffusion} with the FEM method on a $32 \times 32$ uniformly spaced Mesh A~\ref{fig:mesh1}. 
We randomly generated 4000 samples with a FEM  numerical solver. 
The 4000 data are partitioned into 2000 training samples, 1000 validation samples, and 1000 test samples. 
We train a vanilla Message Passing Graph Neural Network (MPGNN) as our DL surrogate. 
The details of the FEM solver and the Graph Neural Network architecture can be found in~\ref{section:reaction-diffusion_solver}. 
We have a rough estimate of $\epsilon=3.91$.
With reference to Theorem~\ref{thm:final_error_estimate} {assuming $C_{\rm num} \approx C_{\rm DL}$}, the correction chain only needs around 26\% numerical samples compared to the long DL-based MCMC chain. 
Yet again this is not a rigorous error rate for the DL model{, where the ratio is an overestimation.}
For a more accurate estimation of numerical samples needed, we include Fig~\ref{fig:nonlinear_reaction_diffusion_uniform_sample_percentage} to show the error with respect to the percentage of numerical samples against DL-based surrogate samples. 
Fig~\ref{fig:nonlinear_reaction_diffusion_uniform_sample_percentage} shows the error reduced significantly even with a small number of numerical samples.

\begin{figure}[H]
    \captionsetup[subfloat]{labelfont=scriptsize,textfont=scriptsize}
    \subfloat[Numerical-only \\ (100000)]{
        \includegraphics[width=0.19\textwidth]{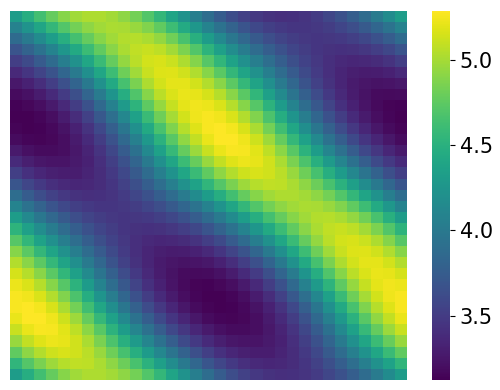}
        \label{fig:reaction-diffusion-reference-uniform}
    }
    \subfloat[{DL-only \\ (100000)}]{
        \includegraphics[width=0.19\textwidth]{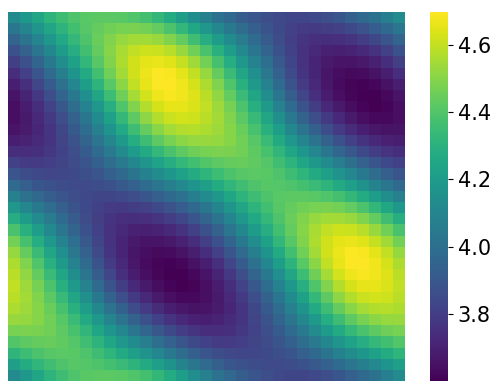}
        \label{fig:reaction-diffusion-ml-posterior-uniform}
    }
    \subfloat[Hybrid \\ (100000+1000)]{
        \includegraphics[width=0.19\textwidth]{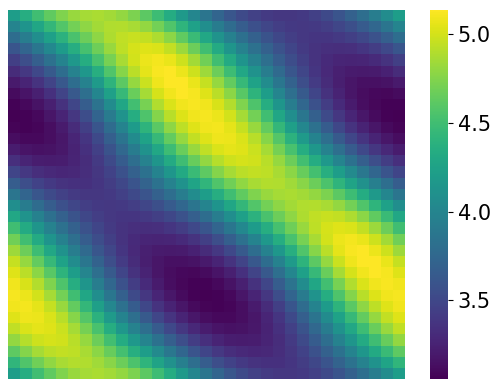}
        \label{fig:reaction-diffusion-hybrid-posterior-1000-uniform}
    }
    \subfloat[Hybrid \\ (100000+5000)]{
        \includegraphics[width=0.19\textwidth]{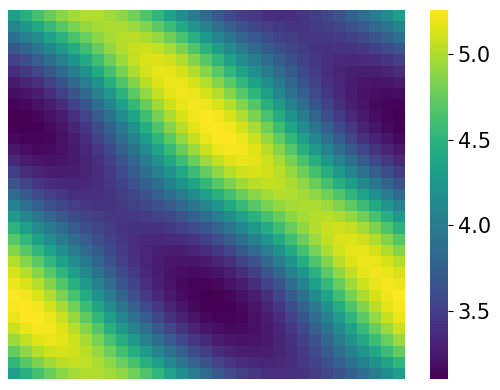}
        \label{fig:reaction-diffusion-hybrid-posterior-5000-uniform}
    }
    \subfloat[Hybrid \\ (100000+10000)]{
        \includegraphics[width=0.19\textwidth]{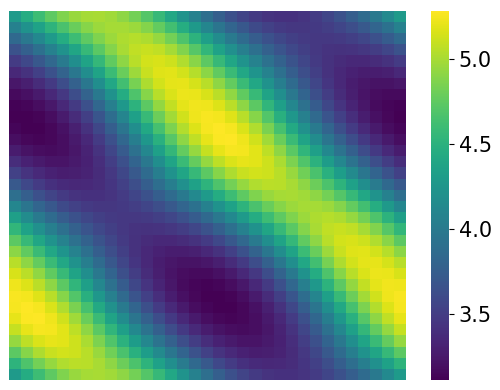}
        \label{fig:reaction-diffusion-hybrid-posterior-10000-uniform}
    }
    \caption{Expected mean of $K(x)$ from eight runs of MCMC with reaction diffusion equation with uniform prior.}
    \label{fig:reaction-diffusion-expected-posterior-uniform}
\end{figure}

We run three experiments: (i) a plain MCMC chain with numerical solver, denoted as \textbf{Numerical MCMC} in Table~\ref{tab:reaction_diffusion_random_field_results_uniform}, (ii) a plain MCMC chain with the MPGNN-based DL surrogate model, denoted as \textbf{DL MCMC} in Table~\ref{tab:reaction_diffusion_random_field_results_uniform}, and the proposed hybrid approach, denoted as \textbf{Hybrid MCMC} in Table~\ref{tab:reaction_diffusion_random_field_results_uniform}, {with the number of samples in the numerical chain being a fraction of the samples in the DL chain. We tested $M_{\rm num} / M_{\rm DL}$ ratios ranging from 1\% to 100\% as shown
in Fig~\ref{fig:nonlinear_reaction_diffusion_uniform_sample_percentage}. 
We show in details three cases ($1\%$, 5\%, and 10\%) in Fig~\ref{fig:reaction-diffusion-expected-posterior-uniform} and Table~\ref{tab:reaction_diffusion_random_field_results_uniform}.
For compactness of the paper we skip trace plots and histogram plots.
We computed PSRF for all QoIs in both MCMC chains. 
The results show maximum PSRF value of 1.00922 which is smaller than the common indicative 1.2 value.
}

The average results of five MCMC runs are presented in Figure~\ref{fig:reaction-diffusion-expected-posterior-uniform}. 
The $L^2$ and $L^{\infty}$ difference between the numerical MCMC results and the DL-based methods (including the DL MCMC and the Hybrid MCMC) are presented in Table~\ref{tab:reaction_diffusion_random_field_results_uniform}. 
These show that the additional numerical samples in the hybrid method reduce both $L^2$ and $L^{\infty}$ difference compared to the DL MCMC method.
The results validate the theoretical conclusion in Theorem~\ref{thm:final_error_estimate}, that the hybrid two-level approach can reach the same accuracy level as the numerical MCMC at a fraction of the computational cost. 
%
\begin{table}[H]
\scriptsize
\centering
\begin{tabular}{l|c|c|c|c|c}
\toprule
\textbf{Method} & 
\begin{tabular}{@{}c@{}} \textbf{Numerical} \\ \textbf{MCMC} \end{tabular} & \begin{tabular}{@{}c@{}} \textbf{DL}        \\ \textbf{MCMC} \end{tabular} & \begin{tabular}{@{}c@{}} \textbf{Hybrid}    \\ \textbf{MCMC} \end{tabular} & \begin{tabular}{@{}c@{}} \textbf{Hybrid}    \\ \textbf{MCMC} \end{tabular} & \begin{tabular}{@{}c@{}} \textbf{Hybrid}    \\ \textbf{MCMC} \end{tabular} \\ 
\midrule
\textbf{Samples} & 100,000 & 100,000 & 
\begin{tabular}{@{}c@{}}100,000 \\ + 1,000  \end{tabular} & 
\begin{tabular}{@{}c@{}}100,000 \\ + 5,000  \end{tabular} & 
\begin{tabular}{@{}c@{}}100,000 \\ + 10,000 \end{tabular} \\ 
\midrule
\textbf{$L^2$ difference} &  - & 2.914E-3 & 9.638E-4 & 4.756E-4 & 3.861E-4 \\
\midrule
\textbf{$L^\infty$ difference} & 
- & 1.772E-1 & 0.716E-2 & 0.351E-2 & 0.285E-2 \\
\midrule
\begin{tabular}{@{}l@{}} \textbf{Compute} \\ \textbf{time [s]} \end{tabular} & 2840.75 & 65.78 & \begin{tabular}{@{}l@{}} 94.19 (serial) \\ 65.78 (parallel) \end{tabular} & \begin{tabular}{@{}l@{}} 207.82 (serial) \\ 142.04 (parallel) \end{tabular} & \begin{tabular}{@{}l@{}} 349.86 (serial) \\ 284.08 (parallel) \end{tabular} \\
\midrule
\textbf{Speed up} & - & 43.19x & \begin{tabular}{@{}l@{}} 30.15x (serial) \\ 43.19x (parallel) \end{tabular} & \begin{tabular}{@{}l@{}} 13.67x (serial) \\ 20x (parallel) \end{tabular} & \begin{tabular}{@{}l@{}} 8.12x (serial) \\ 10x (parallel) \end{tabular} \\ 
\bottomrule
\end{tabular}
\caption{Difference between the posterior expectation results of plain numerical MCMC and the posterior expectation results of the plain and hybrid DL-based MCMC with non-linear reaction-diffusion equation with uniform prior}
\label{tab:reaction_diffusion_random_field_results_uniform}
\end{table}
\color{black}

\subsubsection{Gaussian prior}
\label{sec:reaction-diffusion-gauss}
\begin{figure}
    \centering
    \includegraphics[width=0.8\linewidth]{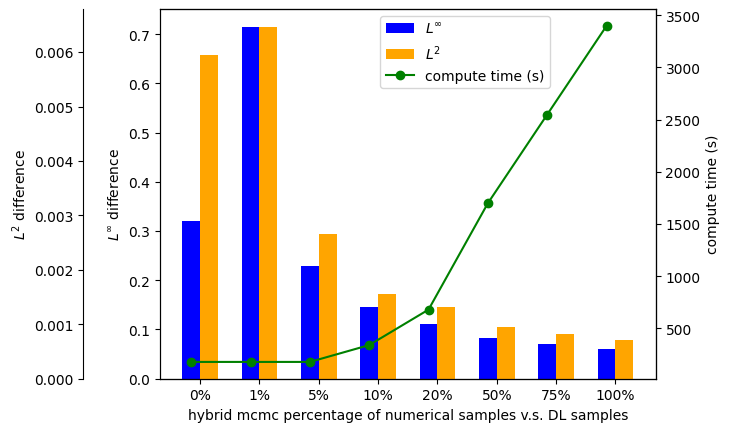}
    \caption{Error in comparison with numerical MCMC results at different percentage of
numerical samples against DL-based surrogate samples for the reaction diffusion experiments with Gaussian prior}
    \label{fig:error_over_different_percentage_reaction_diffusion_gaussian}
\end{figure}

{We follow the same random field setup for $K(x)$ in Section~\ref{sec:poisson-gauss}.}
We randomly generated 4000 samples with the Finite Element solver. 
The 4000 data are partitioned into 2000 training data, 1000 validation data, and 1000 test data to train a U-net neural network. 
The details of the FEM solver and the neural network architecture can be found in~\ref{section:reaction-diffusion_solver}.
We have a rough estimate of $\epsilon=0.93$.
With reference to Theorem~\ref{thm:final_error_estimate} {assuming $C_{\rm num} \approx C_{\rm DL}$}, the correction chain only needs around 1\% numerical samples compared to the long DL-based MCMC chain. 
For a more accurate estimation, we include Fig~\ref{fig:error_over_different_percentage_reaction_diffusion_gaussian} to show the error with respect to the percentage of numerical samples.
We see that the estimation $\epsilon$ we have is not perfect, where in the sweep test we see an increase of the error with 1\% of numerical samples but better results are acheived with 5\% or more numerical samples.
%
We performed numerical experiments on three setups: (i) a plain MCMC chain with numerical solver, denoted as \textbf{Numerical MCMC} in Table~\ref{tab:reaction_diffusion_random_field_results}, (ii) a plain MCMC chain with DL-based surrogate model, denoted as \textbf{DL MCMC} in Table~\ref{tab:reaction_diffusion_random_field_results}, (iii) and the proposed hybrid approach, denoted as \textbf{Hybrid MCMC} in Table~\ref{tab:reaction_diffusion_random_field_results}, {once again testing cases such that the $M_{\rm num} / M_{\rm DL}$ ratio ranges from 1\% to 100\%, as shown in Fig~\ref{fig:error_over_different_percentage_reaction_diffusion_gaussian}. 
We show in detail three cases (1\%, 5\%, and 10\%) in Table~\ref{tab:reaction_diffusion_random_field_results}. 
For compactness of the paper we skip trace plots and histogram plots. We computed PSRF for all QoIs in both MCMC chains. 
The results show maximum PSRF value of 1.0131 which is smaller than the common indicative 1.2 value.}

\begin{figure}[!htb]%
    \captionsetup[subfloat]{labelfont=scriptsize,textfont=scriptsize}
    \subfloat[\centering Numerical MCMC(100,000)]{
        \includegraphics[width=0.19\textwidth]{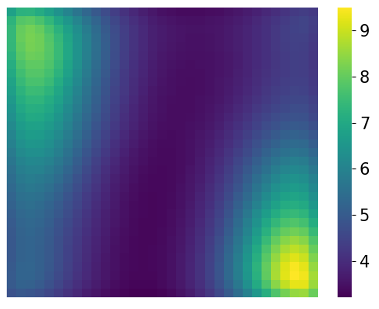}
        \label{fig:reaction-diffusion-posterior-num}
    }
    \subfloat[\centering DL MCMC (100,000)]{
        \includegraphics[width=0.19\textwidth]{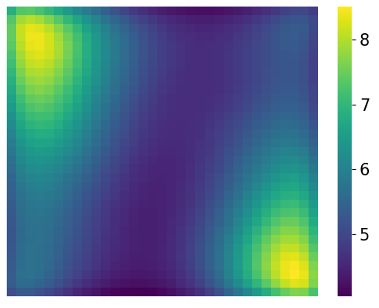}
        \label{fig:reaction-diffusion-posterior-dl}
    }
    \subfloat[\centering Hybrid (100,000+1,000)]{
        \includegraphics[width=0.19\textwidth]{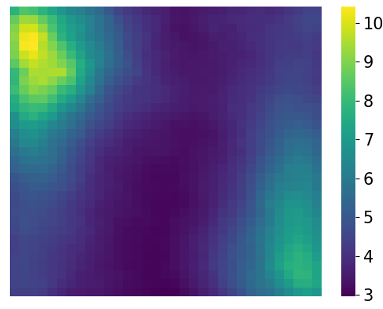}
        \label{fig:reaction-diffusion-posterior-1000}
    }
    \subfloat[\centering Hybrid (100,000+5,000)]{
        \includegraphics[width=0.19\textwidth]{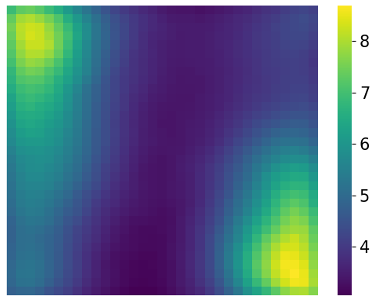}
        \label{fig:reaction-diffusion-posterior-5000}
    }
    \subfloat[\centering Hybrid (100,000+10,000)]{
        \includegraphics[width=0.19\textwidth]{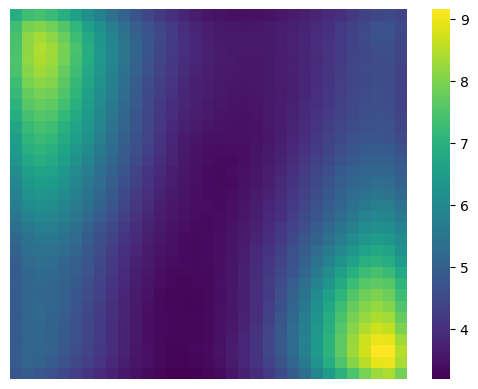}
        \label{fig:reaction-diffusion-posterior-10000}
    }
    \caption{Expected mean of $K(x)$ from eight runs of MCMC for the nonlinear reaction-diffusion experiment with Gaussian prior}
    \label{fig:reaction-diffusion-posterior}
\end{figure}

The results of the average of the eight MCMC run are shown in Figure~\ref{fig:reaction-diffusion-posterior}. 
The $L^2$ and $L^\infty$ difference between MCMC results generated with the DL-based surrogate and plain numerical MCMC are included in Table \ref{tab:reaction_diffusion_random_field_results}. 
The results show that the additional numerical samples in the hybrid method reduce both $L^2$ and $L^\infty$ difference compared to the plain MCMC model. 

\begin{table}[H]
\scriptsize
\centering
\begin{tabular}{l|c|c|c|c|c}
\toprule
\textbf{Method} & 
\begin{tabular}{@{}c@{}} \textbf{Numerical} \\ \textbf{MCMC} \end{tabular} & \begin{tabular}{@{}c@{}} \textbf{DL}        \\ \textbf{MCMC} \end{tabular} & \begin{tabular}{@{}c@{}} \textbf{Hybrid}    \\ \textbf{MCMC} \end{tabular} & \begin{tabular}{@{}c@{}} \textbf{Hybrid}    \\ \textbf{MCMC} \end{tabular} & \begin{tabular}{@{}c@{}} \textbf{Hybrid}    \\ \textbf{MCMC} \end{tabular} \\ 
\midrule
\textbf{Samples} & 100,000 & 100,000 & 
\begin{tabular}{@{}c@{}}100,000 \\ + 1,000  \end{tabular} & 
\begin{tabular}{@{}c@{}}100,000 \\ + 5,000  \end{tabular} & 
\begin{tabular}{@{}c@{}}100,000 \\ + 10,000 \end{tabular} \\ 
\midrule
\textbf{$L^2$ difference} &  - & 5.953E-3 & 6.466E-3 & 2.643E-3 & 1.565E-3 \\
\midrule
\textbf{$L^\infty$ difference} & 
- & 3.203E-1 & 7.151E-1 & 2.296E-2 & 1.448E-2 \\
\midrule
\begin{tabular}{@{}l@{}} \textbf{Compute} \\ \textbf{time [s]} \end{tabular} & 3396.72 & 178.81 & \begin{tabular}{@{}l@{}} 212.78 (serial) \\ 178.81 (parallel) \end{tabular} & \begin{tabular}{@{}l@{}} 348.65 (serial) \\ 178.81 (parallel) \end{tabular} & \begin{tabular}{@{}l@{}} 518.78 (serial) \\ 339.67 (parallel) \end{tabular} \\
\midrule
\textbf{Speed up} & - & 18.99x & \begin{tabular}{@{}l@{}} 15.96x (serial) \\ 18.99x (parallel) \end{tabular} & \begin{tabular}{@{}l@{}} 9.74x (serial) \\ 18.99x (parallel) \end{tabular} & \begin{tabular}{@{}l@{}} 6.55x (serial) \\ 10.0x (parallel) \end{tabular} \\ 
\bottomrule
\end{tabular}
\caption{Difference between the posterior expectation results of plain numerical MCMC and the posterior expectation results of the plain and hybrid DL-based MCMC with non-linear reaction-diffusion equation with Gaussian prior}
\label{tab:reaction_diffusion_random_field_results}
\end{table}

\subsection{Navier Stokes equations}
\label{sec:navier-stokes}

We consider the Bayesian inverse problem with the forward model governed by the two dimensional Navier Stokes equations in the vorticity form in a domain of two-dimensional unit torus $\mathbb{T}^2$, 
\begin{equation}
\label{eq:ns}
\begin{cases}
\displaystyle \frac{\partial \omega(x, t)}{\partial t} + u(x, t) \cdot \nabla \omega(x, t) - \nu \Delta \omega(x, t)  = f, &\quad \text{for } x\in \mathbb{T}^2, \\[0.5em]
\nabla \cdot u(x, t) = 0, &\quad \text{for } x\in\mathbb{T}^2, \\[0.5em]
\omega(x,0) = \omega_0;
\end{cases}
\end{equation}
with periodic boundary conditions and forcing
\begin{equation}
\label{eqn:forcing1}
    f = 0.1 (\sin(2 \pi (x_1+x_2))+\cos(2 \pi (x_1+x_2))). 
\end{equation}
$\omega$ is the vorticity, $u$ is the velocity, and $\nu = 0.001$ is the viscosity. 
Thirty-six equally distanced observations are captured from a random realization of the forward model with additional Gaussian noise $\delta$ with zero mean and variance $\sigma^2 = 1$.

\subsubsection{Uniform prior}
\label{sec:navier-stokes-uniform}
\begin{figure}
    \centering
    \includegraphics[width=0.8\linewidth]{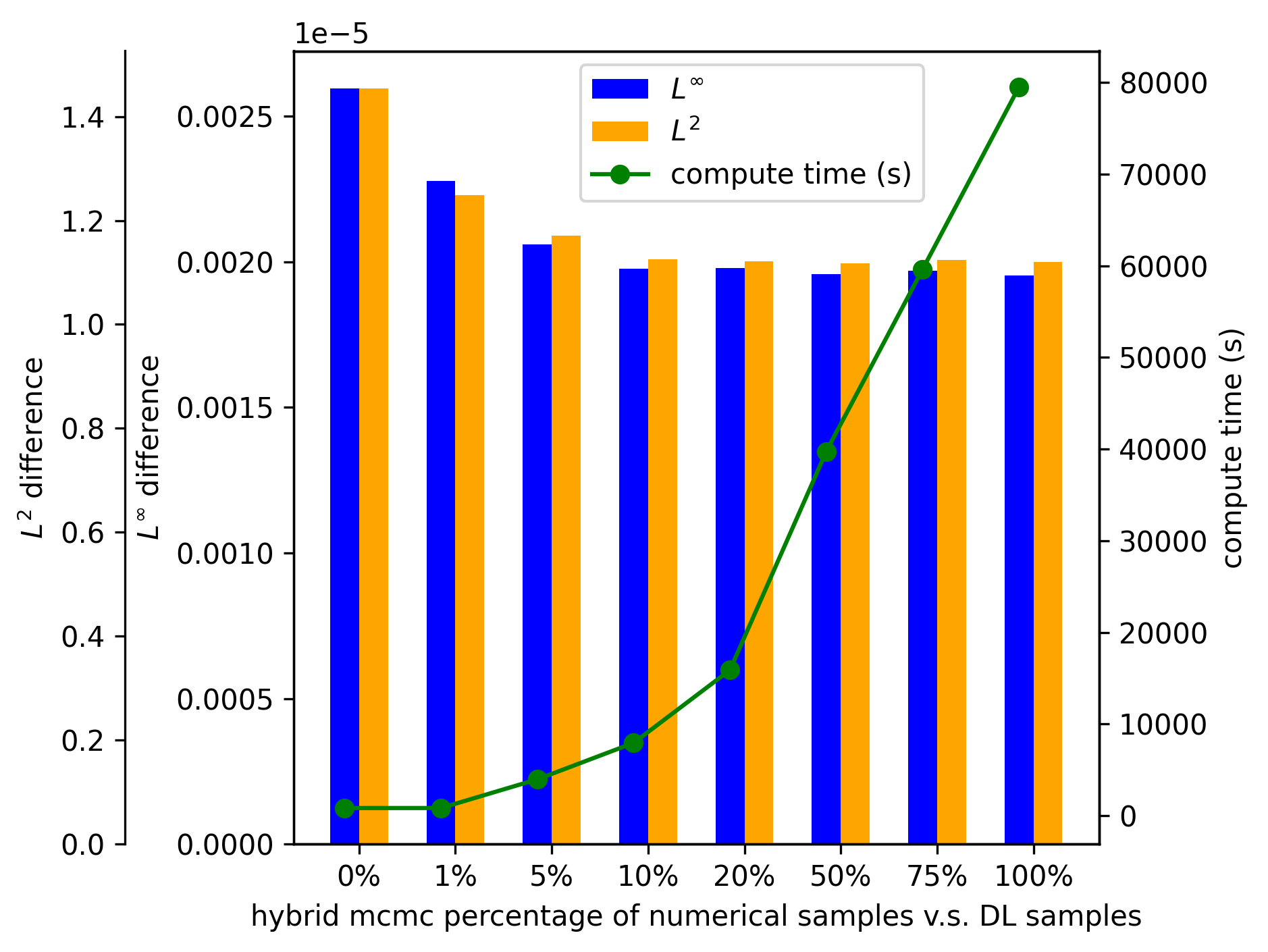}
    \caption{Error in comparison with numerical MCMC results at different percentage of
numerical samples against DL-based surrogate samples for the 2D Navier Stokes experiments with uniform prior}
\label{fig:error_over_different_percentage_ns_uniform}
\end{figure}

In this section, we consider the uniform prior case, using a Fourier expansion of $\omega_0$ with coefficient of each Fourier term uniformly distributed
\begin{equation*}
    Z_{mn} = z_{mn} N^2 \sqrt{2}\cdot 7^{3/2}(4\pi^2 (m^2 + n^2) + 49)^{-5/4},
\end{equation*}
where $N=64$ is the number of modes in each axis, m and n are the mode index in the $x$ and $y$ directions, and $z_{mn}\sim U[-1, 1]$. 
Using the inverse fast Fourier transform (ifft), we have 
\begin{equation*}
    \omega_0(x, y) = \rm{ifft}(Z) = \frac{1}{N^2} \sum_{m=0}^{N-1}\sum_{n=0}^{N-1} Z_{mn}\exp(i2\pi m x/N + i2\pi ny/N).
\end{equation*}
{This particular setup is chosen to match the Gaussian prior $\mathcal{N}(0, 7^{\frac{3}{2}}(-\Delta+49 I)^{-2.5})$ used later in the Gaussian case. 
In the Gaussian prior case, The coefficient $z_{mn}$ follows $N(0, 1)$ instead of an uniform prior.}
We solve the two-dimensional Navier-Stokes equations~\eqref{eq:ns} by using a pseudospectral method with Crank-Nicolson time integration on Mesh B, that is composed of $64 \times 64$ collocation points. 
We randomly generated 8000 samples with the numerical solver, where 4000 samples are  used for training, 2000 for validation, and 2000 for testing. 
DeepONet is used as the DL-based surrogate model.
For more details on the numerical solver and the DL architecture, the interested reader can refer to~\ref{section:ns_solver}.  
{We have a rough estimate of $\epsilon=1.93$.
With reference to Theorem~\ref{thm:final_error_estimate} {assuming $C_{\rm num} \approx C_{\rm DL}$}, the correction chain needs less than 1\% numerical samples compared to the long DL-based MCMC chain.
}
We also include Fig~\ref{fig:error_over_different_percentage_ns_uniform} to show the error with respect to the percentage of numerical samples.
In the sweep test, we see an increase of the error with 1\% of numerical samples but slightly better results are achieved with 5\% or more numerical samples.

%

We run three experiments: (i) a plain MCMC chain with numerical solver, denoted as \textbf{Numerical MCMC} in Table~\ref{tab:nse-uniform-prior}, (ii) a plain MCMC chain with DL-based surrogate model, denoted as \textbf{DL MCMC} in Table~\ref{tab:nse-uniform-prior}, and the proposed hybrid approach, denoted as \textbf{Hybrid MCMC} in Table~\ref{tab:nse-uniform-prior}. 
As usual, we tested rations $M_{\rm num} / M_{\rm DL}$ from 1\% to 100\%, as shown in Fig.~\ref{fig:error_over_different_percentage_ns_uniform}. We show in details three cases (1\%, 5\%, and 10\%) in Table~\ref{tab:nse-uniform-prior}. 
For compactness of the paper we skip trace plots and histogram plots. We computed PSRF for all QoIs in both MCMC chains. 
The results show maximum PSRF value of 1.029 which is smaller than the common indicative 1.2 value.
The average results of eight MCMC runs are included in Table~\ref{tab:nse-uniform-prior}. 
We note that in this case the DL MCMC based on DeepONet already achieves relatively good performance, and it is obviously the fastest method. 
Only marginal accuracy gains are achieved by using our hybrid approach with more than 5\% of numerical samples. 
This can be the case when the DL-surrogate model approximates the numerical forward operator well, such that adding some numerical samples within our hybrid MCMC framework does not significantly improve the results.
\begin{table}[H]
\scriptsize
\centering
\begin{tabular}{l|c|c|c|c|c}
\toprule
\textbf{Method} & 
\begin{tabular}{@{}c@{}} \textbf{Numerical} \\ \textbf{MCMC} \end{tabular} & \begin{tabular}{@{}c@{}} \textbf{DL}        \\ \textbf{MCMC} \end{tabular}  & \begin{tabular}{@{}c@{}} \textbf{Hybrid}    \\ \textbf{MCMC} \end{tabular} & \begin{tabular}{@{}c@{}} \textbf{Hybrid}    \\ \textbf{MCMC} \end{tabular}  & \begin{tabular}{@{}c@{}} \textbf{Hybrid}    \\ \textbf{MCMC} \end{tabular} \\
\midrule
\textbf{Samples} & 100,000 & 100,000  & \begin{tabular}{@{}c@{}}100,000 \\ + 1,000  \end{tabular} & 
\begin{tabular}{@{}c@{}}100,000 \\ + 5,000  \end{tabular} & 
\begin{tabular}{@{}c@{}}100,000 \\ + 10,000 \end{tabular} \\ 
\midrule
$L^2$ difference &  - & 1.453E-05 & 1.249E-05 &  1.171E-05 & 1.125E-5  \\
\midrule
$L^\infty$ difference & 
- & 2.594E-3 & 2.278E-3 &  2.2.061E-3 & 1.976E-3 \\
\midrule
\begin{tabular}{@{}l@{}} \textbf{Compute} \\ \textbf{time [s]} \end{tabular} & 
79490.13 & 846.78 & \begin{tabular}{@{}l@{}} 1641.68 (serial) \\  846.78 (parallel) \end{tabular} & \begin{tabular}{@{}l@{}} 4821.29 (serial) \\ 3974.51 (parallel) \end{tabular} & \begin{tabular}{@{}l@{}} 8795.79 (serial) \\ 7949.01 (parallel) \end{tabular} \\ 
\midrule
\textbf{Speed up} & - & 93.87x & \begin{tabular}{@{}l@{}} 48.42x (serial) \\ 93.87x (parallel) \end{tabular} & \begin{tabular}{@{}l@{}} 16.48x (serial) \\ 20x (parallel) \end{tabular} & \begin{tabular}{@{}l@{}} 9.37x (serial) \\ 10x (parallel) \end{tabular} \\ 
\bottomrule
\end{tabular}
\caption{Difference between the posterior expectation results of plain numerical MCMC and the posterior expectation results of the plain and hybrid DL-based MCMC with Navier Stokes equations with uniform prior}
\label{tab:nse-uniform-prior}
\end{table}

\subsubsection{Gaussian prior}
\label{sec:navier-stokes-gauss}

In this section, we consider the Gaussian prior case, focusing on $\omega_0$. 
We sample the Gaussian random field from the following Gaussian prior with the distribution: $\mathcal{N}(0, 7^{\frac{3}{2}}(-\Delta+49 I)^{-2.5})$.
We solve the above equation with the pseudo-spectral method with the Crank-Nicolson time integration method. In this experiment a resolution of $64 \times 64$ is used. We randomly generated 4000 samples with the numerical solver. The 4000 data are partitioned into 2000 training data, 1000 validation data, and 1000 test data to train the Fourier Neural Operator (FNO) as described in \cite{DBLP:conf/iclr/LiKALBSA21}. Details of the DL-based surrogate model setup can be found in~\ref{section:ns_solver}.
We have a rough estimate of $\epsilon=3.65$.
With reference to Theorem~\ref{thm:final_error_estimate} {assuming $C_{\rm num} \approx C_{\rm DL}$}, the correction chain needs only around 5\% numerical samples if compared to the long DL-based MCMC chain.

%
\begin{figure}
    \centering
    \includegraphics[width=0.8\linewidth]{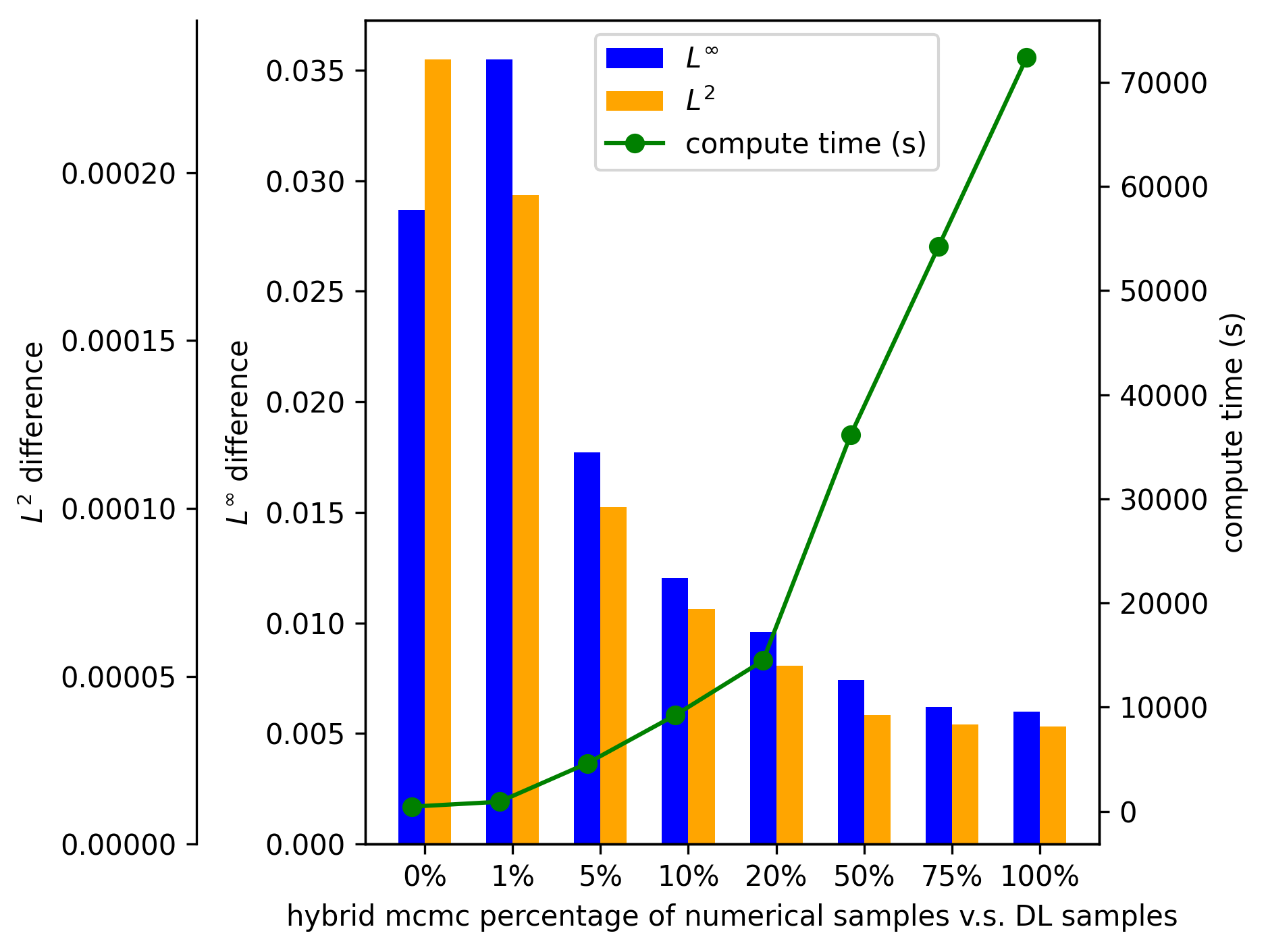}
    \caption{Error in comparison with numerical MCMC results at different percentage of
numerical samples against DL-based surrogate samples for the 2D Navier Stokes experiments with Gaussian prior}
    \label{fig:error_over_different_percentage_ns_gaussian}
\end{figure}
We run three experiments: (i) a plain MCMC chain with numerical solver, denoted as \textbf{Numerical MCMC} in Table~\ref{tab:reaction_nse_field_results}, (ii) a plain MCMC chain with DL-based surrogate model, denoted as \textbf{DL MCMC} in Table~\ref{tab:reaction_nse_field_results}, and the proposed hybrid approach, denoted as \textbf{Hybrid MCMC} in Table~\ref{tab:reaction_nse_field_results}, {where we tested cases with the ratio $M_{\rm num} / M_{\rm DL}$ ranging from 1\% to 100\%, as shown in Fig \ref{fig:error_over_different_percentage_ns_gaussian}. We show in details three cases (1\%, 5\%, and 10\%) in Table~\ref{tab:reaction_nse_field_results}. 
For compactness of the paper we skip trace plots and histogram plots. We computed PSRF for all QoIs in both MCMC chains. The results show maximum PSRF value of 1.089 which is smaller than the common indicative 1.2 value.}

The results of the average of the eight MCMC run are presented in Table \ref{tab:reaction_nse_field_results} where the posterior expectation obtained from the hybrid MCMC algorithm is closer to the posterior expectation obtained from plain numerical MCMC, but with much lower computational cost.
\begin{table}[H]
\scriptsize
\centering
\begin{tabular}{l|c|c|c|c|c}
\toprule
\textbf{Method} & 
\begin{tabular}{@{}c@{}} \textbf{Numerical} \\ \textbf{MCMC} \end{tabular} & \begin{tabular}{@{}c@{}} \textbf{DL}        \\ \textbf{MCMC} \end{tabular} & \begin{tabular}{@{}c@{}} \textbf{Hybrid}    \\ \textbf{MCMC} \end{tabular} & \begin{tabular}{@{}c@{}} \textbf{Hybrid}    \\ \textbf{MCMC} \end{tabular} & \begin{tabular}{@{}c@{}} \textbf{Hybrid}    \\ \textbf{MCMC} \end{tabular} \\ 
\midrule
\textbf{Samples} & 100,000 & 100,000 & 
\begin{tabular}{@{}c@{}}100,000 \\ + 1,000  \end{tabular} & 
\begin{tabular}{@{}c@{}}100,000 \\ + 5,000  \end{tabular} & 
\begin{tabular}{@{}c@{}}100,000 \\ + 10,000 \end{tabular} \\ 
\midrule
$L^2$ difference & 
N.A. & 2.337E-4 & 1.934E-4 & 1.004E-4 & 6.995-5 \\ 
\midrule
$L^\infty$ difference & 
N.A. & 2.869E-2 & 3.548E-2 & 1.772E-2 & 1.205E-2 \\
\midrule
\begin{tabular}{@{}l@{}} \textbf{Compute} \\ \textbf{time [s]} \end{tabular} & 
92343.7 & 456.4 & \begin{tabular}{@{}l@{}} 1379.8 (serial) \\ 923.4 (parallel) \end{tabular} & \begin{tabular}{@{}l@{}} 5082.6 (serial) \\ 4617.2 (parallel) \end{tabular} & \begin{tabular}{@{}l@{}} 9690.8 (serial) \\ 9234.4 (parallel) \end{tabular} \\ 
\midrule
\textbf{Speed up} & - & 202.29x & \begin{tabular}{@{}l@{}} 66.96x (serial) \\ 100x (parallel) \end{tabular} & \begin{tabular}{@{}l@{}} 18.17x (serial) \\ 20x (parallel) \end{tabular} & \begin{tabular}{@{}l@{}} 9.53x (serial) \\ 10x (parallel) \end{tabular} \\ 
\bottomrule
\end{tabular}
    \caption{Difference between the posterior expectation results of plain numerical MCMC and the posterior expectation results of the plain and hybrid DL-based MCMC with 2D Navier-Stokes equation}
    \label{tab:reaction_nse_field_results}
\end{table}
\section{Conclusion}
\label{sec:conclusions}
In this paper, we introduced a novel method, that we named hybrid two-level MCMC approach, to compute the posterior mean of quantities of interest in Bayesian inverse problems. 
In this method, we take advantage of the fast evaluation of DL surrogates and of the high accuracy of numerical models. 
In particular, we have theoretically shown the potential to solve Bayesian inverse problems accurately up to an estimator error $\mathcal{O}(h)$, by coupling one short MCMC chain generated by a high-fidelity numerical solver with mesh size $h$ and another long MCMC chain generated with fast DL surrogates.
We show the complete estimator error analysis and conclude that its theoretical bound is $\mathcal{O}(2^{-L})$ with one long base MCMC chain of $\mathcal{O}(2^{2L})$ number of DL surrogate samples and a short correction MCMC chain of $\mathcal{O}((1+2^\epsilon)^2))$ number of numerical samples, given the numerical forward model has an error rate of $2^{-L}$ and DL-based forward surrogate has an error rate of $2^{-L+\epsilon}$.
In addition, we show that with a surrogate speedup rate $s$ of one forward solve, the overall speedup of our hybrid algorithm is {$\mathcal{O}( 2^{2L}/\max(\frac{1}{s} 2^{2L}, \frac{C_{\rm num}}{C_{\rm DL}}(1+2^{\epsilon})^2))$}.
This implies that the overall speedup depends on the performance of the surrogate model, more speedup can be expected with high accurate surrogate models due to less numerical samples needed for error correction.
To validate the theoretical findings, we performed numerical experiments on a Poisson equation, a nonlinear reaction-diffusion equation, a the Navier-Stokes equation. 
In all numerical experiments, we use both uniform priors and Gaussian priors. 
All results of our numerical experiments {qualitatively validate our theoretical findings.
However, we note that the theoretical result depends on the assumption of knowing the exact DL surrogate error, and the exact constants that appeared in the derivation of} Theorem~\ref{thm:final_error_estimate}, which are rather challenging to obtain analytically.
{Due to this, we see in the numerical experiments that the actual optimal sample ratio from numerical sweeping test may deviate from the theoretical estimates of $(1+2^\epsilon)^2/2^{2L}$ without taking into consideration of the constants $C_{\rm num}$ and $C_{\rm DL}$.}
Yet, with a sweeping test, we show that for almost every experimental case, we can improve the accuracy of the MCMC estimator with an increasing number of numerical samples under our hybrid two-level MCMC method.
We note that the theoretical framework proposed can also be used to understand the feasibility of using DL surrogates only (without hybridizing them with high-fidelity numerical solvers). 
More specifically, if the theoretical error estimates are already within the numerical error of a DL-only surrogate, then the use of the high-fidelity numerical solver will unlikely yield better accuracy. 
This theoretical result is particularly important, given the widespread use of DL surrogates in the field.
As a final note, this paper focuses on a hybrid approach for the MCMC method to compute the posterior mean of quantities of interest in Bayesian inverse problems governed by PDEs; however, the same approach can in principle be applied to other Bayesian inverse problems not necessarily governed by PDEs, e.g. ODE governed Bayesian inverse problems, as well other methods such as filtering algorithms like the ensemble Kalman filter or sequential Monte Carlo methods.

\newpage
\color{black}
\bibliographystyle{amsrefs} 
\bibliography{reference}

\begin{bibdiv}
\begin{biblist}

\bib{pmlr-v97-allen-zhu19a}{inproceedings}{
      author={Allen-Zhu, Zeyuan},
      author={Li, Yuanzhi},
      author={Song, Zhao},
       title={A convergence theory for deep learning via over-parameterization},
        date={201909--15 Jun},
   booktitle={Proceedings of the 36th international conference on machine learning},
      editor={Chaudhuri, Kamalika},
      editor={Salakhutdinov, Ruslan},
      series={Proceedings of Machine Learning Research},
      volume={97},
   publisher={PMLR},
       pages={242\ndash 252},
         url={https://proceedings.mlr.press/v97/allen-zhu19a.html},
}

\bib{MR3285819}{book}{
      author={Aster, Richard~C.},
      author={Borchers, Brian},
      author={Thurber, Clifford~H.},
       title={Parameter estimation and inverse problems},
     edition={Second},
   publisher={Elsevier/Academic Press, Amsterdam},
        date={2013},
        ISBN={978-0-12-385048-5},
         url={https://doi.org/10.1016/B978-0-12-385048-5.00001-X},
      review={\MR{3285819}},
}

\bib{10.1073/pnas.2311878121}{article}{
      author={Bahri, Y.},
      author={Dyer, E.},
      author={Kaplan, J.},
      author={Lee, J.},
      author={Sharma, U.},
       title={Explaining neural scaling laws},
        date={2024},
     journal={Proceedings of the National Academy of Sciences},
      volume={121},
}

\bib{barrata2023dolfinx}{article}{
      author={Barrata, Igor~A},
      author={Dean, Joseph~P},
      author={Dokken, J{\o}rgen~S},
      author={Habera, Michal},
      author={HALE, Jack},
      author={Richardson, Chris},
      author={Rognes, Marie~E},
      author={Scroggs, Matthew~W},
      author={Sime, Nathan},
      author={Wells, Garth~N},
       title={Dolfinx: The next generation fenics problem solving environment},
        date={2023},
}

\bib{MR2322235}{book}{
      author={Braess, Dietrich},
       title={Finite elements},
     edition={Third},
   publisher={Cambridge University Press, Cambridge},
        date={2007},
        ISBN={978-0-521-70518-9; 0-521-70518-5},
         url={https://doi-org.remotexs.ntu.edu.sg/10.1017/CBO9780511618635},
        note={Theory, fast solvers, and applications in elasticity theory, Translated from the German by Larry L. Schumaker},
      review={\MR{2322235}},
}

\bib{MR1665662}{article}{
      author={Brooks, Stephen~P.},
      author={Gelman, Andrew},
       title={General methods for monitoring convergence of iterative simulations},
        date={1998},
        ISSN={1061-8600,1537-2715},
     journal={J. Comput. Graph. Statist.},
      volume={7},
      number={4},
       pages={434\ndash 455},
         url={https://doi-org.remotexs.ntu.edu.sg/10.2307/1390675},
      review={\MR{1665662}},
}

\bib{brooks2011handbook}{book}{
      author={Brooks, Steve},
      author={Gelman, Andrew},
      author={Jones, Galin},
      author={Meng, Xiao-Li},
       title={Handbook of markov chain monte carlo},
   publisher={CRC press},
        date={2011},
}

\bib{MR2767184}{book}{
      author={\c~C\i~nlar, Erhan},
       title={Probability and stochastics},
      series={Graduate Texts in Mathematics},
   publisher={Springer, New York},
        date={2011},
      volume={261},
        ISBN={978-0-387-87858-4},
         url={https://doi-org.remotexs.ntu.edu.sg/10.1007/978-0-387-87859-1},
      review={\MR{2767184}},
}

\bib{CANTWELL2015205}{article}{
      author={Cantwell, C.D.},
      author={Moxey, D.},
      author={Comerford, A.},
      author={Bolis, A.},
      author={Rocco, G.},
      author={Mengaldo, G.},
      author={{De Grazia}, D.},
      author={Yakovlev, S.},
      author={Lombard, J.-E.},
      author={Ekelschot, D.},
      author={Jordi, B.},
      author={Xu, H.},
      author={Mohamied, Y.},
      author={Eskilsson, C.},
      author={Nelson, B.},
      author={Vos, P.},
      author={Biotto, C.},
      author={Kirby, R.M.},
      author={Sherwin, S.J.},
       title={Nektar++: An open-source spectral/hp element framework},
        date={2015},
        ISSN={0010-4655},
     journal={Computer Physics Communications},
      volume={192},
       pages={205\ndash 219},
         url={https://www.sciencedirect.com/science/article/pii/S0010465515000533},
}

\bib{CAO2023112104}{article}{
      author={Cao, Lianghao},
      author={O'Leary-Roseberry, Thomas},
      author={Jha, Prashant~K.},
      author={Oden, J.~Tinsley},
      author={Ghattas, Omar},
       title={Residual-based error correction for neural operator accelerated infinite-dimensional bayesian inverse problems},
        date={2023},
        ISSN={0021-9991},
     journal={Journal of Computational Physics},
      volume={486},
       pages={112104},
         url={https://www.sciencedirect.com/science/article/pii/S0021999123001997},
}

\bib{chen1995universal}{article}{
      author={Chen, Tianping},
      author={Chen, Hong},
       title={Universal approximation to nonlinear operators by neural networks with arbitrary activation functions and its application to dynamical systems},
        date={1995},
     journal={IEEE Transactions on Neural Networks},
      volume={6},
      number={4},
       pages={911\ndash 917},
}

\bib{Christen01122005}{article}{
      author={Christen, J.~Andrés},
      author={and, Colin~Fox},
       title={Markov chain monte carlo using an approximation},
        date={2005},
     journal={Journal of Computational and Graphical Statistics},
      volume={14},
      number={4},
       pages={795\ndash 810},
      eprint={https://doi.org/10.1198/106186005X76983},
         url={https://doi.org/10.1198/106186005X76983},
}

\bib{MR2558668}{article}{
      author={Cotter, S.~L.},
      author={Dashti, M.},
      author={Robinson, J.~C.},
      author={Stuart, A.~M.},
       title={Bayesian inverse problems for functions and applications to fluid mechanics},
        date={2009},
        ISSN={0266-5611},
     journal={Inverse Problems},
      volume={25},
      number={11},
       pages={115008, 43},
         url={https://doi.org/10.1088/0266-5611/25/11/115008},
      review={\MR{2558668}},
}

\bib{cotter2013mcmc}{article}{
      author={Cotter, Simon~L},
      author={Roberts, Gareth~O},
      author={Stuart, Andrew~M},
      author={White, David},
       title={Mcmc methods for functions: modifying old algorithms to make them faster},
        date={2013},
}

\bib{cui2011bayesian}{article}{
      author={Cui, Tiangang},
      author={Fox, Colin},
      author={O'sullivan, MJ},
       title={Bayesian calibration of a large-scale geothermal reservoir model by a new adaptive delayed acceptance metropolis hastings algorithm},
        date={2011},
     journal={Water Resources Research},
      volume={47},
      number={10},
}

\bib{doi:10.1137/130915005}{article}{
      author={Dodwell, T.~J.},
      author={Ketelsen, C.},
      author={Scheichl, R.},
      author={Teckentrup, A.~L.},
       title={A hierarchical multilevel markov chain monte carlo algorithm with applications to uncertainty quantification in subsurface flow},
        date={2015},
     journal={SIAM/ASA Journal on Uncertainty Quantification},
      volume={3},
      number={1},
       pages={1075\ndash 1108},
      eprint={https://doi-org.remotexs.ntu.edu.sg/10.1137/130915005},
         url={https://doi-org.remotexs.ntu.edu.sg/10.1137/130915005},
}

\bib{MR2421969}{article}{
      author={Dorn, Oliver},
      author={Villegas, Rossmary},
       title={History matching of petroleum reservoirs using a level set technique},
        date={2008},
        ISSN={0266-5611,1361-6420},
     journal={Inverse Problems},
      volume={24},
      number={3},
       pages={035015, 29},
         url={https://doi-org.remotexs.ntu.edu.sg/10.1088/0266-5611/24/3/035015},
      review={\MR{2421969}},
}

\bib{du2019gradient}{inproceedings}{
      author={Du, Simon},
      author={Lee, Jason},
      author={Li, Haochuan},
      author={Wang, Liwei},
      author={Zhai, Xiyu},
       title={Gradient descent finds global minima of deep neural networks},
organization={PMLR},
        date={2019},
   booktitle={International conference on machine learning},
       pages={1675\ndash 1685},
}

\bib{MR2231730}{article}{
      author={Efendiev, Y.},
      author={Hou, T.},
      author={Luo, W.},
       title={Preconditioning {M}arkov chain {M}onte {C}arlo simulations using coarse-scale models},
        date={2006},
        ISSN={1064-8275,1095-7197},
     journal={SIAM J. Sci. Comput.},
      volume={28},
      number={2},
       pages={776\ndash 803},
         url={https://doi-org.remotexs.ntu.edu.sg/10.1137/050628568},
      review={\MR{2231730}},
}

\bib{efendiev2005efficient}{article}{
      author={Efendiev, Yalchin},
      author={Datta-Gupta, Akhil},
      author={Ginting, Victor},
      author={Ma, Xiang},
      author={Mallick, Bani},
       title={An efficient two-stage markov chain monte carlo method for dynamic data integration},
        date={2005},
     journal={Water Resources Research},
      volume={41},
      number={12},
}

\bib{MR3372290}{article}{
      author={Efendiev, Yalchin},
      author={Jin, Bangti},
      author={Presho, Michael},
      author={Tan, Xiaosi},
       title={Multilevel {M}arkov chain {M}onte {C}arlo method for high-contrast single-phase flow problems},
        date={2015},
        ISSN={1815-2406,1991-7120},
     journal={Commun. Comput. Phys.},
      volume={17},
      number={1},
       pages={259\ndash 286},
         url={https://doi-org.remotexs.ntu.edu.sg/10.4208/cicp.021013.260614a},
      review={\MR{3372290}},
}

\bib{MR2050138}{book}{
      author={Ern, Alexandre},
      author={Guermond, Jean-Luc},
       title={Theory and practice of finite elements},
      series={Applied Mathematical Sciences},
   publisher={Springer-Verlag, New York},
        date={2004},
      volume={159},
        ISBN={0-387-20574-8},
         url={https://doi.org/10.1007/978-1-4757-4355-5},
      review={\MR{2050138}},
}

\bib{gelman1992inference}{article}{
      author={Gelman, Andrew},
      author={Rubin, Donald~B},
       title={Inference from iterative simulation using multiple sequences},
        date={1992},
     journal={Statistical science},
      volume={7},
      number={4},
       pages={457\ndash 472},
}

\bib{MR2436856}{article}{
      author={Giles, Michael~B.},
       title={Multilevel {M}onte {C}arlo path simulation},
        date={2008},
        ISSN={0030-364X},
     journal={Oper. Res.},
      volume={56},
      number={3},
       pages={607\ndash 617},
         url={https://doi.org/10.1287/opre.1070.0496},
      review={\MR{2436856}},
}

\bib{MR2429876}{article}{
      author={He, Yinnian},
       title={The {E}uler implicit/explicit scheme for the 2{D} time-dependent {N}avier-{S}tokes equations with smooth or non-smooth initial data},
        date={2008},
        ISSN={0025-5718,1088-6842},
     journal={Math. Comp.},
      volume={77},
      number={264},
       pages={2097\ndash 2124},
         url={https://doi-org.remotexs.ntu.edu.sg/10.1090/S0025-5718-08-02127-3},
      review={\MR{2429876}},
}

\bib{hestness2017deep}{article}{
      author={Hestness, Joel},
      author={Narang, Sharan},
      author={Ardalani, Newsha},
      author={Diamos, Gregory},
      author={Jun, Heewoo},
      author={Kianinejad, Hassan},
      author={Patwary, Md},
      author={Ali, Mostofa},
      author={Yang, Yang},
      author={Zhou, Yanqi},
       title={Deep learning scaling is predictable, empirically},
        date={2017},
     journal={arXiv preprint arXiv:1712.00409},
}

\bib{MR1043610}{article}{
      author={Heywood, John~G.},
      author={Rannacher, Rolf},
       title={Finite-element approximation of the nonstationary {N}avier-{S}tokes problem. {IV}. {E}rror analysis for second-order time discretization},
        date={1990},
        ISSN={0036-1429},
     journal={SIAM J. Numer. Anal.},
      volume={27},
      number={2},
       pages={353\ndash 384},
         url={https://doi.org/10.1137/0727022},
      review={\MR{1043610}},
}

\bib{MR4069815}{article}{
      author={Hoang, Viet~Ha},
      author={Quek, Jia~Hao},
      author={Schwab, Christoph},
       title={Analysis of a multilevel {M}arkov chain {M}onte {C}arlo finite element method for {B}ayesian inversion of log-normal diffusions},
        date={2020},
        ISSN={0266-5611},
     journal={Inverse Problems},
      volume={36},
      number={3},
       pages={035021, 46},
         url={https://doi.org/10.1088/1361-6420/ab2a1e},
      review={\MR{4069815}},
}

\bib{MR4246090}{article}{
      author={Hoang, Viet~Ha},
      author={Quek, Jia~Hao},
      author={Schwab, Christoph},
       title={Multilevel {M}arkov chain {M}onte {C}arlo for {B}ayesian inversion of parabolic partial differential equations under {G}aussian prior},
        date={2021},
     journal={SIAM/ASA J. Uncertain. Quantif.},
      volume={9},
      number={2},
       pages={384\ndash 419},
         url={https://doi.org/10.1137/20M1354714},
      review={\MR{4246090}},
}

\bib{hoang2016convergence}{inproceedings}{
      author={Hoang, Viet~Ha},
      author={Schwab, Christoph},
       title={Convergence rate analysis of mcmc-fem for bayesian inversion of log-normal diffusion problems},
organization={ETH Zurich},
        date={2016},
   booktitle={Research reports/seminar for applied mathematics},
      volume={2016},
}

\bib{MR3084684}{article}{
      author={Hoang, Viet~Ha},
      author={Schwab, Christoph},
      author={Stuart, Andrew~M.},
       title={Complexity analysis of accelerated {MCMC} methods for {B}ayesian inversion},
        date={2013},
        ISSN={0266-5611},
     journal={Inverse Problems},
      volume={29},
      number={8},
       pages={085010, 37},
         url={https://doi.org/10.1088/0266-5611/29/8/085010},
      review={\MR{3084684}},
}

\bib{iakovlevlearning}{inproceedings}{
      author={Iakovlev, Valerii},
      author={Heinonen, Markus},
      author={L{\"a}hdesm{\"a}ki, Harri},
       title={Learning continuous-time pdes from sparse data with graph neural networks},
        date={2020},
   booktitle={International conference on learning representations},
}

\bib{jameson1977finite}{inproceedings}{
      author={Jameson, Antony},
      author={Caughey, D},
       title={A finite volume method for transonic potential flow calculations},
        date={1977},
   booktitle={3rd computational fluid dynamics conference},
       pages={635},
}

\bib{10.1093/imanum/drz055}{article}{
      author={Jentzen, Arnulf},
      author={Kuckuck, Benno},
      author={Neufeld, Ariel},
      author={von Wurstemberger, Philippe},
       title={{Strong error analysis for stochastic gradient descent optimization algorithms}},
        date={202005},
        ISSN={0272-4979},
     journal={IMA Journal of Numerical Analysis},
      volume={41},
      number={1},
       pages={455\ndash 492},
      eprint={https://academic.oup.com/imajna/article-pdf/41/1/455/35970895/drz055.pdf},
         url={https://doi.org/10.1093/imanum/drz055},
}

\bib{jin2020quantifying}{article}{
      author={Jin, Pengzhan},
      author={Lu, Lu},
      author={Tang, Yifa},
      author={Karniadakis, George~Em},
       title={Quantifying the generalization error in deep learning in terms of data distribution and neural network smoothness},
        date={2020},
     journal={Neural Networks},
      volume={130},
       pages={85\ndash 99},
}

\bib{MR2102218}{book}{
      author={Kaipio, Jari},
      author={Somersalo, Erkki},
       title={Statistical and computational inverse problems},
      series={Applied Mathematical Sciences},
   publisher={Springer-Verlag, New York},
        date={2005},
      volume={160},
        ISBN={0-387-22073-9},
      review={\MR{2102218}},
}

\bib{karniadakis2005spectral}{book}{
      author={Karniadakis, George},
      author={Sherwin, Spencer~J},
       title={Spectral/hp element methods for computational fluid dynamics},
   publisher={Oxford University Press, USA},
        date={2005},
}

\bib{doi:10.1137/17M1135566}{article}{
      author={Krumscheid, S.},
      author={Nobile, F.},
       title={Multilevel monte carlo approximation of functions},
        date={2018},
     journal={SIAM/ASA Journal on Uncertainty Quantification},
      volume={6},
      number={3},
       pages={1256\ndash 1293},
      eprint={https://doi.org/10.1137/17M1135566},
         url={https://doi.org/10.1137/17M1135566},
}

\bib{laloy2013efficient}{article}{
      author={Laloy, Eric},
      author={Rogiers, Bart},
      author={Vrugt, Jasper~A},
      author={Mallants, Dirk},
      author={Jacques, Diederik},
       title={Efficient posterior exploration of a high-dimensional groundwater model from two-stage markov chain monte carlo simulation and polynomial chaos expansion},
        date={2013},
     journal={Water Resources Research},
      volume={49},
      number={5},
       pages={2664\ndash 2682},
}

\bib{LESHNO1993861}{article}{
      author={Leshno, Moshe},
      author={Lin, Vladimir~Ya.},
      author={Pinkus, Allan},
      author={Schocken, Shimon},
       title={Multilayer feedforward networks with a nonpolynomial activation function can approximate any function},
        date={1993},
        ISSN={0893-6080},
     journal={Neural Networks},
      volume={6},
      number={6},
       pages={861\ndash 867},
         url={https://www.sciencedirect.com/science/article/pii/S0893608005801315},
}

\bib{MR1925043}{book}{
      author={LeVeque, Randall~J.},
       title={Finite volume methods for hyperbolic problems},
      series={Cambridge Texts in Applied Mathematics},
   publisher={Cambridge University Press, Cambridge},
        date={2002},
        ISBN={0-521-81087-6; 0-521-00924-3},
         url={https://doi-org.remotexs.ntu.edu.sg/10.1017/CBO9780511791253},
      review={\MR{1925043}},
}

\bib{levin2017markov}{book}{
      author={Levin, David~A},
      author={Peres, Yuval},
       title={Markov chains and mixing times},
   publisher={American Mathematical Soc.},
        date={2017},
      volume={107},
}

\bib{DBLP:conf/iclr/LiKALBSA21}{inproceedings}{
      author={Li, Zongyi},
      author={Kovachki, Nikola~Borislavov},
      author={Azizzadenesheli, Kamyar},
      author={Liu, Burigede},
      author={Bhattacharya, Kaushik},
      author={Stuart, Andrew~M.},
      author={Anandkumar, Anima},
       title={Fourier neural operator for parametric partial differential equations},
        date={2021},
   booktitle={9th international conference on learning representations, {ICLR} 2021, virtual event, austria, may 3-7, 2021},
   publisher={OpenReview.net},
         url={https://openreview.net/forum?id=c8P9NQVtmnO},
}

\bib{lu2021learning}{article}{
      author={Lu, Lu},
      author={Jin, Pengzhan},
      author={Pang, Guofei},
      author={Zhang, Zhongqiang},
      author={Karniadakis, George~Em},
       title={Learning nonlinear operators via deeponet based on the universal approximation theorem of operators},
        date={2021},
     journal={Nature Machine Intelligence},
      volume={3},
      number={3},
       pages={218\ndash 229},
}

\bib{LU2022114778}{article}{
      author={Lu, Lu},
      author={Meng, Xuhui},
      author={Cai, Shengze},
      author={Mao, Zhiping},
      author={Goswami, Somdatta},
      author={Zhang, Zhongqiang},
      author={Karniadakis, George~Em},
       title={A comprehensive and fair comparison of two neural operators (with practical extensions) based on fair data},
        date={2022},
        ISSN={0045-7825},
     journal={Computer Methods in Applied Mechanics and Engineering},
      volume={393},
       pages={114778},
         url={https://www.sciencedirect.com/science/article/pii/S0045782522001207},
}

\bib{MR4537564}{article}{
      author={Lykkegaard, M.~B.},
      author={Dodwell, T.~J.},
      author={Fox, C.},
      author={Mingas, G.},
      author={Scheichl, R.},
       title={Multilevel delayed acceptance {MCMC}},
        date={2023},
        ISSN={2166-2525},
     journal={SIAM/ASA J. Uncertain. Quantif.},
      volume={11},
      number={1},
       pages={1\ndash 30},
         url={https://doi-org.remotexs.ntu.edu.sg/10.1137/22M1476770},
      review={\MR{4537564}},
}

\bib{doi:10.1137/110845598}{article}{
      author={Martin, James},
      author={Wilcox, Lucas~C.},
      author={Burstedde, Carsten},
      author={Ghattas, Omar},
       title={A stochastic newton mcmc method for large-scale statistical inverse problems with application to seismic inversion},
        date={2012},
     journal={SIAM Journal on Scientific Computing},
      volume={34},
      number={3},
       pages={A1460\ndash A1487},
      eprint={https://doi-org.remotexs.ntu.edu.sg/10.1137/110845598},
         url={https://doi-org.remotexs.ntu.edu.sg/10.1137/110845598},
}

\bib{marzouk2009stochastic}{article}{
      author={Marzouk, Youssef},
      author={Xiu, Dongbin},
       title={A stochastic collocation approach to bayesian inference in inverse problems},
        date={2009},
}

\bib{maulik2022efficient}{article}{
      author={Maulik, Romit},
      author={Rao, Vishwas},
      author={Wang, Jiali},
      author={Mengaldo, Gianmarco},
      author={Constantinescu, Emil},
      author={Lusch, Bethany},
      author={Balaprakash, Prasanna},
      author={Foster, Ian},
      author={Kotamarthi, Rao},
       title={Efficient high-dimensional variational data assimilation with machine-learned reduced-order models},
        date={2022},
     journal={Geoscientific Model Development},
      volume={15},
      number={8},
       pages={3433\ndash 3445},
}

\bib{mengaldo2021industry}{article}{
      author={Mengaldo, Gianmarco},
      author={Moxey, David},
      author={Turner, Michael},
      author={Moura, Rodrigo~Costa},
      author={Jassim, Ayad},
      author={Taylor, Mark},
      author={Peir{\'o}, Joaquim},
      author={Sherwin, Spencer~J},
       title={Industry-relevant implicit large-eddy simulation of a high-performance road car via spectral/hp element methods},
        date={2021},
     journal={SIAM Review},
      volume={63},
      number={4},
}

\bib{MR2897628}{article}{
      author={Mishra, S.},
      author={Schwab, Ch.},
      author={\v{S}ukys, J.},
       title={Multi-level {M}onte {C}arlo finite volume methods for nonlinear systems of conservation laws in multi-dimensions},
        date={2012},
        ISSN={0021-9991},
     journal={J. Comput. Phys.},
      volume={231},
      number={8},
       pages={3365\ndash 3388},
         url={https://doi.org/10.1016/j.jcp.2012.01.011},
      review={\MR{2897628}},
}

\bib{parolini2005mathematical}{article}{
      author={Parolini, Nicola},
      author={Quarteroni, Alfio},
       title={Mathematical models and numerical simulations for the america’s cup},
        date={2005},
     journal={Computer Methods in Applied Mechanics and Engineering},
      volume={194},
      number={9-11},
       pages={1001\ndash 1026},
}

\bib{petersen2018optimal}{article}{
      author={Petersen, Philipp},
      author={Voigtlaender, Felix},
       title={Optimal approximation of piecewise smooth functions using deep relu neural networks},
        date={2018},
     journal={Neural Networks},
      volume={108},
       pages={296\ndash 330},
}

\bib{peyret2002spectral}{book}{
      author={Peyret, Roger},
       title={Spectral methods for incompressible viscous flow},
   publisher={Springer},
        date={2002},
      volume={148},
}

\bib{MR3881695}{article}{
      author={Raissi, M.},
      author={Perdikaris, P.},
      author={Karniadakis, G.~E.},
       title={Physics-informed neural networks: a deep learning framework for solving forward and inverse problems involving nonlinear partial differential equations},
        date={2019},
        ISSN={0021-9991},
     journal={J. Comput. Phys.},
      volume={378},
       pages={686\ndash 707},
         url={https://doi.org/10.1016/j.jcp.2018.10.045},
      review={\MR{3881695}},
}

\bib{reddy2024accelerating}{article}{
      author={Reddy, Sohail},
      author={Fairbanks, Hillary},
       title={Accelerating multilevel markov chain monte carlo using machine learning models},
        date={2024},
     journal={Physica Scripta},
}

\bib{10.1007/978-3-319-24574-4_28}{inproceedings}{
      author={Ronneberger, Olaf},
      author={Fischer, Philipp},
      author={Brox, Thomas},
       title={U-net: Convolutional networks for biomedical image segmentation},
        date={2015},
   booktitle={Medical image computing and computer-assisted intervention -- miccai 2015},
      editor={Navab, Nassir},
      editor={Hornegger, Joachim},
      editor={Wells, William~M.},
      editor={Frangi, Alejandro~F.},
   publisher={Springer International Publishing},
     address={Cham},
       pages={234\ndash 241},
}

\bib{skamarock2019description}{techreport}{
      author={Skamarock, William~C},
      author={Klemp, Joseph~B},
      author={Dudhia, Jimy},
      author={Gill, David~O},
      author={Liu, Zhiquan},
      author={Berner, Judith},
      author={Wang, Wei},
      author={Powers, Jordan~G},
      author={Duda, Michael~G},
      author={Barker, Dale~M},
      author={Huang, Xiang-Yu},
       title={A description of the advanced research wrf version 4},
        type={NCAR Technical Note},
 institution={National Center for Atmospheric Research},
        date={2019},
      number={NCAR/TN-556+STR},
}

\bib{MR2652785}{article}{
      author={Stuart, A.~M.},
       title={Inverse problems: a {B}ayesian perspective},
        date={2010},
        ISSN={0962-4929},
     journal={Acta Numer.},
      volume={19},
       pages={451\ndash 559},
         url={https://doi.org/10.1017/S0962492910000061},
      review={\MR{2652785}},
}

\bib{MR2130010}{book}{
      author={Tarantola, Albert},
       title={Inverse problem theory and methods for model parameter estimation},
   publisher={Society for Industrial and Applied Mathematics (SIAM), Philadelphia, PA},
        date={2005},
        ISBN={0-89871-572-5},
         url={https://doi-org.remotexs.ntu.edu.sg/10.1137/1.9780898717921},
      review={\MR{2130010}},
}

\bib{tay2025optimization}{article}{
      author={Tay, Wee-Beng},
      author={Liao-Yang, Tian-Chun},
      author={Luo, Hong-Rui},
      author={Chen, Kai-Peng},
      author={Chen, Si-Run},
      author={Yang, Jun-Tao},
      author={See, Simon},
      author={Khoo, Boo-Cheong},
       title={Optimization of two-element airfoils using nvidia modulus, a physics-informed neural network solver},
        date={2025},
     journal={Journal of Aircraft},
       pages={1\ndash 9},
}

\bib{VillaPetraGhattas21}{article}{
      author={Villa, Umberto},
      author={Petra, Noemi},
      author={Ghattas, Omar},
       title={{HIPPYlib: An Extensible Software Framework for Large-Scale Inverse Problems Governed by PDEs: Part I: Deterministic Inversion and Linearized Bayesian Inference}},
        date={2021-04},
        ISSN={0098-3500},
     journal={ACM Trans. Math. Softw.},
      volume={47},
      number={2},
         url={https://doi.org/10.1145/3428447},
}

\bib{aadebffb88b448c89c654fcdda52c02b}{article}{
      author={Xu, Jinchao},
       title={Finite neuron method and convergence analysis},
        date={2020},
        ISSN={1991-7120},
     journal={Communications in Computational Physics},
      volume={28},
      number={5},
       pages={1707\ndash 1745},
         url={http://global-sci.org/intro/article_detail/cicp/18394.html},
}

\bib{MR4523340}{article}{
      author={Yang, Juntao},
      author={Hoang, Viet~Ha},
       title={Multilevel {M}arkov {C}hain {M}onte {C}arlo for {B}ayesian inverse problem for {N}avier-{S}tokes equation},
        date={2023},
        ISSN={1930-8337},
     journal={Inverse Probl. Imaging},
      volume={17},
      number={1},
       pages={106\ndash 135},
         url={https://doi.org/10.3934/ipi.2022033},
      review={\MR{4523340}},
}

\end{biblist}
\end{bibdiv}

\newpage
\appendix

\section{Hybrid two-level MCMC with Gaussian Prior}
\label{app:two-level-mcmc-gaussian}
In Section~\ref{sec:two-level-mcmc-uniform}, we introduced our new hybrid two-level MCMC for Bayesian inverse problems with uniform priors.
However, in several instances, it may be convenient to work with Gaussian priors. 
For completeness, in this section we discuss the hybrid two-level MCMC method for Gaussian prior. 
Similar to the case of uniform priors, we consider a forward model that predicts the states $u$ of a physical system given parameter $\mathbf{z}$. 
In this case, the prior is Gaussian.
Following the KL expansion~\eqref{eq:KL-expansion} setup, we assume $\mathbf{b} := (\|\psi_j\|_{L^\infty(D)})_{j=1,2,...,n} \in \ell^1$ and $\mathbf{\Bar{b}} := (\|\psi_j\|_{W^{1, \infty}(D)})_{j=1,2,...,n} \in \ell^1$. Then we define the measurable space $(U, \Gamma_b)$, with $\Gamma_b:= \{z \in \mathbb{R}^n, \quad \sum_{j=1}^n b_j |z_j| < \infty\} \in \mathcal{B}(\mathbb{R}^n)$, and $U$ is the parameter space.
We denote the standard Gaussian measure in $\mathbb{R}$ by $\gamma_1$.
Hence, the prior can be defined as $\gamma = \bigotimes_{j=1}^n \gamma_1$ on $(\mathbb{R}^n, \mathcal{B}(\mathbb{R}^n))$, and it completes the probability space $(U,\Gamma,\gamma)$, noting that $\Gamma_b$ has full Gaussian measure, i.e. $\gamma(\Gamma_b) = 1$.

\begin{assumption}
\label{assumption:numerical-error-gaussian}
Let $u$ be the solution of the the forward problem in equation~\eqref{eq:bayesian-setup}. We assume that $u \in V$, where $V$ is a suitable vector spece, e.g., a Sobolev space. The FEM approximation gives 
\begin{equation}
\label{eq:forward_estimator_gaussian}
\|u(z) - u^{\ell}(z)\|_V \leq C \exp(c\sum_{j=1}^n b_j |z_j|)(1+\sum_{j=1}^n \Bar{b}_j |z_j|) 2^{-l}.
\end{equation}
\end{assumption}
The numerical error estimate in Assumption \ref{assumption:numerical-error-gaussian} might not be true for all problems with Gaussian prior. 
It is problem-dependent. 
Yet, it is a typical error rate found in problems such as elliptic equations, diffusion problems, and parabolic problems with unknown coefficients, as investigated in~\cites{hoang2016convergence, MR4246090}. 
The right-hand side  of Equation~\eqref{eq:forward_estimator_gaussian} differs from Assumption~\ref{assumption:numerical-error} made for a uniform prior, having an additional exponential term that depends on $\mathbf{z}$.

This specific form of approximation error is of significant interest for our two-level hybrid MCMC approach, because the exponential term in the error will lead to divergence of the method (in contrast to the uniform prior introduced in Section~\ref{sec:two-level-mcmc-uniform}). However, with Fernique's theorem and an additional indication function to be presented below, we still can reach a similar posterior estimation error rate as the case in uniform prior.
Using Assumption~\ref{assumption:numerical-error-gaussian}, we can write the error between the DL-based surrogate and the numerical discretization as follows
\begin{align}
&\|u^{\rm num}(z) - u^{\rm DL}(z)\|_V \nonumber\\
&\;\;\;\;\leq C(1+2^{\epsilon}) \exp\bigg(c\sum_{j=1}^n b_j |z_j|)(1+\sum_{j=1}^n \Bar{b}_j |z_j|\bigg) 2^{-l}, \label{eq:error-gaussian-1}
\end{align}
\begin{align}
&|\Phi^{\rm num}(z;  y) - \Phi^{\rm DL}(z;  y)| \nonumber\\
&\;\;\;\;\leq C(1+2^{\epsilon}) \exp\bigg(c\sum_{j=1}^n b_j |z_j|)(1+\sum_{j=1}^n \Bar{b}_j |z_j|\bigg) 2^{-l}. \label{eq:error-gaussian-2}
\end{align}
Next we derive the two level MCMC approach for Gaussian prior. In order to avoid the unboundedness from the exponential term, we make use of the following switching function
\begin{equation}\label{eq:switching}
    S(z) = 
    \begin{cases}
1, \text{ if } \Phi^{\rm num}(z, y) - \Phi^{\rm DL}(z, y) \leq 0,\\
0, \text{ otherwise}.
\end{cases}
\end{equation}  
With the switching function~\eqref{eq:switching}, we can write the expected QoI as follows
\begin{align}
\nonumber
&\left(\mathbb{E}^{\gamma^{\rm num}}-\mathbb{E}^{\gamma^{\rm DL}}\right)[Q]  \\\nonumber
&= \frac{1}{N^{\rm num}}\int_{\Gamma_{\mathbf{b}}} \exp(-\Phi^{\rm num}) Q S(z) d \gamma - \frac{1}{N^{\rm DL}}\int_{\Gamma_{\mathbf{b}}} \exp(-\Phi^{\rm DL}) Q S(z) d \gamma  \\\nonumber
& +\frac{1}{N^{\rm num}}\int_{\Gamma_{\mathbf{b}}} \exp(-\Phi^{\rm num}) Q (1-S(z)) d \gamma - \frac{1}{N^{\rm DL}}\int_{\Gamma_{\mathbf{b}}} \exp(-\Phi^{\rm DL}) Q (1-S(z)) d \gamma \\ \nonumber
&= \frac{1}{N^{\rm num}}\int_{\Gamma_{\mathbf{b}}} \exp(-\Phi^{\rm num})(1-\exp(\Phi^{\rm num}-\Phi^{\rm DL})) Q S(z) d\gamma  \\ \nonumber
&+ \Big(\frac{1}{N^{\rm num}}-\frac{1}{N^{\rm DL}}\Big)\int_{\Gamma_{\mathbf{b}}} \exp(-\Phi^{\rm DL}) Q S(z) d \gamma \\ \nonumber
&+\frac{1}{N^{\rm DL}}\int_{\Gamma_{\mathbf{b}}} \exp(-\Phi^{\rm DL})(\exp(\Phi^{\rm DL}-\Phi^{\rm num})-1) Q (1-S(z))d\gamma \\ \nonumber
&+ \Big(\frac{1}{N^{\rm num}} - \frac{1}{N^{\rm DL}}\Big) \int_{\Gamma_{\mathbf{b}}} \exp(-\Phi^{\rm num}) Q (1-S(z))d \gamma. \label{eq:two-level-gaussian} \\
\end{align}  
The constant $(1/N^{\rm num}-1/N^{\rm DL})$, can be estimated via
\begin{align}
\nonumber
&\Big(\frac{1}{N^{\rm num}}-\frac{1}{N^{\rm DL}}\Big) = \\\nonumber
&=\frac{1}{N^{\rm num} N^{\rm DL}} \int_{\Gamma_{\mathbf{b}}}\left(\exp \left(-\Phi^{\rm DL}(z, y)\right)-\exp \left(-\Phi^{\rm num}(z, y)\right)\right)\left(S(z)+1-S(z)\right) d \gamma(z) \\\nonumber
&=\frac{1}{N^{\rm num} N^{\rm DL}} \int_{\Gamma_{\mathbf{b}}} \exp \left(-\Phi^{\rm num}(z, y)\right)\left(\exp \left(\Phi^{\rm num}(z, y)-\Phi^{\rm DL}(z, y)\right)-1\right) S(z) d \gamma(z)\\\nonumber
&+\frac{1}{N^{\rm num} N^{\rm DL}} \int_{\Gamma_{\mathbf{b}}} \exp \left(-\Phi^{\rm DL}(z, y)\right)\left(1-\exp \left(\Phi^{\rm DL}(z, y)-\Phi^{\rm num}(z, y)\right)\right)\left(1-S(z)\right) d \gamma(z) \\\nonumber
&=\frac{1}{N^{\rm DL}} \mathbb{E}^{\gamma^{\rm num}}\left[\left(\exp \left(\Phi^{\rm num}(z, y)-\Phi^{\rm DL}(z, y)\right)-1\right) S(z)\right]\\\nonumber
&+\frac{1}{N^{\rm num}} \mathbb{E}^{\gamma^{\rm DL}}\left[\left(1-\exp \left(\Phi^{\rm DL}(z, y)-\Phi^{\rm num}(z, y)\right)\right)\left(1-S(z)\right)\right]. \label{eq:coefficient gaussian}\\
\end{align}  
Combining equations~\eqref{eq:two-level-expansion}, \eqref{eq:two-level-gaussian} and \eqref{eq:coefficient gaussian}, we can derive the overall estimator of our hybrid two-level MCMC approach
\begin{align*}
    \mathbf{E}^{hybrid}(Q) &= \mathbf{E}^{\gamma^{\rm num}} [A_1] + \mathbf{E}^{\gamma^{\rm num}} [A_3] \cdot \mathbf{E}^{\gamma^{\rm DL}}[A_4+A_8] \\
    &+ \mathbf{E}^{\gamma^{\rm DL}}[A_2] + \mathbf{E}^{\gamma^{\rm DL}}[A_5] \cdot \mathbf{E}^{\gamma^{\rm num}} [A_6+A_7] + \mathbf{E}^{\gamma^{\rm DL}}[Q],
\end{align*}
where the terms $A_1, A_2, A_3, A_4, A_5, A_6, A_7$ and $A_8$ are defined as follows
\begin{align*}
A_1 &= (1-\exp(\Phi^{\rm num}(z)-\Phi^{\rm DL}(z))) Q(z) S(z), \\
A_2 &= (\exp(\Phi^{\rm DL}(z)-\Phi^{\rm num}(z))-1) Q(z)(1-S(z)), \\
A_3 &= (\exp(\Phi^{\rm num}(z)-\Phi^{\rm DL}(z))-1) S(z), \\
A_4 &= Q(z) \cdot S(z), \\
A_5 &= (1-\exp(\Phi^{\rm DL}(z)-\Phi^{\rm num}(z)))(1-S(z)) \\
A_6 &= \exp(\Phi^{\rm num}-\Phi^{\rm DL})Q(z) S(z), \\
A_7 &= Q(z) (1-S(z)), \\
A_8 &= \exp(\Phi^{\rm DL}(z)-\Phi^{\rm num}(z))Q(z) (1-S(z)).
\end{align*}
We now perform the error analysis of our method, under the assumption of Gaussian priors. 
In analogy with what we have done for uniform priors in Section~\ref{sec:two-level-mcmc-uniform}, we decompose the error in three terms: 
\begin{subequations}
\begin{align}
&\mathbb{E}^{\gamma^y}[Q] - \mathbf{E}^{\rm hybrid}[Q] = \mathrm{I} + \mathrm{II} + \mathrm{III}, \\
&\mathrm{I} := \mathbb{E}^{\gamma^y}[Q] -  \mathbb{E}^{\gamma^{\rm num}}[Q], \label{eq:gauss-error-1} \\[0.5em]
&\mathrm{II} := \mathbb{E}^{\gamma^{DL}}[Q] - \mathbf{E}^{\gamma^{DL}}_{M_{\rm num}}[Q], \label{eq:gauss-error-2} \\[0.5em]
&\mathrm{III} := \mathbb{E}^{\gamma^{\rm num}} [A_1] - \mathbf{E}^{\gamma^{\rm num}}_{M_{\rm num}} [A_1] + \mathbb{E}^{\gamma^{\rm num}}[A_3] \cdot \mathbb{E}^{\gamma^{\rm DL}}[A_4+A_8] \nonumber \\[0.1em]
& \;\;\;\; - \mathbf{E}^{\gamma^{\rm num}}_{M_{\rm num}} [A_3] \cdot \mathbf{E}^{\gamma^{\rm DL}}_{M_{\rm num}}[A_4+A_8] + \mathbb{E}^{\gamma^{\rm DL}}[A_2] -\mathbf{E}^{\gamma^{\rm DL}}_{M_{\rm num}}[A_2] \nonumber \\[0.1em]
& \;\;\;\; + \mathbb{E}^{\gamma^{\rm DL}}[A_5] \cdot \mathbb{E}^{\gamma^{\rm num}} [A_6+A_7] - \mathbf{E}^{\gamma^{\rm DL}}_{M_{\rm num}}[A_5] \cdot \mathbf{E}^{\gamma^{\rm num}}_{M_{\rm num}} [A_6+A_7]. \label{eq:gauss-error-3}
\end{align}
\end{subequations}
Similarly to Section~\ref{sec:two-level-mcmc-uniform}, for each error term, we have the following error bounds
\begin{subequations}
\begin{align}
&|\mathrm{I}| < C 2^{-L}, \\
&|\mathrm{II}| < M_{\rm DL}^{-1/2}, \\
&\mathcal{E}[|\mathrm{III}|^2] < C (1+2^{\epsilon})^2 M_{\rm num}^{-1} 2^{-2L}.
\end{align}
\end{subequations}
Therefore, choosing $M_{\rm DL} = C_{\rm DL} 2^{2L}$ and $M_{\rm num} = C_{\rm num}(1+2^{\epsilon})^2$, allows us to obtain a theorem for the overall error estimate.
\begin{thm}
\label{thm:gaussian_estimator}
With $M_{\rm DL} = C_{\rm DL} 2^{2L}$ and $M_{\rm num} = C_{\rm num}(1+2^{\epsilon})^2$, we have the following theoretical error estimate  of our hybrid two-level MCMC approach under Gaussian priors
\begin{equation}
\mathcal{E}_{\rm hybrid}[|\mathbb{E}^{\gamma^y}[Q] - \mathbf{E}^{\rm hybrid}[Q]|] \leq C_{\rm hybrid} 2^{-L}.
\end{equation}
\end{thm}
From Theorem~\ref{thm:gaussian_estimator}, we see the same results as the one from Theorem~\ref{thm:final_error_estimate} derived from uniform prior setup.
The theorem and the preceding assumptions are typically valid for log-normal priors with elliptic, diffusion, and parabolic equations. The proof for the two-dimensional Navier-Stokes equation is not available to the best of our knowledge. However, some experimental results also show the theorem for multilevel MCMC with Gaussian prior works for the two-dimensional Navier-Stokes equation~\cite{MR4523340}.  

\newpage
\section{Forward numerical solvers and deep learning surrogate models}
We present the details of the forward numerical solvers used for solving each of the problems presented in Section~\ref{sec:experiments}, along with the corresponding DL surrogate.

\subsection{Poisson equation}
\label{section:poisson_solver}
The Poisson equation~\ref{eq:poisson} in Section~\ref{subsec:poisson} is solved with the finite element package Fenicsx~\cite{barrata2023dolfinx}, whereby first-order Lagrange finite elements are used to discretize the equation.
The mesh adopted is constituted of $32 \times 32$ elements as shown in Figure~\ref{fig:mesh1}.
A direct LU solver from MUMPS backend is used to solve the assembled linear system on CPU.

For the simple uniform prior setup in Section~\ref{sec:poisson-uniform}, we choose a fully connected ReLU neural network to learn the forward mapping from the 4000 generated training samples. 
The input field of $33 \times 33$ is flattened and fed in as input data. 
Two hidden layers are included, each with 512 nodes. 
The neural network is implemented with PyTorch. 
The trained neural network is used as the DL-based surrogate model in the experiment.

For the Gaussian prior setup in Section~\ref{sec:poisson-gauss}, we choose a convolutional neural network (CNN) to be our DL-based surrogate model. 
The CNN model consists of 3 encoding layers, 1 fully connected layer, and 3 decoding layers. 
There are 8 kernels in each convolutional layer, and the kernel size is $(3, 3)$ with stride size $(2, 2)$.

\subsection{Nonlinear reaction-diffusion equation}
\label{section:reaction-diffusion_solver}
To solve the reaction diffusion equation~\eqref{eqn:reaction_diffusion} in Section~\ref{sec:reaction-diffusion}, we also use the finite element package Fenicsx~\cite{barrata2023dolfinx}, whereby first-order Lagrange finite elements are used to discretize the equation.
The resulting nonlinear system is solved with the Newton solver from the PETSC backend. 
Each linearized Newton iteration step is solved with LU direct solver with MUMPS backend. 

For the DL-based surrogate model, we choose message passing graph neural network (MPGNN) for the uniform prior case. 
We refer to the Graph-PDE architecture~\cite{iakovlevlearning} as our reference. 
Instead of training the MPGNN for a time dependent problem, here we train the neural network for a time independent problem. 
Two fully connected multilayer neural networks are used for the message passing and state update. 
There are two hidden layers each with 64 nodes in both the message passing and state update neural networks. 
The hidden state vector output from the message passing neural network is of size 64. 
A Dirichlet boundary condition is also imposed on the MPGNN. 
Five layers of MPGNN are stacked in the model used in the experiments. 
The MPGNN is trained with 2000 training samples with Adam optimizer and trained for 10000 epoches. 

We then choose U-net proposed in~\cite{10.1007/978-3-319-24574-4_28} for the Gaussian prior experiment. 
It is well-known for its outstanding performance in multi-scale physical problems. 
The U-net consists of 3 layers, each with dimensions of $32 \times 32$, $ 16 \times 16$, and $8 \times 8$. 
For each convolutional layer we have 8 kernel with size $(3, 3)$ and stride size $(1, 1)$.

\subsection{Navier-Stokes equations in the vorticity form}
\label{section:ns_solver}
To solve the Navier-Stokes equations~\ref{eq:ns} in Section~\ref{sec:navier-stokes}, we coded a simple numerical pseudo-spectral solver with PyTorch, which is accelerated on GPU. 
A total of $64 \times 64$ collocation points are used for the experiments. 
The interested reader can refer to~\cite{peyret2002spectral} for details of the spectral method implemented. 

For the DL-based surrogate model, we first choose the DeepONet~\cite{lu2021learning} as the deep learning model for the uniform prior experiment. 
Two fully connected multilayer neural networks are used as the branch net and trunk net. 
Both the branch net and trunk net have two hidden layers with 64 nodes in the neural network. 
The DeepONet model is trained with 4000 samples of numerical data for 10000 epoches. 

Then we choose the Fourier Neural Operator for the Gaussian prior experiment, as it demonstrated its ability to learn the dynamics of the Navier-Stokes equation in~\cite{DBLP:conf/iclr/LiKALBSA21} and has also been shown to be used in a Bayesian inversion problem setup with MCMC. 
Specifically, the two-dimensional Fourier neural operator with tensor layers is used. 
In this experiment we used 12 modes for height and width, 8 hidden channels, and 4 layers in the FNO. 
The neural network is trained with 2000 training data for 10000 epochs. 

\newpage
\section{Diagnostics details of numerical experiment ~\ref{sec:poisson-uniform}}
\label{sec:diagnostics}
\begin{figure}[bhtp!]
    \centering
    \includegraphics[width=\linewidth]{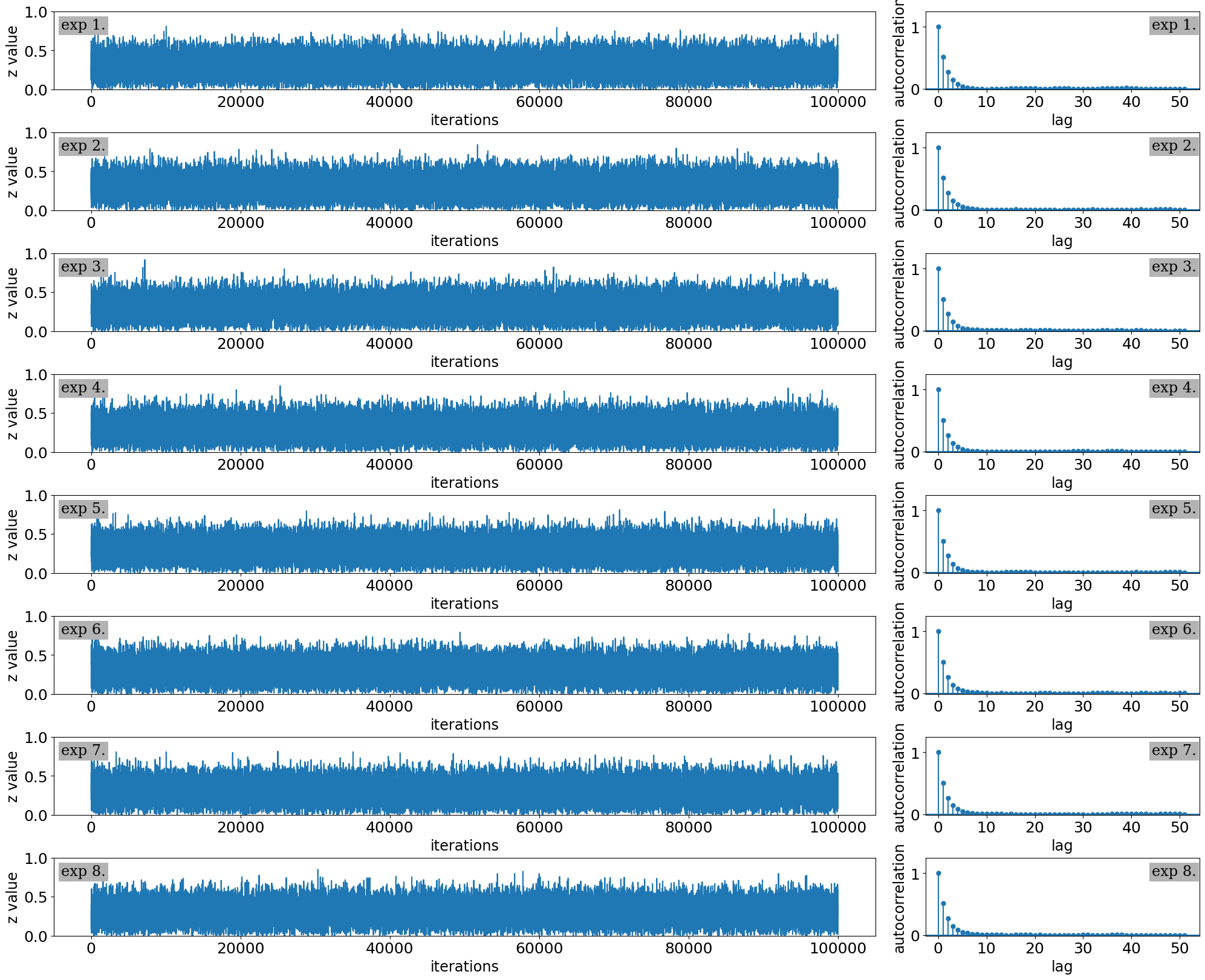}
    \caption{Left) Trace plots of MCMC chains with DL-based surrogate models from eight experiments with different initial sample showing good mixing of samples. Right) Autocorrelation Function(ACF) plot from the same eight MCMC chains showing rapid decay, which indicates good mixing and a high number of Effective Sample Size (ESS).}
    \label{fig:elliptic_coef_uniform_MLchain_traceplot}
\end{figure}
\begin{figure}
    \centering
    \includegraphics[width=\linewidth]{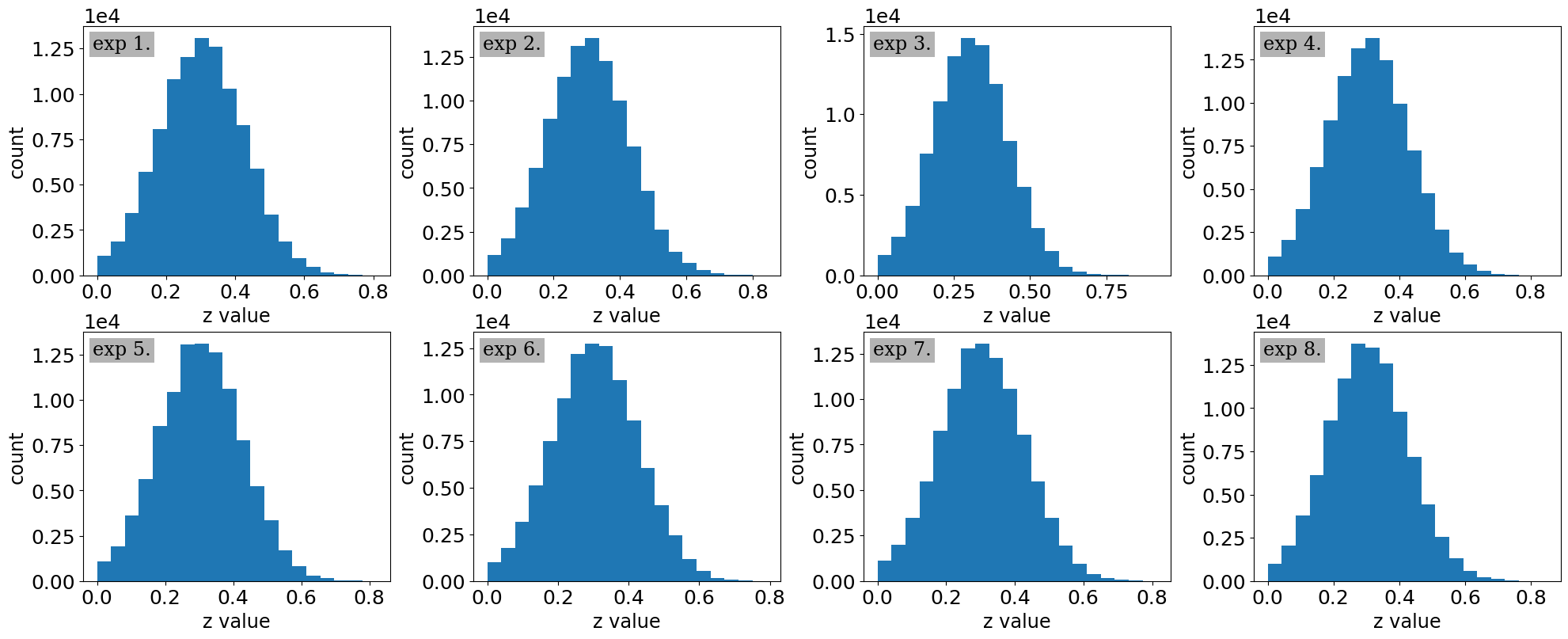}
    \caption{Histogram of MCMC chains with DL-based surrogate models from eight experiments with different initial sample showing similar posterior distribution of the samples with good mixing}
    \label{fig:elliptic_coef_uniform_MLchain_histogram}
\end{figure}
\begin{figure}
    \centering
    \subfloat[Sample mean of eight MCMC chains]{
        \includegraphics[width=0.45\textwidth]{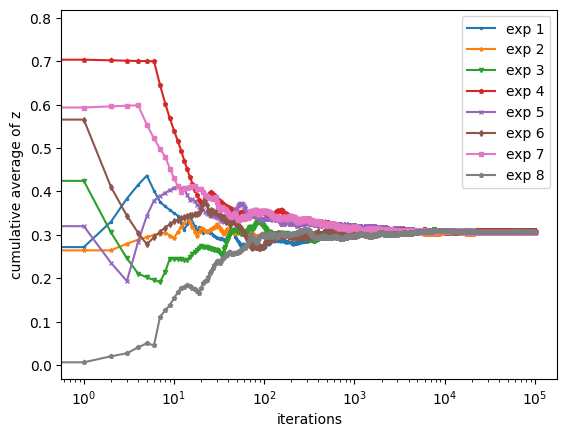}
        \label{fig:sample mean of elliptic uniform}
    }
    \hfill
    \subfloat[Potential Scale Reduction Factor]{
        \includegraphics[width=0.48\textwidth]{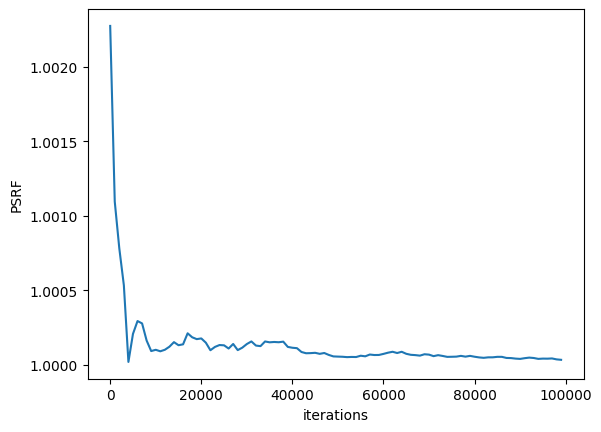}
        \label{fig:PSRF elliptic uniform}
    }
    \caption{Left) Cumulative average of samples of eight independent DL-based MCMC chain showing good inter-chain convergence. Right) PSRF of eight independent DL-based surrogate MCMC chains showing rapid decay and stable small PSRF values indicating good inter-chain convergence}
    \label{fig:sample mean and psrf uniform elliptic}
\end{figure}
\begin{figure}
    \centering
    \includegraphics[width=\linewidth]{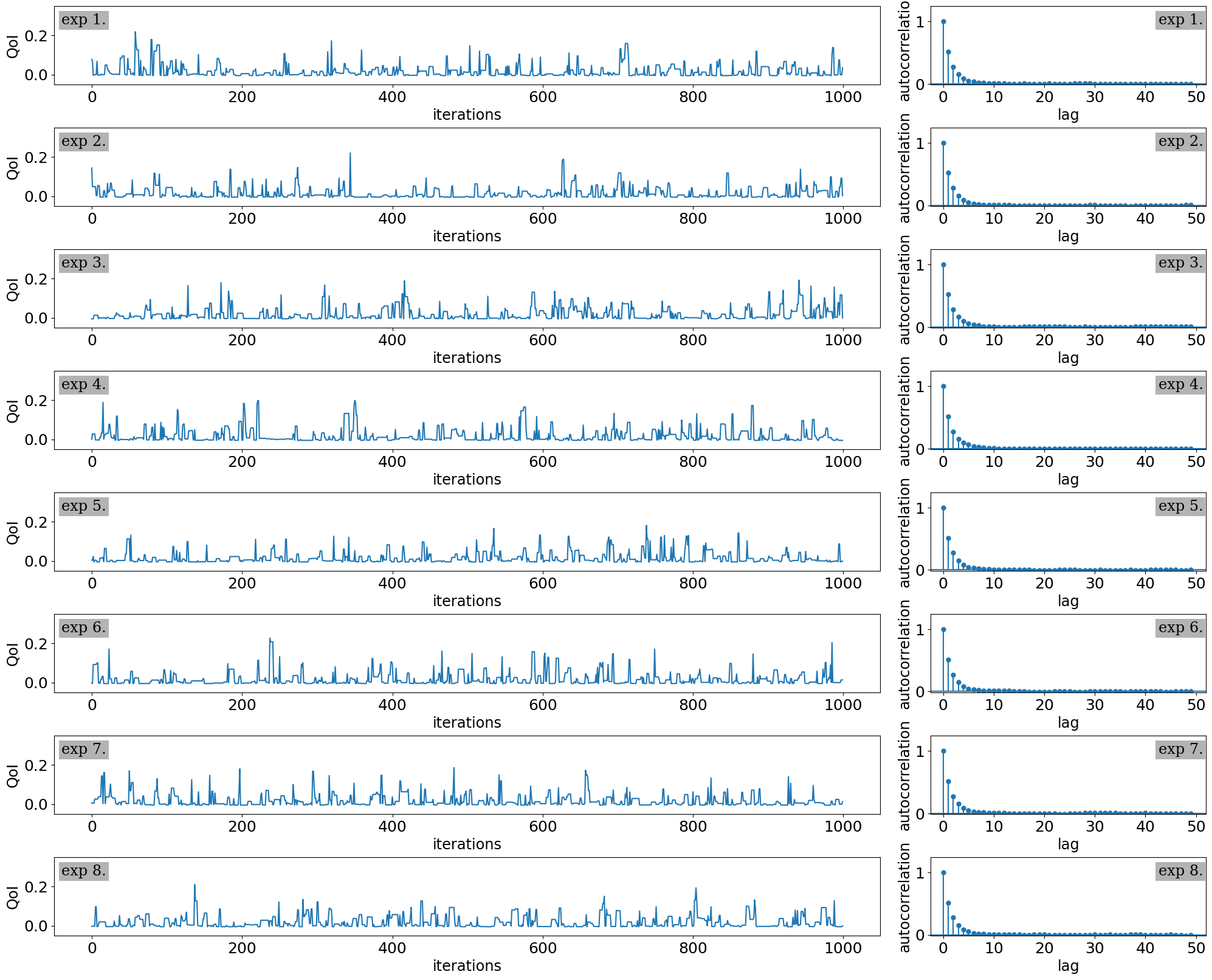}
    \caption{Left) Trace plot of QoI $(1-\exp(\Phi^{\rm num}-\Phi^{\rm DL}))z$ from 1000 numerical samples (1\% of total DL-based surrogate samples) from eight independent experiments with different initial sample showing good mixing of samples. Right) Autocorrelation Function(ACF) plot from the same eight MCMC chains showing rapid decay indicating good mixing and high number of Effective Sample Size(ESS)}
    \label{fig:elliptic_uniform_qoia_traceplot_acf}
\end{figure}
\begin{figure}
    \centering
    \includegraphics[width=\linewidth]{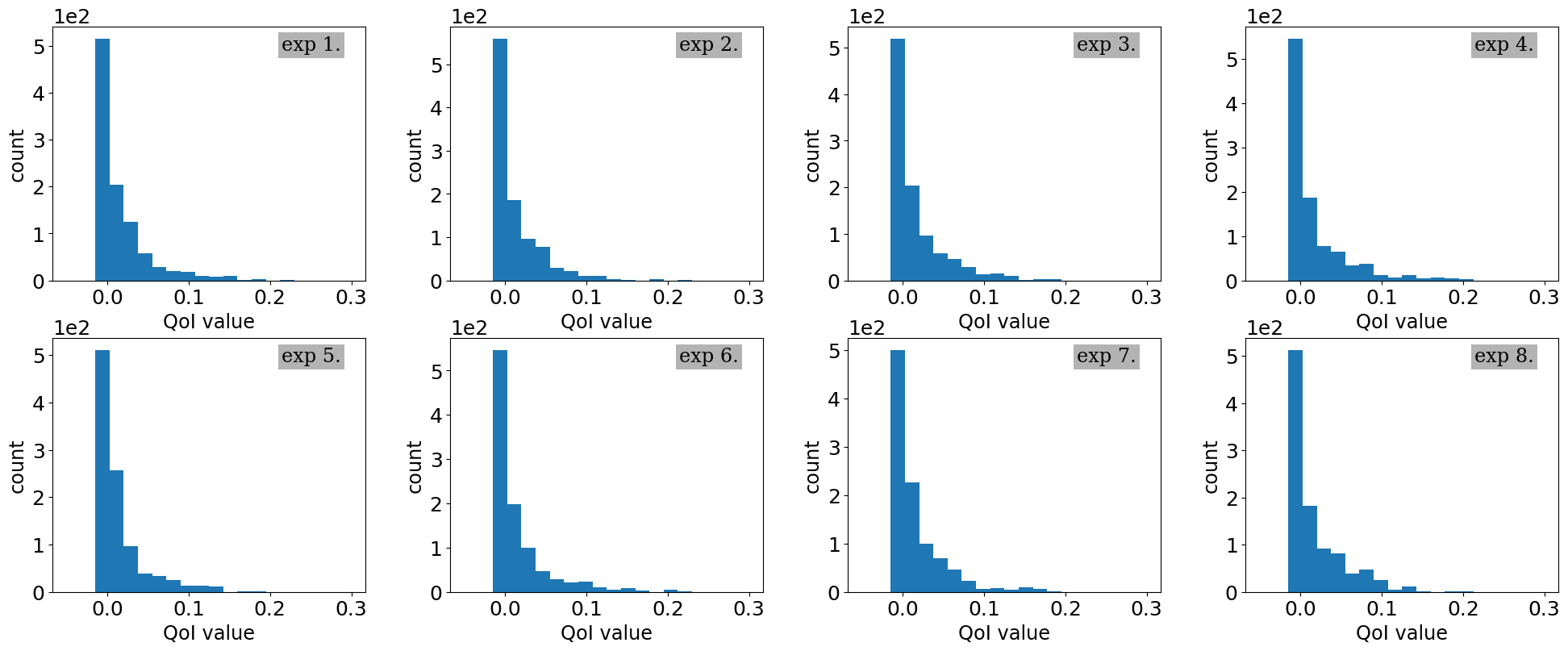}
    \caption{Histogram of of QoI $(1-\exp(\Phi^{\rm num}-\Phi^{\rm DL}))z$ from 1000 numerical samples (1\% of total DL-based surrogate samples) from eight independent MCMC chain with different initial sample showing similar posterior distribution of QoI with good mixing}
    \label{fig:qoia elliptic uniform hist}
\end{figure}
\begin{figure}
    \centering
    \includegraphics[width=\linewidth]{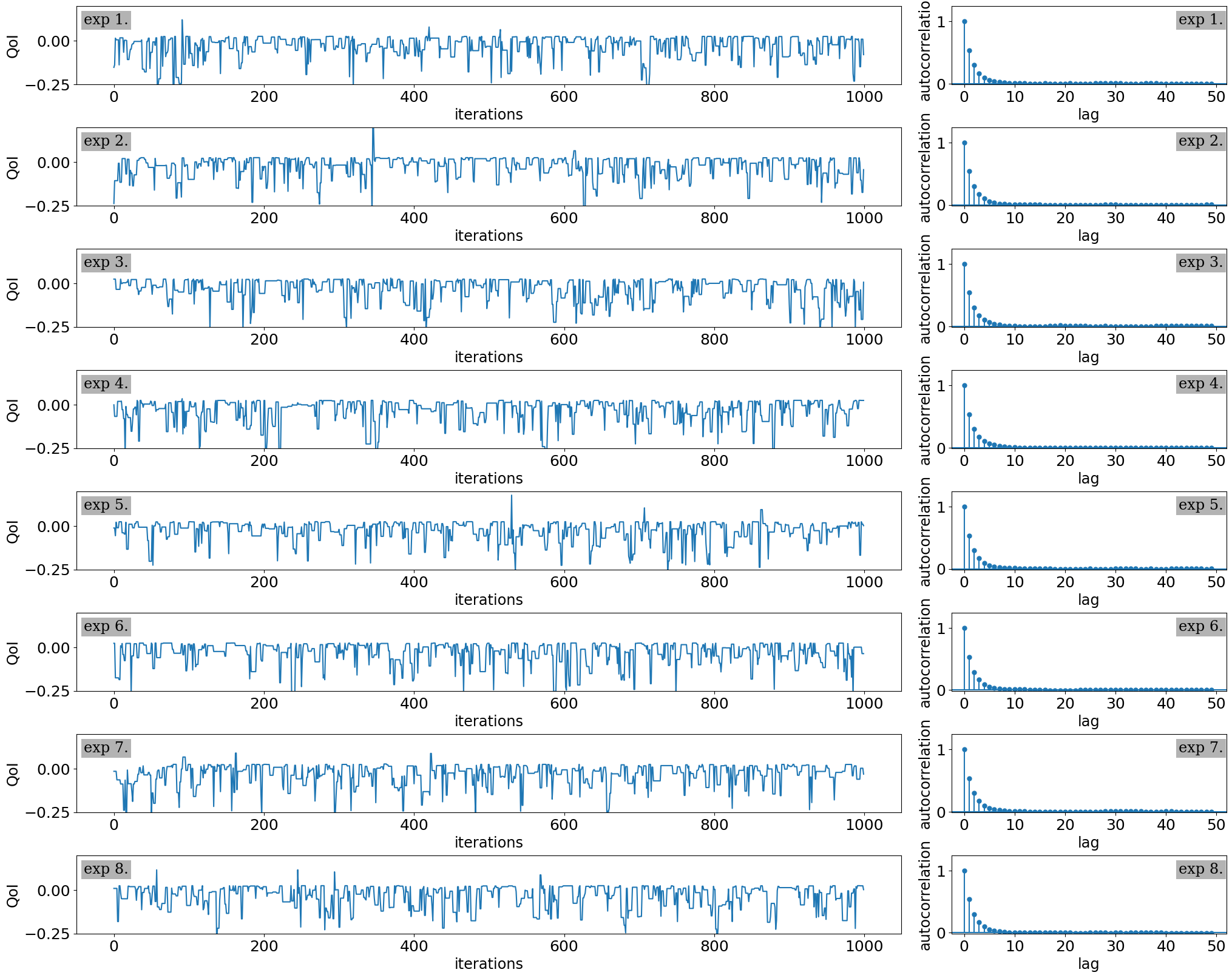}
    \caption{Left) Trace plot of QoI $\exp(\Phi^{\rm num}-\Phi^{\rm DL})-1$ from 1000 numerical samples (1\% of total DL-based surrogate samples) from eight independent MCMC chain with different initial sample showing good mixing of samples. Right) Autocorrelation Function (ACF) plot from the same eight MCMC chains showing rapid decay indicating good mixing and high number of Effective Sample Size(ESS)}
    \label{fig:elliptic_uniform_qoib_traceplot_acf}
\end{figure}
\begin{figure}
    \centering
    \includegraphics[width=\linewidth]{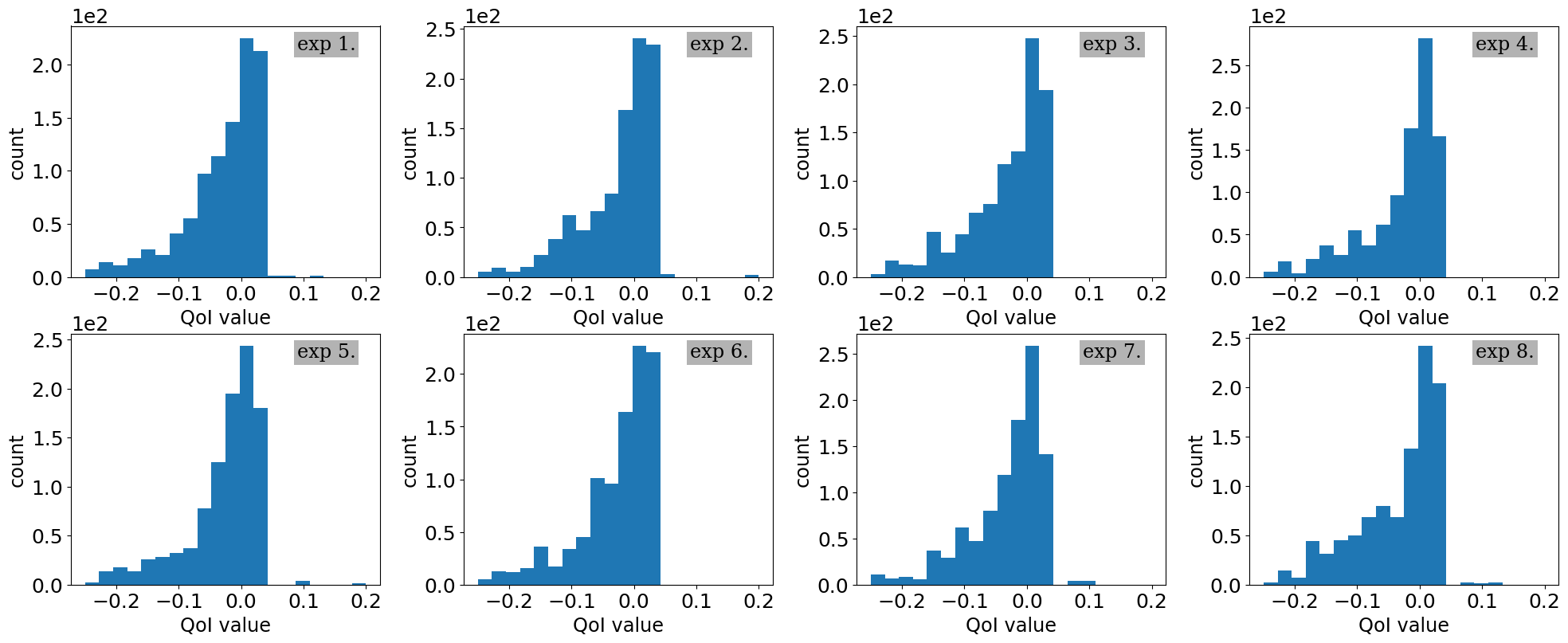}
    \caption{Histogram of QoI $(1-\exp(\Phi^{\rm num}-\Phi^{\rm DL}))z$ from 1000 numerical samples (1\% of total DL-based surrogate samples) from eight independent experiments with different initial sample showing similar posterior distribution of the QoI with good mixing.}
    \label{fig:qoib elliptic uniform hist}
\end{figure}
\begin{figure}
    \centering
    \subfloat[Sample mean of eight MCMC chains]{
        \includegraphics[width=0.46\textwidth]{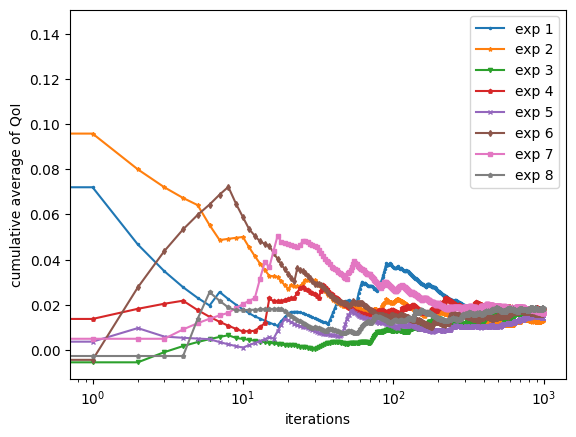}
    }
    \hfill
    \subfloat[Potential Scale Reduction Factor]{
        \includegraphics[width=0.46\textwidth]{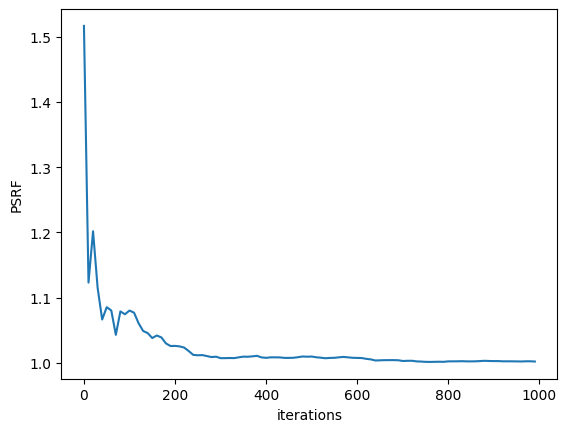}
    }
    \caption{Left) Cumulative average of QoI $(1-\exp(\Phi^{\rm num}-\Phi^{\rm DL}))\cdot z$ from eight independent numerical MCMC chains showing good inter-chain convergence. Right) PSRF of the same eight MCMC chains showing rapid decay and stable small PSRF values indicating a good inter-chain convergence. }
    \label{fig:qoia sample mean and psrf uniform elliptic}
\end{figure}
\begin{figure}
    \centering
    \subfloat[Sample mean of eight MCMC chains]{
        \includegraphics[width=0.47\textwidth]{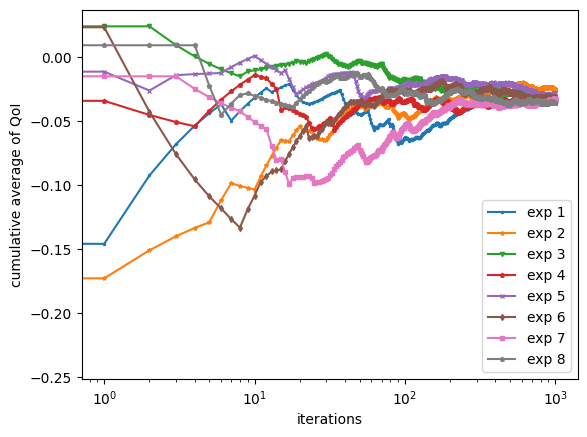}
    }
    \hfill
    \subfloat[Potential Scale Reduction Factor]{
        \includegraphics[width=0.46\textwidth]{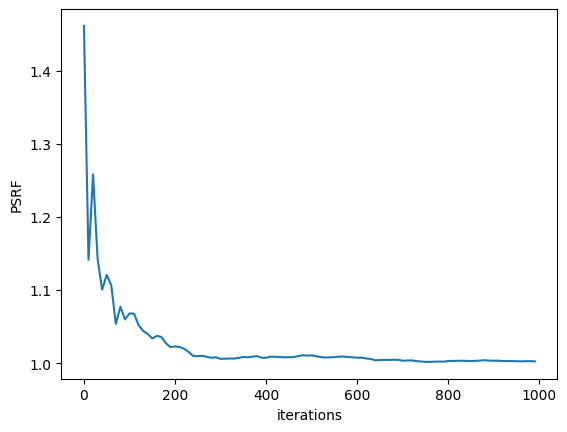}
    }
    \caption{Left) Cumulative average of QoI $\exp(\Phi^{\rm num}-\Phi^{\rm DL})-1$ from eight independent numerical MCMC chains with different initial sample showing good inter-chain convergence. Right) PSRF of the same eight MCMC chains showing rapid decay and stable small PSRF values indicating a good inter-chain convergence}
    \label{fig:qoib sample mean and psrf uniform elliptic}
\end{figure}
\end{document}